\documentclass[12pt,reqno]{amsart}

\usepackage{amsmath}
\usepackage{enumerate}
\usepackage{amsfonts}
\usepackage{amssymb}
\usepackage{amsthm}
\usepackage{bbm}
\usepackage{cite}
\usepackage{xcolor}
\usepackage{tabularx}
\usepackage{graphicx}
\usepackage{xr}
\usepackage{xspace}
\usepackage{thmtools}
\usepackage{thm-restate}
\usepackage[multiple]{footmisc}
\usepackage{microtype}
\usepackage[font=small,format=plain,labelfont=sf,up,textfont=sf,up]{caption}
\usepackage[tmargin=1.2in,bmargin=1.2in,lmargin=1.2in,rmargin=1.2in]{geometry}
\usepackage{hyperref}
\usepackage[textsize=small]{todonotes}
\usepackage{soul}

\definecolor{Red}{rgb}{1,0,0}
\definecolor{halfgray}{gray}{0.55}
\definecolor{webgreen}{rgb}{0,.5,0}
\definecolor{webbrown}{rgb}{.6,0,0}
\definecolor{Maroon}{cmyk}{0, 0.87, 0.68, 0.32}
\definecolor{RoyalBlue}{cmyk}{1, 0.50, 0, 0}
\definecolor{Black}{cmyk}{0, 0, 0, 0}

\hypersetup{
    colorlinks=true, linktocpage=true, pdfstartpage=1, pdfstartview=FitV,
    breaklinks=true, pdfpagemode=UseNone, pageanchor=true, pdfpagemode=UseOutlines,
    plainpages=false, bookmarksnumbered, bookmarksopen=true, bookmarksopenlevel=1,
    hypertexnames=true, pdfhighlight=/O,
    urlcolor=webbrown, linkcolor=RoyalBlue, citecolor=webgreen,
}

\declaretheorem[name=Theorem, numberwithin=section]{thm}
\declaretheorem[sibling=thm, name=Lemma]{lem}
\declaretheorem[sibling=thm, style=definition, name=Remark] {rk}
\declaretheorem[sibling=thm, name=Definition]{defi}

\numberwithin{equation}{section}

\newcommand{\sss}           { \scriptscriptstyle }

\newcommand{\e}{\mathrm{e}}
\newcommand{\indi}{\mathbbm{1}}

\renewcommand{\P}{\mathbb{P}}

\newcommand{\Tcal}{\mathcal{T}}
\newcommand{\Rcal}{\mathcal{R}}
\newcommand{\Ccal}{\mathcal{C}}

\renewcommand{\P}   {\mathbb{P}}

\newcommand{\E}     {\mathbb{E}}

\newcommand{\Rd}    {\mathbb{R}^d}
\newcommand{\Zd}    {\mathbb{Z}^d}

\newcommand{\Pp}        {\mathbb{P}_p}

\newcommand{\Bb}{\mathsf{Bb}}
\newcommand{\trap}{\mathsf{trap}}
\newcommand{\sib}{\mathsf{sib}}
\newcommand{\anc}{\mathsf{anc}}

\newcommand{\eqn}[1]{\begin{equation}#1\end{equation}}
\newcommand{\eqan}[1]{\begin{align}#1\end{align}}

\usepackage[textsize=small]{todonotes}       

\begin{document}
\title{Random walk on barely supercritical
branching~random~walk}
\author{Remco van der Hofstad}
\author{Tim Hulshof}
\author{Jan Nagel}
\address{Department of Mathematics and Computer Science, Eindhoven University of Technology, PO Box 513, 5600 MB Eindhoven, the Netherlands.}
\email{r.w.v.d.hofstad@tue.nl, w.j.t.hulshof@tue.nl, j.h.nagel@tue.nl}
\date{\today}
\begin{abstract} 
Let $\Tcal$ be a supercritical Galton-Watson tree with a bounded offspring distribution that has mean $\mu >1$, conditioned to survive. Let $\varphi_\Tcal$ be a random embedding of $\Tcal$ into $\Zd$ according to a simple random walk step distribution. Let $\Tcal_p$ be percolation on $\Tcal$ with parameter $p$, and let $p_c = \mu^{-1}$ be the critical percolation parameter. 
We consider a random walk $(X_n)_{n \ge 1}$ on $\Tcal_p$ and investigate the behavior of the embedded process $\varphi_{\Tcal_p}(X_n)$ as $n\to \infty$ and simultaneously, $\mathcal{T}_p$ becomes critical, that is, $p=p_n\searrow p_c$. We show that when we scale time by $n/(p_n-p_c)^3$ and space by $\sqrt{(p_n-p_c)/n}$, the process $(\varphi_{\Tcal_p}(X_n))_{n \ge 1}$ converges to a $d$-dimensional Brownian motion. We argue that this scaling can be seen as an interpolation between the scaling of random walk on a static random tree and the anomalous scaling of processes in critical random environments. 
\end{abstract}
\maketitle

\vspace{1em}
{\small

}
\vspace{1em}
\hrule
\vspace{.5em}

\section{Introduction}
Percolation has long stood as a simple and tractable model for random media in physics. Starting with the infinite integer lattice $\Zd$ with nearest-neighbor edges, each edge is kept independently with probability $p \in (0,1)$. The resulting random subgraph already has many properties of real random media, see e.g. \cite{Stau85,Sahi94}. Moreover, the model has a rich mathematical structure. Its main feature is the existence of a \emph{phase transition}: there exists a $p_c (\Zd) \in (0,1)$ such that when $p > p_c$ the model has an infinite connected component or cluster, while it does not when $p < p_c$. The behavior for $p$ close to and at the critical point $p_c$, moreover, has many remarkable features.

Random walk on percolation clusters of the lattice $\Zd$ at or near the critical point has been a central model in both physics and modern probability ever since de Gennes proposed it more than forty years ago \cite{deGe76} (see \cite{BenFri16} for an overview). Computer simulations and non-rigorous studies suggest that the model, which de Gennes dubbed ``the ant in the labyrinth'', has many intriguing features, such as the observation that random walk on large critical clusters exhibits anomalous diffusion. This fact has since been rigorously verified in two dimensions \cite{Kest86b} and in high dimensions \cite{KozNac09}. Yet many questions still remain unanswered. For instance, almost all rigorous results we have are either for the two-dimensional lattice or when the dimension is ``high enough'' (depending on the problem this is typically above either six or ten dimensions \cite{HarSla90a,FitHof17}). Dimensions three to six are terra incognita. And even in high dimensions, where the picture is perhaps most complete, there are still many big open problems. For instance, it is currently unknown what the scaling limit of random walk on large critical clusters is (although there are good conjectures \cite{HarSla00a,Slad02}, on which much progress has been made recently \cite{BenCabFri16b,BenCabFri16a}). 

One subject that is now, after a long line of research, rather well understood for percolation in any dimension, is random walk on the infinite supercritical percolation cluster, i.e., when $p > p_c$. When the temporal and spatial scaling factors are $n$ and $1/\sqrt{n}$, respectively, it has a Brownian scaling limit, just as for random walk on $\Zd$ \cite{KipVar86,DeMetal85,DeMetal89,SidSzn04,BerBis07,MatPia07}. That is, writing $(\tilde X_n)_{n \ge 0}$ for a random walk on an infinite percolation cluster $\Ccal_\infty$ at supercriticality, one obtains a scaling limit of the form
	\begin{align} \label{percolationCLT}
	\big( n^{-1/2} \tilde X_{\lfloor tn\rfloor } \big)_{t\geq 0} \xrightarrow[n \to \infty ]{\mathrm{d}} (\tilde\sigma(p) B_t)_{t\geq 0},
	\end{align}
where $(B_t)_{t\geq 0}$ is a standard $d$-dimensional Brownian motion. 
The diffusion constant $\tilde \sigma (p)$, which describes the typical fluctuations of this Brownian motion, depends on $p$ (and the dimension $d$). Remarkably, it is at this time not even known whether $\tilde \sigma(p)$ is monotonically increasing for $p> p_c$. There is, however, clear evidence that the above scaling limit cannot hold for (an infinite version of) critical clusters, because there it is known that $(n^{-1/6}|\tilde X_n|)_{n \ge 1}$ is a tight sequence of random variables \cite{HeyHofHul14a}. This suggests that either $\tilde \sigma(p) \to 0$ as $p \searrow p_c$, or that there exists a discontinuity for $\tilde \sigma(p)$. Following heuristic arguments (see \cite{biskup2011recent}, Problem 4.21, or \cite{HeyHof17}, Chapter 14), it is conjectured that $\tilde\sigma(p) \sim (p-p_c)$ for $p$ slightly above $p_c$ in high dimensions.

In the current paper we investigate this \emph{barely supercritical regime} for a different model that is believed to be in the same universality class as the high-dimensional setting, namely for random walk on a branching random walk (defined below). The barely supercritical regime here refers to the case where the underlying branching process becomes critical, which we can also parametrize in a natural way as $p$ approaching $p_c$. In Theorem \ref{thm:speed}, we show that this model indeed has $\sigma(p) \sim (p-p_c)$, the scaling conjectured for the high-dimensional percolation case. Furthermore, the underlying tree structure allows for a finer analysis of the process, and as the main result in Theorem \ref{thm:main}, we establish a scaling limit of the random walk after a large number of steps when $p\searrow p_c$ and the underlying branching process thus becomes more and more critical. This requires a bounds that have a high level of uniformity in $p$, which is one of the main challenges in this paper. We show that the Brownian scaling remains observable in the limit where we apply a temporal and spatial scaling with factors $n/(p-p_c)^3$ and $\sqrt{(p-p_c)/n}$, respectively. In a way, this scaling interpolates between the diffusive rescaling by $n^{-1/2}$ of the strictly supercritical case, and the subdiffusive rescaling by $n^{-1/6}$ of the critical case, by letting $p_n$ tend to $p_c$ slower or faster, respectively. To our knowledge this is the first result of such a scaling limit on a sequence of increasingly critical graphs, interpolating between the supercritical and the critical regime.

\subsection{Random walk on a randomly embedded random tree}
Before we proceed with the main results, let us define the model.
\medskip

\textbf{Percolation on Galton-Watson trees.} 
Let $\mathcal{T}$ be a Galton-Watson tree rooted at $\varrho$ with law $\mathbf{P}$, and let $\xi$ be the random number of offspring of the root. Suppose that the tree is supercritical, i.e., $\mu:=\mathbf{E}[\xi]>1$. We write $\Delta$ for the maximal number of children that a single individual in the (unpercolated) tree can have, i.e.,
\[
	\Delta := \sup \{n \, : \, \mathbf{P}(\xi = n) > 0\}.
\] 
We consider percolation on the edges of the tree and let $\mathcal{T}_p$ denote the connected component of the root $\varrho$ in $\mathcal{T}$, when each edge is deleted with probability $1-p\in (0,1)$, independently of every other edge and of $\mathcal{T}$. Then $\mathcal{T}_p$ is again a Galton-Watson tree, whose distribution we denote by $\mathbf{P}_p$. The root of a tree $\Tcal$ chosen according to $\mathbf{P}_p$ now has $\xi_p$ offspring, such that, conditioned on $\{\xi =n\}$, $\xi_p$ is a binomial random variable with parameters $n$ and $p$. 

The percolated tree has mean number of offspring $p \mu$ and is supercritical if and only if $p>p_c :=1/\mu$. In this setting, denote by $\bar{\mathbf{P}}_p={\mathbf{P}}_p(\, \cdot\, |\, |\mathcal{T}|=\infty)$ the distribution of the percolated tree, conditioned on non-extinction. Given a tree $T$ rooted at $\varrho$, let $(X_n)_{n \ge 0}$ be a simple random walk on $T$ started at $\varrho$. Given $v \in T$, we write $P^v_{T}$ for the law of $(X_n)_{n \ge 0}$ with $X_0 =v$, and write $P_T$ if $v =\varrho$. Given a realization of a Galton-Watson tree $\Tcal$, we call $P_{\Tcal}$ the \emph{quenched law} of the random walk, and we call $\mathbb{P}_p := \bar{\mathbf{P}}_p \times P_{\mathcal{T}}$ the \emph{annealed law} of the random walk on the tree (writing $\P^v_p := \bar{\mathbf{P}}_p \times P^v_{\mathcal{T}}$, with the convention that $P_T^v = \delta_{\{(X_n)_{n \ge 0} = (v)_{n \ge 0}\}}$ if $v \notin T$). 
\medskip

\textbf{Spatial embedding: branching random walk.} 
We embed a Galton-Watson tree $\mathcal{T}$ into $\mathbb{Z}^d$ by means of a \emph{branching random walk,} which we will define now:
Let $D$ denote a non-degenerate probability distribution on $\Zd $.
Given a tree $T$ rooted at $\varrho$, set $Z(\varrho)=0$ and assign to each vertex $v\neq \varrho$ of $\mathcal{T}$ an independent random variable $Z(v)$ with law $D$. For any $v \in T$ there exists a unique path $\varrho=v_0, v_1,\dots, v_m=v$. The branching random walk embedding $\varphi = \varphi_T$ is defined such that $\varphi(v) := Z(v_0)+\dots +Z(v_m)$. If we write $(Y_n)_{n \ge 1}$ for a random walk on $\Zd$ with step distribution $D$ started at $0$, then two vertices $v,w \in T$ with the unique path $v = v_0,v_1,\dots,v_k=w$ between them are mapped such that $(\varphi(v_0),\dots,\varphi(v_k))$ has the same law as the translation of $(Y_1,\dots,Y_k)$ by $\varphi(v_0)$, i.e., the marginal of the branching random walk embedding of any simple path in $T$ is a random walk path. This type of branching random walk is also sometimes referred to in the literature as ``random walk indexed by a tree''.
\medskip

\textbf{Random walk on barely supercritical BRW.} 
The process that we consider in this paper is that of $(\varphi_{\Tcal}(X_i))_{1 \le i \le a_n}$ under $\P_{p_n}$ as $p_n \searrow p_c$, i.e., the first $a_n$ steps of a random walk on a barely supercritical infinite tree, embedded by a branching random walk $\varphi$ into $\Zd$.\footnote{Note that this is a different process than if we were to consider a simple random walk on the subgraph of $\Zd$ traced out by the branching random walk, as is for instance the topic of \cite{BenCabFri16b} (for a different kind of trees): random walk on the trace is, in our setting, a vacuous complication, because $\Tcal$ is supercritical, and thus grows at an exponential rate, while $\Zd$ has polynomial growth, so that $\varphi_{\Tcal}(\Tcal) = \Zd$ $\Pp$-almost surely.} 

The main result of this paper is the following scaling limit. In its statement, we regard $\left( \varphi (X_{\lfloor t n\rfloor}) \right)_{t\geq0}$ as a random element in the space $\mathcal{D}_{\mathbb{R}^d}[0,\infty)$, endowed with the Skorokhod topology (see \cite{Bill99}) and the Borel $\sigma$-algebra. 

\begin{thm} \label{thm:main}
Consider random walk $(X_m)_{m \ge 0}$ a percolated Galton-Watson tree with $\mu > 1$, $\Delta < \infty$, and percolation parameter $p_n$ satisfying $p_n\searrow p_c=1/\mu$,
embedded into $\Zd$ with $d \ge 1$ by $\varphi$ whose one-step distribution $D$ satisfies $\sum_{x \in \Zd} x D(x) = 0$ and $\sum_{x \in \Zd} \e^{c \|x\|} D(x) <\infty$ for all $c>0$. For any such sequence $(p_n)_{n\geq 1}$ satisfying $\mathrm{e}^{\delta\sqrt{n}} (p_n -p_c) \to \infty$ for any $\delta>0$ as $n \to \infty$,
	\begin{align*}
	\left( \sqrt{\tfrac{p_n-p_c}{n}} \varphi (X_{\lfloor t n(p_n-p_c)^{-3}\rfloor}) \right)_{t\geq0} \xrightarrow[n\to \infty]{\mathrm{d}} \left( (\kappa \Sigma)^{1/2}  B_t \right)_{t\geq 0}
	\end{align*}
under $\mathbb{P}_{p_n}$, with $(B_t)_{t \ge 0}$ a standard Brownian motion in $\mathbb{R}^d$, $\Sigma$ the covariance matrix of $D$, and $\kappa$ is an explicit constant given in Theorem~\ref{thm:speed} below.
\end{thm}

Note that the process in Theorem \ref{thm:main} is effectively a sequence pairs indexed by $n$: The first element of the pair is the random embedding of a tree percolated at $p_n$ such that, as $n$ tends to infinity, it becomes more and more critical. The second element of the pair is a time-rescaled random walk $\left(X_{\lfloor t n(p_n-p_c)^{-3}\rfloor}) \right)_{t\geq0}$ on this tree, that, as $n$ tends to infinity, takes more steps per unit of time $t$.
The increasing number of steps is needed to see the Brownian motion in the limit, because the closer to criticality the tree is, the more steps it takes before the ballistic nature of the random walk becomes visible.
In Section~\ref{sec:timerescaling} we give a heuristic explanation for why this is the correct order, and also discuss the implications of observing the random walk on different time scales. In Section~\ref{sec:ass} we explain the assumptions of the theorem.

One of the main features of the convergence in Theorem \ref{thm:main} is that it, in some sense,  {\em interpolates} between a random walk in a strictly supercritical and critical environment. Indeed, for $p=p_n>p_c$ fixed, the scaling is the same as in \eqref{percolationCLT} (and also in \eqref{supercritconv} below), while in the barely supercritical regime where $p_n\searrow p_c$, the additional factors $(p_n-p_c)$ bring the scaling close to that on a critical tree. We discuss the latter point in more detail in Section \ref{sec:IICscaling}.

To understand the problem of random walk on a barely supercritical BRW, and to give a clearer context to the main results, let us first discuss the behavior of random walk on a strictly supercritical BRW. 
\medskip

\textbf{Random walk on supercritical BRW.} 
Consider for the moment random walk on BRW under $\mathbb{P}_p$, with $p > p_c$ fixed. Given $v \in T$, let $|v|$ denote the distance of a vertex $v$ to the root.
It is shown in \cite{LyoPemPer96a} that there exists an infinite sequence of \emph{regeneration times} $(\tau_k)_{ k\geq 1}$ such that $|X_n|<|X_{\tau_k}|$ for any $n<\tau_k$ and $|X_n|>|X_{\tau_k}|$ for any $n>\tau_k$, $\P_p$-almost surely. These regeneration times decompose the trajectory of the random walk into independent increments with good moment bounds. Since $\big(\varphi(X_{\tau_{k+1}})-\varphi(X_{\tau_{k}})\big)_{k\geq 1}$ are then independent as well, and this increment is given by the displacement of the embedding random walk after $\tau_{k+1}-\tau_k$ steps, standard arguments allow us to conclude the convergence
	\begin{align} \label{supercritconv}
	\big( n^{-1/2} \varphi(X_{\lfloor tn\rfloor }) \big)_{t\geq 0} \xrightarrow[n \to \infty ]{\mathrm{d}} (\Sigma(p)^{1/2}B_t)_{t\geq 0}
	\end{align}
in distribution under $\mathbb{P}_p$, with $(B_t)_{t \ge 0}$ a standard $d$-dimensional Brownian motion and $\Sigma(p)$ denoting a covariance matrix. (The proof of Theorem~\ref{thm:main} in Section~\ref{sec:mainpf} below also establishes \eqref{supercritconv}.) Further writing $\Sigma$ for the covariance matrix associated with $D$, i.e., $\Sigma_{ij}:=\sum_{x} x_ix_j D(x)$, and $v^{\mathsf{T}}$ for the transpose of $v \in \Zd$, we can determine that the covariance matrix $\Sigma(p)$ is given by
	\begin{align} \label{speed-var}
	\Sigma(p)= \lim_{n\to \infty} \frac{1}{n} \mathbb{E}_p[  \varphi(X_{n})\varphi(X_{n})^{\mathsf{T}} ] = \lim_{n\to \infty} \frac{1}{n} \mathbb{E}_p[ |X_n| ] \, \Sigma,
	\end{align} 
where the last equality follows by conditioning on the distance of $X_n$ from the root, so that $\varphi(X_n)$ can be written as a sum of i.i.d.\ increments with covariance $\Sigma$.

It is a well-known result \cite{LyoPemPer96a} that the \emph{effective speed}
	\begin{align} \label{speed}
	v(p) := \lim_{n\to \infty} \frac{|X_n|}{n} 
	\end{align} 
exists $\mathbb{P}_p$-almost surely and is non-zero when $p > p_c$. This implies that the covariance matrix in the supercritical setting is given by
	\begin{align}\label{supercritvar}
	\Sigma(p)= v(p) \Sigma.
	\end{align}
We know that $v(p)>0$ for any fixed $p > p_c$. It is not known what happens to $v(p)$ when $p \searrow p_c$. The following result establishes that $v(p) \to 0$, and at what rate it does so. 

For $k \in \mathbb{N}$, define the $k$-th factorial moment of the (unpercolated) offspring distribution as
	\[
	m_k := \textbf{E}\left[\prod_{i=0}^{k-1} (\xi - i)\right].
	\]
\begin{thm}[The speed of random walk on a barely supercritical tree]\label{thm:speed} Consider a simple random walk $(X_n)_{n \ge 0}$ with $X_0 = \varrho$ on a Galton-Watson tree with  $m_1 >1$ and $m_3 < \infty$. Then,
	\[
	\lim_{p\searrow p_c} \frac{v(p)}{(p - p_c)^2} =  \frac{m_1^4}{3 m_2} =: \kappa .
	\]
\end{thm}
We prove this theorem in Section~\ref{sec:speed} by analyzing the Taylor expansion of the extinction probability in terms of the generating function. Note that \eqref{supercritvar} and Theorem~\ref{thm:speed} imply that the variance of the limit in Theorem \ref{thm:main} is given by the limit of the supercritical variance, appropriately rescaled. That is,
	\begin{align}
	\lim_{p\searrow p_c} \frac{\Sigma(p)}{(p-p_c)^2} = \kappa \Sigma,
	\end{align}
where the scaling in $(p-p_c)$ is thus same as the conjectured scaling for random walk on a high-dimensional barely supercritical infinite percolation cluster mentioned above.

\subsection{On the scaling factors in Theorem \ref{thm:main}}
\label{sec:timerescaling}

Heuristically, the fact that rescaling space by $\sqrt{(p-p_c)/n}$ and time by $n/(p-p_c)^3$ yields a Brownian limit can be understood as follows: The duality principle tells us that any infinite supercritical GW-tree can be decomposed into a supercritical GW-tree with no leaves (which we call the ``backbone tree'') to which, at each vertex, are attached a random number of i.i.d.\ subcritical GW-trees (which we call ``traps''), see Section~\ref{sec:prelim0} below for details. 

The furcations or branch points in the backbone tree are then separated by path-like segments with lengths of order $(p-p_c)^{-1}$. Simple random walk on a GW-tree with no leaves is transient because each time the walk reaches a furcation, it is more likely to take a step away from the root than towards it. 
If we look at the sequence of pairs of an $n$ step walk $(X_{k})_{0 \le k \le n}$ on the tree $\Tcal_{p_n}$, then if the walk sees a growing number of furcations as $n$ tends to infinity, a limit of the form \eqref{supercritconv} can be expected.
Together, the rate at which $p_n$ tends to $p_c$, and the length of the long line segments in the backbone tree, thus imply that we need to rescale the tree by $(p_n-p_c)/n$ for the walk to see a linearly growing number of furcations. Since the embedding of the tree into $\Zd$ is diffusive, this means that we need to rescale space by $\sqrt{(p_n-p_c)/n}$.

If the walk were restricted to the backbone tree, then it would take the walk order $(p-p_c)^{-2}$ steps to visit the next furcation, because that is how long it takes for a random walk to travel distance $O((p-p_c)^{-1})$ on a line.
Between two furcations, we can expect a tight number of large traps of size $O((p-p_c)^{-2})$, each visited $O((p-p_c)^{-1})$ times. An  application of electrical network theory for Markov chains tells us that the time to exit such a large trap is of order $(p-p_c)^{-2}$. The time spent in traps between furcations thus accumulates to $O((p-p_c)^{-3})$. Since the above spatial scaling implies that the number of furcations for the walk restricted to the backbone grows linearly in $n$, to see a scaling limit like \eqref{supercritconv}, we need to rescale time by a factor $n/(p-p_c)^{3}$ for the same to be true about the walk on the full tree.

We can expect rather different behaviour when the number of furcation points visited by the random walk remains bounded, i.e., in the setting where we consider the walk at a time scale $a_n$ such that $a_n(p_n-p_c)^3$ remains bounded. In this case, the walk will barely observe that the BRW is supercritical. In particular, when $a_n(p_n-p_c)^3\rightarrow 0$, the random walk does not observe the supercritical nature at all, and basically observes a {\em critical} BRW conditioned to survive forever. As discussed in more detail in the next section, we do not expect that the limit in Theorem \ref{thm:main} holds for these time scales.

\subsection{Interpolation of supercritical and critical scaling}
\label{sec:IICscaling}
As mentioned, we view the limit in Theorem~\ref{thm:main} as interpolating between the supercritical and the critical regime. 
This behavior becomes more visible when a nonstandard scaling is applied to the random walk $(X_n)_{n \ge 0}$: 
Writing $a_n := n(p_n-p_c)^{-3}$, an equivalent representation of the convergence in Theorem \ref{thm:main} then reads
	\begin{align}\label{e:thmrestate}
	\left( (p_n-p_c)^{-1} a_n^{-1/2} \varphi (X_{\lfloor t a_n\rfloor}) \right)_{t\geq0} \xrightarrow[n\to \infty]{\mathrm{d}} \left( (\kappa \Sigma)^{1/2}  B_t \right)_{t\geq 0} .
	\end{align}

Consider a sequence $(p_n)_{n \ge 1}$ converging to the critical value $p_c$ very slowly. The prefactor $(p_n-p_c)$ then plays a negligible role compared to the prefactor $a_n^{-1/2}$. Indeed, completely omitting $(p_n-p_c)^{-1}$ yields exactly the  same diffusive scaling factors as those of the random walk on the supercritical percolation cluster as in \eqref{percolationCLT}, or of the random walk on a supercritical BRW as in \eqref{supercritconv}.  In this setting, the random walk on the BRW behaves ``almost supercritical''. 

Now consider a sequence $(p_n)_{n \ge 1}$ converging very quickly to $p_c$. This gives rise to behavior of the random walk that is ``almost critical''. To observe this, note that in this setting the dominant factor in the time rescaling by $a_n$ is $(p_n-p_c)^{-3}$ and the  
 convergence in Theorem \ref{thm:main} may be rewritten as 
	\begin{align}
	\label{almost-critical}
	\left(n^{-1/3} a_n^{-1/6} \varphi (X_{\lfloor t a_n\rfloor}) \right)_{t\geq0} \xrightarrow[n\to \infty]{\mathrm{d}} \left( (\kappa \Sigma)^{1/2}  B_t \right)_{t\geq 0} .
	\end{align}
This time, the prefactor $n^{-1/3}$ plays a negligible role compared to $a_n^{1/6}$ when $n\rightarrow \infty$. The rescaling we see here, of time by $a_n$ and of space by $a_n^{-1/6}$, is characteristic for processes on high-dimensional critical random objects, as we discuss now.

A key example of such critical processes is a random walk on a critical BRW conditioned to survive forever. There, it is known that $a_n^{-1/6} \varphi (X_{\lfloor t a_n \rfloor})$ is tight (\cite{Kest86b} proves for random walk on the critical Galton-Watson tree conditioned to survive that $a_n^{-1/3} |X_{\lfloor t a_n \rfloor}|$ is tight, and tightness of the associated BRW then easily follows from the diffusive nature of the random walk embedding $\varphi$).

For percolation a similar result is proved in \cite{HeyHofHul14a} (which builds upon \cite{BarJarKumSla08, KozNac09}). We expect that the scaling limit of critical BRW conditioned to survive forever is superBrownian motion (SBM) conditioned to survive forever (or equivalently, SBM with an immortal particle, see \cite{Evan93}). See \cite{Hofs06a} for several examples of critical models that converge to this superprocess, including the incipient infinite cluster of oriented percolation. For a random walk $(Y_n)_{n \ge 0}$ on such a critical high-dimensional structure, we further expect that the limit
	\begin{equation}\label{e:IICBM}
	\left(a_n^{-1/6} Y_{\lfloor t a_n \rfloor} \right)_{t \ge 0} \xrightarrow[n\to \infty]{\mathrm{d}} \left( \sqrt{\bar{\kappa}}  B^{\mathrm{IIC}}_{\bar \nu t} \right)_{t\geq 0} 
	\end{equation}
\color{black}
holds, where $(B^{\mathrm{IIC}}_{\bar \nu t})_{t\geq 0}$ is a Brownian motion on the trace of SBM conditioned to survive forever (mind that this combined process has not yet been defined rigorously, as far as we are aware). 
 
A full scaling limit for random walk on a critical Galton-Watson tree conditioned to be infinite was recently proved by Athreya, L\"ohr, and Winter \cite{athreya2017invariance}, with a rescaling that, after embedding the tree, exactly corresponds to the one in \eqref{e:IICBM}. 
Also recently, Ben-Arous, Cabezas, and Fribergh \cite{BenCabFri16a} proved that the same scaling limit holds in high dimensions (with different $\bar \kappa$ and $\bar \nu$) for the model where the random walk moves on the trace of the critical BRW in $\Zd$ (rather than directly on the tree, as in the current paper). 

Croydon \cite{Croy09} showed that a random walk $(Y_m)_{m \ge 1}$ on a sequence of critical trees $(\Tcal_n)_{n \ge 1}$ conditioned to have size $n$, and randomly embedded by $\varphi$ (with $\Sigma = \mathrm{Id}$), has the scaling limit
	\begin{equation}\label{e:ISEBM}
	\left(b_n^{-1/6} \varphi (Y_{\lfloor t b_n \rfloor}) \right)_{t \ge 0} \xrightarrow[n\to \infty]{\mathrm{d}} \left( \sqrt{\bar{\kappa}}  B^{\mathrm{ISE}}_{\bar \nu t} \right)_{t\geq 0}, 
	\end{equation}
with $b_n = n^{3/2}$ and where $B^{\mathrm{ISE}}$ is a Brownian motion on the \emph{Integrated superBrownian Excursion,} or ISE. The ISE is now understood to be a canonical random object. It can be viewed as a Brownian embedding of a continuum random tree of a fixed size (cf.\ \cite{Aldo91}), or as the scaling limit of high-dimensional lattice trees, cf.\ \cite{DerSla98}. Here, the time scale $b_n = n^{3/2}$ is asymptotically the minimal time scale that ensures that the random walk explores a non-vanishing fraction of the tree, and thus notices that it moves on a finite structure. We expect that if we set $b_n \ll n^{3/2}$ instead, using the same scaling as in \eqref{e:ISEBM} will yield a limit equal to that in \eqref{e:IICBM}. This is because the random walk only explores a small part of the embedded tree that is close to its starting point, and this local perspective is the same as the local perspective of a critical branching process conditioned to survive forever as in \eqref{e:IICBM}.

It would be interesting to consider the regime when $a_n(p_n-p_c)^3$ converges to a constant.
We expect that when $a_n \to \infty$ but $a_n (p_n - p_c)^3 \to 0$, the limit is the same as in \eqref{e:IICBM}.  In the case where $a_n(p_n-p_c)^3$ converges to a positive constant, it is not clear to us what the scaling limit should be (but we expect that it will not be equal to that in \eqref{e:IICBM}).

Heuristically, the case where $a_n (p_n-p_c)^3$ converges can be viewed as the setting in which the random walk sees only a finite number of furcation points of the tree.
Returning to the point of view of \eqref{e:thmrestate}, where the time scale is set to be $a_n = n (p_n - p_c)^{-3}$, which diverges as $n \to \infty$ for any $p_n \to p_c$, we see that this choice can be understood as a time scale such that the random walk sees a growing number of furcation points.
This means that it can ``feel'' the drift that arises from the fact that the barely supercritical tree is in fact supercritical.
The factor $n$ in $a_n$ does not seem necessary for that, and we thus expect that any diverging $a_n$ leads to a Brownian limit. In particular, we do not see the need for a restriction on $p_n$ here. The reason we do have restrictions on $(p_n)$ in the assumptions of Theorem \ref{thm:main} is different: we need to control the displacement in traps, which requires $n=a_n(p_n-p_c)^3$ to grow sufficiently fast compared to $(p_n-p_c)^{-1}$, as formulated in the condition that $\mathrm{e}^{\delta\sqrt{n}} (p_n -p_c) \to \infty$ for any $\delta>0$ as $n \to \infty$.

\subsection{About the assumptions in Theorem~\ref{thm:main}}\label{sec:ass}
(1) We only assume that there exists a finite maximal degree $\Delta$ to simplify the proofs below. We believe that the proof can be modified (albeit with some lengthy computations) to admit any offspring distribution whose generating function $f(s)  = \mathbf{E}[s^\xi]$ has derivatives that satisfy $f^{(n+1)}(s) \le C nf^{(n)}(s)$ for all $n \ge 1$ and $0 \le s \le 1$ (we do not know if this condition is necessary). Examples of distributions that satisfy this condition are Poisson and Geometric.

(2) The assumption on $D$ that $\sum_{x \in \Zd} x D(x) = 0$ is necessary, otherwise the BRW would have a drift, which would require a different analysis and yield a different limit. 

(3) The assumption that $\sum_{x \in \Zd} \e^{c\|x\|} D(x) <\infty$ prevents the BRW from making very large jumps. The strength of the assumption is for techical reasons. We believe that the result should remain valid when $\sum_{x \in \Zd} \|x\|^4 D(x) <\infty$. It is known that the behavior of BRW (and other statistical mechanical models) alters dramatically when this restriction is relaxed further, see e.g.\ \cite{Kest95, JanMar05,Huls15}. In particular, macroscopically long edges in the BRW embedding are not sufficiently rare that the random walk can avoid them w.h.p. When the random walk crosses such a long edge, it will likely only spend a very short time at the other end of the edge, which should not affect the random walk's diffusivity, but we expect that it will affect the continuity of the limiting process. We thus conjecture that when $D$ is such that $\alpha=\sup\{a:\, \sum_{x \in \Zd} \|x\|^a D(x) <\infty\} $ satisfies $2 \leq \alpha < 4$, we may still be able to see convergence to a Brownian motion in finite-dimensional distributions, but we would not be able to achieve such convergence in the Skorokhod topology. When $\alpha < 2$ the long edges become so common that the random walk will cross them frequently enough that also the finite-dimensional distributions will be affected, see e.g. \cite{HeyHofSak08,CheSak15}. When $\alpha < 2$, we expect to see convergence to an $\alpha$-stable motion, but again only for the finite-dimensional distributions.

(4) We require that $(p_n)_{n \ge 1}$ satisfies that for any $\delta>0$ we have $ \e^{\delta\sqrt{n}}(p_n -p_c) \to \infty$, so that we may apply a result of Neuman and Zheng regarding the maximal displacement of subcritical BRW \cite{NeuZhe17}, which we crucially use to estimate the size of the traps after embedding into $\Zd$. We expect the result to remain valid when thus condition is relaxed (see Section \ref{sec:IICscaling} below).

\subsection{About the proofs}
The proof of Theorem~\ref{thm:main} is, at its core, close to the classical proofs of Brownian limits for random walk in random environment. In particular, we construct a sequence of regeneration times to find (almost) independent increments, and follow the standard approach from there. 

What is new about our proof is that we require all the necessary moment bounds to hold uniformly for $p \in (p_c, p_d)$ for some $p_d > p_c$. To obtain such uniform bounds we need to perform a careful analysis of the structure of the backbone tree. In particular, we need to take into account that furcations occur only on a length scale of order $(p-p_c)^{-1}$ and the influence of the traps grow as $p$ tends to $p_c$. 

Our regeneration structure is based on that of \cite{SznZer1999}, but is specifically designed to account for the length rescaling of the tree, due to furcation points growing far apart. Classical regeneration times would require the random walk to never visit a parent vertex, but this probability tends to zero as the tree becomes critical. 
Instead, inspired by the regenerations of \cite{GanMatPia2012}, we construct a sequence of regeneration times $(\tau_k)_{k \ge 1}$ that allows the random walk to backtrack a distance of order $(p-p_c)^{-1}$ after a regeneration time, to increase the density of the regeneration times to the right scale. This relaxation comes at the cost of (1) having a small and localized intersection between the past and future of the walk at $\tau_k$, and (2) making $(\tau_{k+1} -\tau_k)_{k\geq 1}$ a stationary and $1$-dependent sequence, rather than an i.i.d.\ sequence (but these are not serious complications). From this new construction of regeneration times we are then able to derive all the necessary moment bounds on the regeneration times and distances uniformly for $p \in (p_c,p_d)$. We stress the fact that the distribution of $(\tau_{k+1} -\tau_k)_{k\geq 1}$ depends sensitively on $p_n$, and thus it might be appropriate to write $(\tau_{k+1}^{\sss(n)} -\tau_k^{\sss(n)})_{k\geq 1}$ instead. We omit it to simplify notation. Recalling from \eqref{e:thmrestate} that we investigate the random walk on a time scale $a_n=n(p_n-p_c)^{3}$ is it worthwhile however to note that $n$ can be interpreted as the order of the number of regenerations before time $a_n$.
The fact that $n\rightarrow \infty$ allows us to use classical laws of large numbers for the sequence $(\tau_k)_{k\geq 1}$.

Besides controlling the regeneration structure, we also need exponentially tight control over the displacement of the random walk in the traps (Lemma~\ref{lem:nzlemma}), and a bound on all moments of the size of the trace of the random walk on the backbone tree (Lemma~\ref{lem:BBTbounds}). Uniformity in $p$ is again an important requirement here, since we apply these bounds to $p=p_n\searrow p_c$. If, rather than looking at $\varphi(X_{\lfloor a_nt \rfloor})$, we would consider the position of the walker only at regeneration times, i.e., $\varphi(X_{\tau_n})$, then Theorem \ref{thm:main} would holds without any assumption on $p_n$ except for $p_n\searrow p_c$. This follows from the first three steps of the proof in Section \ref{sec:mainpf}. The restriction $\mathrm{e}^{\delta\sqrt{n}} (p_n -p_c) \to \infty$ is needed to control the maximal displacement between regeneration times.
\medskip

The proof of Theorem~\ref{thm:speed} in Section~\ref{sec:speed} is a straightforward Taylor expansion of a generating function for the effective speed of simple random walk on a supercritical GW-tree derived in \cite{LyoPemPer96a}, applied to our setting. It is, however, interesting to note that in two seemingly far removed parts of the proof we see a term involving the third moment of $\xi$ arise, but these terms cancel perfectly, leaving us with an expression of the asymptotic speed that only involves the first two moments. It is unclear to us why the third moment should drop out like this.

An alternative proof of the scaling as in Theorem  \ref{thm:speed} may be obtained as a byproduct of the regeneration structure. Although this method does not give the precise limit, it is more robust and could be applied to more general models (see Remark \ref{rem:speedbounds}).

\subsection{The structure of the paper}
We start by stating some preliminary lemmas about the structure of slightly supercritical trees in Section \ref{sec:prelim0}. We conclude that the typical length scale of the trees is $(p - p_c)^{-1}$, and introduce a rescaling by this length that we use throughout the paper. In Section \ref{prelim}, we state some further lemmas with preliminary bounds on escape probabilities of random walk on such trees. 
The proofs of these lemmas can be found in Section~\ref{sec:aprioriproofs}.

In Section~\ref{sec:regen} we define a new regeneration structure, which decomposes the trajectory into one-dependent increments. We state moment bounds on the regeneration distances in Lemma~\ref{lem:moments} (which we prove in Section~\ref{sec:regdist}), and moment bounds on times between regenerations in Lemma~\ref{lem:moments2} (which we prove in Section~\ref{sec:regtime}).

In Section~\ref{sec:mainpf} we give the main steps of the proof of Theorem~\ref{thm:main}. Utilizing the decoupling effect of the regenerations, we first prove that the limit holds if we consider the increments of the walk between regeneration times. For this we need the moment bounds on the regeneration distances of Lemma~\ref{lem:moments}, and also bounds on the inter-regeneration times of Lemma~\ref{lem:moments2}. Because this process does not take into account where the random walk is between regeneration times, and the regeneration times by definition never occur in the traps of the tree, the final step of the proof is to show that the random walk is not able to walk great distances in traps. For this we use a recent result by Neuman and Zheng \cite{NeuZhe17} on the maximal displacement of subcritical BRW, which yields a bound on the maximal displacement of the traps in Lemma~\ref{lem:nzlemma}. The proof of Lemma \ref{lem:nzlemma} can be found in Section \ref{sec:maxdisp}. We also need bounds on the size and shape of the trace of the backbone of the BRW, in Lemma~\ref{lem:BBTbounds} (which we prove in Section~\ref{sec:BBT}). 

Finally, Section \ref{sec:speed} contains the expansion of the speed on the tree and thus proves Theorem \ref{thm:speed}.

\section{Preliminaries: the shape of the tree}\label{sec:prelim0}
In this section we establish some facts about the percolated trees, such as their growth rate. We also discuss a useful decomposition of the tree.

To simplify notation we frequently drop the subscript $p$. We write $c$ and $C$ for generic constants whose value may change from line to line. Central bounds appearing in the a-priori estimates in this section are denoted by $a_i$ and $c_i$. The constants $a_i,c,C,$ and $c_i$ may depend on the original offspring distribution and on $d$, but we stress that they are independent of $p$. Since we are interested in $p$ close to $p_c$, we will state many results only for $p\in (p_c,p_d)$ for some unspecified $p_d>p_c$. 

Given a tree $T$ rooted at $\varrho$ and $v,w \in T$, $v\neq w$, we say $w$ is a \emph{descendant} of $v$ if the unique self-avoiding path from $\varrho$ to $w$ passes through $v$, and we call $w$ an \emph{ancestor} of $v$ if the path from $\varrho$ to $v$ passes through $w$. We write $T_v$ for the subtree of $T$ rooted at $v$ induced by $v$ and all descendants of $v$. 
We say a vertex $v$ has an \emph{infinite line of descent} if $T_v$ is an infinite tree.
We can decompose any infinite tree $T$ into its \emph{backbone tree} $T^{\Bb}$, the tree induced by all vertices with an infinite line of descent, and the forest $T \setminus T^{\Bb}$, which consists only of finite trees. Here, we denote by the difference $T_1\setminus T_2$ of a tree $T_1$ and a subgraph $T_2$ of $T_1$ the graph obtained by removing from $T_1$ all edges of $T_2$ and all then isolated vertices. This is equivalent to removing from $T_1$ all vertices that are only part of edges in $T_2$. Note that by this definition, $T\setminus T_v$ is the subtree of $T$ with all descendants of $v$ removed, but still containing $v$. Given $v \in T^{\Bb}$, we call the connected component of $v$ in $T \setminus T^{\Bb}$ the \emph{trap at $v$} and denote it by $T_v^{\trap}$.

Lyons \cite[Proposition 4.10]{Lyon92} gives an explicit description of this decomposition for Galton-Watson trees using the ``duality principle'' (cf.\ \cite[Chapter 12]{AthNey72}). For $f$ a generating function of an offspring distribution of a Galton-Watson tree with extinction probability $q\in (0,1)$, set 
\begin{align} \label{treedecomposition}
\hat f (s) := \frac{f((1-q)s+q)-q}{1-q} \qquad \text{ and } \qquad  f^*(s) := \frac{f(qs)}{q} .
\end{align}
A supercritical Galton-Watson tree $\mathcal T$ with generating function $f$ conditioned on survival can be generated by sampling a tree $\mathcal{T}^{\Bb}$ with generating function $\hat f$, and then adding to every vertex $v\in \mathcal{T}^{\Bb}$ a random number $U_v$ of edges, and to the other ends of those edges independent Galton-Watson trees $\Tcal^*_{i}$, $i \in \{1,\dots,U_v\}$ with generating function $f^*$. While $\Tcal^{\Bb}$ has no leaves, the trees generated by $f^*$ go extinct with probability 1. Conditionally on $\mathcal{T}^{\Bb}$, the $(U_v)_{v\in \mathcal{T}^{\Bb}}$ are independent (but not identically distributed). The marginal distribution of $U_v$ can be characterized as follows: for $v \in \Tcal^{\Bb}$ write $\delta_v =  \deg_{\Tcal^\Bb}(v)-1$ and write $f^{(n)}(s)$ for the $n$-th derivative of $f(s)$. Then,
\begin{equation}\label{e:Uvdef}
	\mathbf{E}[s^{U_v}] = \frac{f^{(\delta_v)}(qs)}{f^{(\delta_v)}(q)}.
\end{equation}
Then $\mathcal{T}^{\trap}_v$ is the random subtree consisting of $v$, the $U_v$ edges attached to $v$, and the finite trees generated by $f^*$ attached to the edges. 
\begin{lem}[Properties of the generating functions]\label{lem:gf} 
Consider a Galton-Watson tree with generating function $f$ and $m_1 >1$ and $m_3 < \infty$. There exist a $p_d > p_c$ and constants $c_0 \in (0,1)$ and $c_1> 0$ such that for any $p \in (p_c,p_d)$,
	\eqn{
	\label{equal-lem2.1}
	f'_p(0) \geq c_0, \qquad \mu_p = \hat \mu_p = f'_p(1) = p/p_c,
	}
and, for $p\searrow p_c$,
	\begin{align}
	\label{asy-lem2.1}											
	q_p 				& = 1- c_1(p-p_c)(1+o(1)), \notag\\
	\mu_p^* = f'_p(q_p) 	&=  1-m_1(p-p_c)(1+o(1)),\\
	\hat f_p''(0) 			& = 2m_1 (p-p_c)(1+o(1)).\notag
	\end{align}
\end{lem}
The proof of this lemma is standard and can be found in Section~\ref{sec:speed}.

\begin{rk}[Topology of barely supercritical tree, backbone tree and traps as $p_n \searrow p_c$] \label{rem:treeshape}
From the above lemma we can infer that, heuristically speaking, the percolated GW-tree $\Tcal  \sim \bar{\mathbf{P}}_p$,  conditioned on survival, and its decomposition look as follows for $p$ close to $p_c$: Both $\Tcal$ and $\Tcal^{\Bb}$ grow at rate $f_p'(1) = p / p_c$, but only a fraction $p-p_c$ of the vertices in $\Tcal$ is contained in $\Tcal^{\Bb}$. By the asymptotics above,  
\begin{align*}
\mathbf{P}(\mathrm{deg}_{\mathcal{T}^{\Bb}}(\rho) \geq 3) = 1- \hat f'_p(0)-\tfrac{1}{2}\hat f''_p(0) = o(p-p_c) , 
\end{align*}
so we are unlikely to see any vertices in $\Tcal^{\Bb}$ with out-degree three or more up to a height of order $(p - p_c)^{-1}$. Up to this height, the tree looks like a binary tree where each edge has been replaced by a path whose length is distributed as an independent geometric random variable with a parameter of order $p - p_c$. The trap at any vertex $v$ consists of a random number of i.i.d.\ subcritical Galton-Watson trees with mean offspring distribution $\mu_{p}^* \approx 1-c_1 (p-p_c)$, whose expected size and depth are both known to be of order $(p -p_c)^{-1}$. The distribution of the traps, however, is such that $\P(|\Tcal^\trap_v| \ge A (p-p_c)^{-2} ) \approx \frac{p-p_c}{A} \e^{- A/2}$, which implies that among $O((p - p_c)^{-1})$ typical traps, almost all traps will be very small (or non-existent), while a tight number of them are macroscopically large, having a size of the order of $(p-p_c)^{-2}$  and depth of the order $(p-p_c)^{-1}$. Although we do not use any of the computations of this heuristic directly in the proofs that follow, we do use them often as guiding principles.
\end{rk}

For a tree $ T$ with backbone tree $T^{\Bb }$, let $G'_m(T) = \{v\in T^{\Bb }\, : \, |v|=m\}$ be the backbone-tree vertices in generation $m\geq 0$. As discussed above, the relevant spatial scaling factor of the random tree will be of order $p-p_c$, so throughout this paper we will often consider the backbone tree only at generations that are multiples of $L/(p-p_c)$, for some $L\geq 1$ to be determined later. We call these generations the \emph{$L$-levels of $T$,} and for $m \in \mathbb{N}$ write
\begin{align} \label{modgen}
G_{[m]} := G_{[m]}(T) := G'_{m \lfloor L/(p-p_c)\rfloor} ( T )
\end{align}
for the $m$-th $L$-level, and write $G_{[0]} := \{\varrho\}$, the root of the tree.
We tacitly ignore the dependency on $L$ (and $p$) in the definition of $G_{[m]}$, the reason being that we soon fix the value of $L$. Roughly speaking, the choice of $L$ will be such that with good probability (uniform in $p-p_c$), the random walk starting from a vertex $v \in G_{[m]}$ never hits $G_{[m-1]}$. 
For $m<n$, let 
\begin{align} \label{bbsegment}
T^{\Bb }_{[m,n]} := \{v\in {T}^{\Bb }\,:\, m \lfloor L/(p-p_c)\rfloor \leq |v|\leq n \lfloor L/(p-p_c)\rfloor \}
\end{align} 
denote the forest segment of the backbone tree between the $m$-th and $n$-th $L$-level. Let
\begin{align} \label{treesegment}
T_{[m,n]} := \bigcup_{v\in T^{\Bb}_{[m,n]}} {T}^{\trap}_v
\end{align}  
denote the same forest segment with all the traps attached. 

The next lemma shows that on the scale of a single $L$-level, the backbone tree $\mathcal{T}^{\Bb }$ looks like a non-degenerate tree:
\begin{lem}[The size of the first $L$-level]\label{lem:bb}
For any $L > 0$ there exists $p_d>0$ and constants $0<a_1\leq a_2<1$ and $a_3<\infty$, such that for all $p\in (p_c,p_d)$
\begin{align*}
a_1\leq \bar{\mathbf{P}}_p(|G_{[1]}(\mathcal{T})|=1)\leq a_2 \quad \text{ and }\quad \bar{\mathbf{E}}_p[|G_{[1]}(\mathcal{T})| ]\leq a_3 .
\end{align*}
\end{lem}

\proof
By \eqref{treedecomposition} we know that the backbone tree $\Tcal^{\Bb}$ is a Galton-Watson tree with generating function $\hat{f}(s)$,
so the probability that a vertex in the backbone tree has exactly one child in $\Tcal^{\Bb}$ is $\hat p_1 = f'(q_p)$. By Lemma~\ref{lem:gf}, $f'(q_p) = 1-c_2(p-p_c) (1+o(1))$ as $p \searrow p_c$.
Therefore,
\begin{align}\label{e:unifprob}
\lim_{p\searrow p_c} \bar{\mathbf{P}}_p(|G_{[1]}(\mathcal{T})|=1) = \lim_{p\searrow p_c} \left( 1-c_2(p-p_c)(1+o(1))\right)^{L/(p-p_c)} = \mathrm e^{-c_2 L}.
\end{align}
This proves the first part of the lemma. 

For the second part, recall that $\Tcal_p$ has offspring mean $\mu_p = \mu p$ and that $p_c = 1/\mu$, so that
\begin{align*}
\bar{\mathbf{E}}_p[|G_{[1]}(\mathcal{T})| ] = \mu_p^{\lfloor L/(p-p_c) \rfloor} = (1+\mu(p-p_c))^{\lfloor L/(p-p_c)\rfloor } ,
\end{align*}
which is uniformly bounded as $p\searrow p_c$. \qed


\section{Preliminaries: escape probabilities}\label{prelim}
In this section we establish three useful bounds on the probability that the random walk escapes to infinity before returning to the previous $L$-level. Recall that $(X_n)_{n\geq 0}$ denotes the random walk on $\mathcal{T}$.

For $A$ a set of vertices of $\mathcal{T}$, let $\eta(A)=\inf \{n\geq 0\,:\, X_n \in A\}$ be the \emph{hitting time} of the set $A$. For
$m \ge 0$, we denote by $\eta_m := \eta(G_{[m]}(\mathcal{T}))$ the hitting time of the $m$th $L$-level. 
A crucial step in defining the regeneration times is the following uniform bound for the probability that the random walk backtracks $\lfloor L/(p-p_c)\rfloor$ steps in tree distance: 

\begin{lem}[Annealed escape probability] \label{lem:backtrackprob}
There exists an $L_0\geq 1$ depending only on the original offspring distribution, such that for any $L\geq L_0$, for any tree $T$ with $G_{[1]}(T)\neq \varnothing$ and any $v\in G_{[1]}(T)$, 
\begin{align*}
\mathbb{P}_p^v(\eta_0=\infty \mid \mathcal{T}_{[0,1]} = T_{[0,1]} ) \geq \tfrac{1}{3} .
\end{align*} 
\end{lem}
We prove this lemma in Section~\ref{sec:backtrack}.

In what follows, fix an $L$ such that the estimate in Lemma \ref{lem:backtrackprob} holds. 

It would simplify our argument substantially if we could get a quenched version of Lemma \ref{lem:backtrackprob}. Unfortunately, such a quenched version does not hold. We do, however, have the following statement, which shows that the backtracking probabilities are small with a high probability, even when we additionally delete all outgoing edges at the starting point $v \in G_{[1]}(\Tcal)$, so that the walker has to escape to infinity via a different vertex in $G_{[1]}(\mathcal{T})$:

\begin{lem}[Quenched indirect escape probability] \label{lem:backtrackprob2}
There exists a function $h$ independent of $p$ with $h(\alpha)\to 0$ as $\alpha \to 0$, such that 
\begin{align*}
 \bar{\mathbf{E}}_p \left[ \sum_{v\in G_{[1]}(\mathcal{T})} \mathbbm{1}_{\{ P^{v}_{\mathcal{T}\setminus \mathcal{T}_{v}}(\eta_{0}=\infty) < \alpha  \}}\right] \leq a_2 + h(\alpha) ,
\end{align*}
with $a_2$ as in Lemma \ref{lem:bb}. 
 \end{lem}
We prove this lemma in Section~\ref{sec:backtrack2}.

For a vertex $v\in G_{[m]}(T)$ with $m\geq 1$, let $\anc(v)\in G_{[m-1]}(T)$ denote the \emph{$L$-ancestor} of $v$, its ancestor in $L$-level $m-1$, and define the number of \emph{$L$-siblings} of $v$ as
\begin{align}
\sib(v) := |G_{[1]}(T_{\anc(v)})|-1.
\end{align}
Note that $\anc(v)$ and $\sib(v)$ are only defined for vertices on the backbone tree, and that $\sib(v)$ only counts the siblings on the backbone tree. The following bound shows that we have a uniformly bounded probability of hitting vertices that have no $L$-siblings:

\begin{lem}[Escape probability on thin parts of the backbone tree] \label{lem:onlychild}
There exists a constant $a_4>0$ such that for any tree $T$ with $G_{[1]}(T)\neq \varnothing$ and any $v\in G_{[1]}(T)$, 
\begin{align*}
\mathbb{P}^v_p( \eta_2<\eta_0, \sib(X_{\eta_2})=0 \mid \mathcal{T}_{[0,1]} = T_{[0,1]} ) \geq a_4 .
\end{align*} 
\end{lem}
We prove this lemma in Section~\ref{sec:onlychild}. 

\section{The regeneration structure}\label{sec:regen}

Regeneration times are a classical tool to decouple the increments of random walks in random environments. For random walks on GW-trees, they were utilized already in the papers \cite{LyoPemPer96a,LyoPemPer96b} or \cite{PerZei2008}. In our definition we follow the formulation of \cite{SznZer1999}.  
There are two main changes from the classical regeneration time structure of the above mentioned papers: First, we have to allow the random walk to backtrack a distance of order $(p-p_c)^{-1}$, similar to the construction in \cite{GanMatPia2012,GanGuoNag2017}. Second, to obtain a stationary sequence even with backtracking, we have to control the environment where the walker regenerates. For the first point, we can rely on Lemma~\ref{lem:backtrackprob} to see that we have a good probability of not backtracking too far. For the second point, we want to ensure that a regeneration point has no siblings, such that the random walk can backtrack only on a branch where the backbone has no furcations.

\begin{figure}[h]
\begin{centering}
\includegraphics[width=0.9\textwidth]{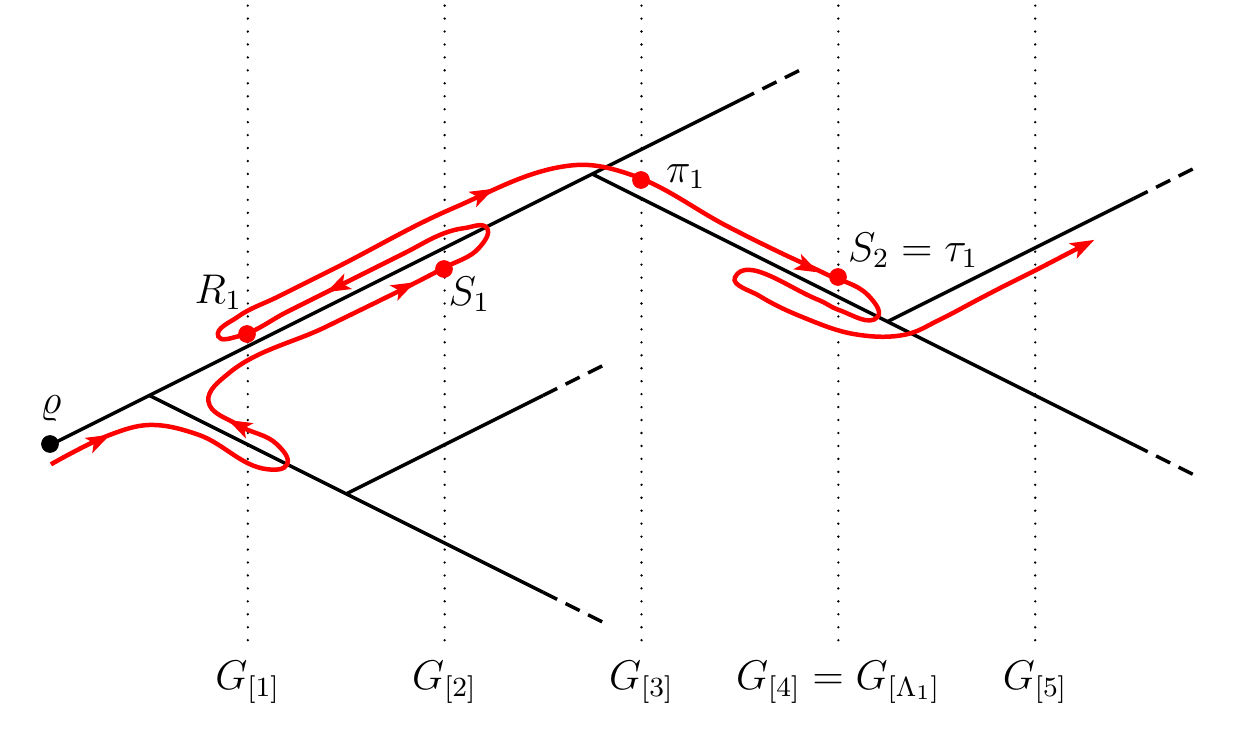}
\caption{An idealized sketch of the sample path of the random walk (in red) on a portion of the backbone tree (in black). The labels correspond to hitting times of vertices, $L$-levels are indicated by dashed lines. The hitting time $\eta_1$ of the first $L$-level is not the first potential regeneration time $S_1$, since $\sib(X_{\eta_1}) =1$. Observe that $\eta_2$ is the first potential regeneration time $S_1$ because $X_{S_1}$ has no $L$-siblings, but that the random walk then backtracks more than one $L$-level. Finally, $\eta_4=S_2$ satisfies the condition for a regeneration time and $\pi_1=\eta_3$.\label{fig:1}}
\end{centering}
\end{figure}

Intuitively, the regeneration times are constructed as follows: We wait until the random walker for the first time reaches a vertex with the potential for regeneration, namely a vertex in an $L$-level with no siblings. We call this time $S_1$ and the associated $L$-level $M_1$. If the previous $L$-level $M_1-1$ is never visited again, then we call $S_1$ the first regeneration time $\tau_1$. But suppose that the walker does revisit $L$-level $M_1-1$ again at time $R_1$. Then $S_1$ is not a regeneration time. Instead, we denote by $N_1$ the highest $L$-level that the walker visited between $S_1$ and $R_1$ on the backbone tree, and we wait until the walker first reaches an $L$-level with no siblings with generation greater than $N_1+1$ (the additional generation is required to guarantee that we see a new part of the tree, so that we can apply annealed estimates). We call this potential regeneration time $S_2$ and the corresponding $l$-level $M_2$. If $M_2-1$ is never visited again, then we set $\tau_1 = S_2$. If it is revisited again, then we repeat the above procedure. Because the walk on the tree is transient, the first regeneration time $\tau_1$ is finite almost surely. We repeat the entire procedure to construct the sequence $(\tau_k)_{k \ge 1}$ of regeneration times. See Figure~\ref{fig:1} for a sketch of this construction.

To define the above construction formally, we start with some more notation. Let $(z_n)_{n \ge 0}$ be an infinite path, and let $\theta_m$ be the time shift on the path such that $(\theta_mz)_n=z_{m+n}$ and such that for any function $f$ that takes an infinite path as its argument, $f \circ \theta_m$ denotes the same function applied to the time-shifted path. Define the \emph{backtracking time} $\eta'=\eta'((X_n)_{n\geq 1}) := \inf\{n \ge 0 \, : \, |X_n| = |X_0| - \lfloor L / (p -p_c)\rfloor\}$. We need $\eta' = \infty$ to regenerate. The regeneration structure that we rely on is then defined as follows:

\begin{defi}\label{def:regen}
Given a tree $T$ and a random walk $(X_n)_{n \ge 0}$ on $T$, we define a sequence of stopping times $S_1\leq R_1\leq S_2\leq R_2\leq \dots$ and distances $M_k$, $N_k$, beginning with 
\begin{align} \label{regdef1}
M_1& :=\min\{m\geq 1\, : \, \sib(X_{\eta_m})=0 \} ,\qquad  S_1 :=\eta_{M_1},\notag  \\
 R_1& :=S_1+\eta' \circ \theta_{S_1}, \\
 N_1& := \max \{m \ge 1 \,:\, \eta_m<R_1\} \notag ,
\end{align}
and recursively, for $k\geq 1$, 
\begin{align} \label{regdef2}
M_{k+1}& :=\min\{m\geq  N_k+2\,:\, \sib(X_{\eta_m})=0 \},\qquad  S_{k+1}:=\eta_{M_{k+1}}, \notag \\
R_{k+1}& :=S_{k+1}+\eta'\circ \theta_{S_{k+1}}, \\
N_{k+1}& := \max \{m\,:\, \eta_m<R_{k+1}\} \notag .
\end{align}
These definitions make sense until $R_k=\infty$ for some $k$ and we set 
\begin{align*}
K := \inf\{k\,:\, R_k=\infty\} , \qquad \tau_1 :=S_K.
\end{align*}
We call $\tau_1$ the \emph{first regeneration time.} 
By Lemma \ref{lem:backtrackprob}, $K<\infty$ almost surely so that $\tau_1$ is well-defined. We set $\tau_0 :=0$ and for $k\geq 2$ we define the subsequent regeneration times as
\begin{align}
\tau_{k} :=\tau_{k-1}+\tau_1 \circ \theta_{\tau_{k-1}}.
\end{align}
Finally, denote by $\pi_k := \eta(\anc(X_{\tau_k}))$ the times when the $L$-ancestors of the regeneration points 
are visited for the first time, and by $\Lambda_k$ the $L$-generation at which the $k$-th regeneration time occurs, i.e., $\Lambda_k := |X_{\tau_k}| / (\lfloor L / (p-p_c) \rfloor)$. 
\end{defi}

\begin{rk}[Decoupling property of the regenerations] \label{rem:regenerationconsequences}
We can make the following observations: 
\begin{itemize}
\item[(a)] Up to time $\pi_k$, the random walk walks on the tree up to $L$-generation $\Lambda_k-1$. From time $\tau_k$ on, the walk never visits $L$-generation $\Lambda_k-1$ again. This means that $(X_n)_{n\leq {\pi_k}}$ and $(X_n)_{n\geq {\tau_k}}$ visit disjoint parts of the tree. 
\item[(b)] The only part of the environment that is visited by both $(X_n)_{n\leq {\tau_k}}$ and $(X_n)_{n\geq {\tau_k}}$ is the tree segment consisting of the tree rooted at $X_{\pi_k}$ with all descendants of $X_{\tau_k}$ removed. Since we require $\sib(X_{\tau_k}) =0$, we know that the corresponding backbone-segment $\Tcal^{\Bb}_{X_{\pi_k}} \setminus \Tcal^{\Bb}_{X_{\tau_k}}$ is isomorphic to a line graph of length $\lfloor L /(p-p_c)\rfloor-1$. 
\item[(c)] The definition above allows for the random walk to move distance at most $\lfloor L / (p-p_c) \rfloor-1$ towards the root after a regeneration time. Since the backbone structure in this tree segment is fixed, and because the inter-regeneration distances $|X_{\tau_{k}}|-|X_{\tau_{k-1}}|$ are independent of the traps, we may conclude that the inter-regeneration distances are independent. With the exception of $k=1$, they are identically distributed as well. (This observation is formalized in Lemma~\ref{lem:stationarity} below.)
\item[(d)] The inter-regeneration times $\tau_{k}-\tau_{k-1}$ are not independent, since two subsequent time intervals both depend on the traps in the tree segment between $X_{\pi_{k-1}}$ and $X_{\tau_{k}}$. The inter-regeneration times are, however, stationary and $1$-dependent. (This is formalized in Lemma~\ref{lem:stationarity2} below.)
\end{itemize}
\end{rk}

To state the above observations generally and unambiguously, we introduce the $\sigma$-fields
\begin{align*}
\mathcal{G}_k := \sigma \big( \mathcal{T}\setminus \mathcal{T}_{X_{\pi_{k}}}, (X_n)_{0\leq n< \pi_{k}}, \pi_{k}  \big)
\end{align*}
and the time shift $\bar\theta_m$ defined for a set
\begin{align*}
B = \big\{  \mathcal{T} \in B_1, (X_n)_{n\geq 0}\in B_2 \big\}
\end{align*}
by
\begin{align*}
B\circ \bar\theta_k  =  \big\{  \mathcal{T}_{X_{\pi_{k}}} \in B_1, (X_n-X_{\tau_k})_{n\geq  \tau_{k}}\in B_2 \big\} .
\end{align*}
We say that a sequence $(Y_n)_{n\geq 0}$ is $m$-dependent, if $(Y_n)_{n\leq k}$ and $(Y_n)_{n\geq k+\ell}$ are independent whenever $\ell >m$.  
The following two lemmas are standard, so we omit their proofs (see e.g.\ \cite[Lemma~17 and Proposition~18]{Guo2016}):  

\begin{lem}[Stationarity of the tree and walk] \label{lem:stationarity}
For any measurable set $B = \big\{  \mathcal{T} \in B_1, (X_n)_{n\geq 0}\in B_2 \big\}$
and $k\geq 1$,
\begin{align*}
\mathbb{P}_p(B\circ\bar\theta_k \mid \mathcal{G}_k) = \mathbb{P}_p^v(B \mid G_{[1]}(\mathcal T)=\{v\} , \eta_0=\infty ) .  
\end{align*}
\end{lem}

\begin{lem}[Stationarity and $1$-dependence of the regeneration times] \label{lem:stationarity2}
Under $\mathbb{P}_p$, the sequence
\begin{align*}
\big( \mathcal{T}_{X_{\pi_k}}\setminus \mathcal{T}_{X_{\tau_{k+1}}}, (X_n-X_{\tau_k})_{\tau_k\leq n< \tau_{k+1}}, \tau_{k+1}-\tau_k \big)_{k\geq 1}
\end{align*}
is stationary and 1-dependent. Furthermore, the marginal distribution of this sequence is given by
\begin{multline*}
\mathbb{P}_p\big( \mathcal{T}_{X_{\pi_k}}\setminus \mathcal{T}_{X_{\tau_{k+1}}}\in B_1, (X_n-X_{\tau_k})_{\tau_k\leq n< \tau_{k+1}}\in B_2, \tau_{k+1}-\tau_k\in B_3 \big) \\
= \mathbb{P}_p^v\big( \mathcal{T}\setminus \mathcal{T}_{X_{\tau_{1}}}\in B_1, (X_n)_{0\leq n< \tau_{1}}\in B_2, \tau_{1}\in B_3 \mid  G_{[1]}(\mathcal T)=\{v\} , \eta_0=\infty \big) . 
\end{multline*}
\end{lem}

The moment bounds for the regeneration distances and the regeneration times in the two following lemmas are crucial ingredients for the proof of the main result:

\begin{lem}[Moment bounds on regeneration distances]\label{lem:moments}
For any $a\geq 1$, there exists a finite constant $C_a$ such that 
\begin{align*}
(p-p_c)^a \mathbb{E}_p\left[|X_{\tau_1}|^a\right] \leq C_a, \qquad (p-p_c)^a\mathbb{E}_p\left[(|X_{\tau_{2}}|-|X_{\tau_1}|)^a \right] \leq C_a,
\end{align*}
for any $p>p_c$. 
\end{lem}
We prove this lemma in Section~\ref{sec:regdist} below. Note that we also have the trivial lower bound $|X_{\tau_1}|\geq c(p-p_c)^{-1}$. 

\begin{lem}[Moment bounds on regeneration times]\label{lem:moments2}
There exists a constant $C$ such that 
\begin{align*}
(p-p_c)^6 \mathbb{E}_p\left[({\tau_1})^2\right] \leq C, \qquad (p-p_c)^6\mathbb{E}_p \left[({\tau_{2}}-{\tau_1})^2 \right] \leq C,
\end{align*}
for any $p>p_c$. Furthermore, there exists a constant $c_0>0$, such that
\begin{align*}
(p-p_c)^3 \mathbb{E}_p[({\tau_{2}}-{\tau_1})] \geq c_0.
\end{align*}
\end{lem}
We prove this lemma in Section~\ref{sec:regtime} below.
 
\begin{rk}[Robust bounds on the effective speed near criticality]\label{rem:speedbounds} 
Given the moment bounds of Lemmas \ref{lem:moments} and \ref{lem:moments2}, we may apply the Law of Large Numbers to obtain an expression for the effective speed in terms of regeneration times, as
\begin{align} \label{speedformula}
v(p) = \lim_{n\to\infty} \frac{|X_n|}{n} = \lim_{n\to \infty} \frac{|X_{\tau_n}|}{\tau_n} = \frac{\mathbb{E}_p[|X_{\tau_2}|-|X_{\tau_1}|]}{\mathbb{E}_p[\tau_2-\tau_1]} . 
\end{align}
In particular, the $p$-independence of the moment bounds implies  
\begin{align} \label{speedbounds}
0< \liminf_{p_n \searrow p_c} (p_n-p_c)^{-2} v(p_n) \leq  \limsup_{p_n \searrow p_c}(p_n-p_c)^{-2} v(p_n) <\infty .
\end{align}
This gives the order of the effective speed $v(p)$ close to criticality. The result of Theorem~\ref{thm:speed} is stronger and shows that $\liminf$ and $\limsup$ in \eqref{speedbounds} agree, but the proof of Theorem~\ref{thm:speed} in Section~\ref{sec:speed} below relies on an explicit formula for the effective speed of SRW on GW-trees obtained in \cite{LyoPemPer96a}, while the above arguments are quite robust against changes to tree or the behavior of the random walk.  
\end{rk}

\section{The scaling limit: proof of Theorem \ref{thm:main}}\label{sec:mainpf}

In this section we prove Theorem~\ref{thm:main}, subject to the proofs of Theorem~\ref{thm:speed}, the lemmas in Sections~\ref{sec:prelim0}, \ref{prelim}, and \ref{sec:regen}, and subject to the proofs of two further lemmas, Lemmas \ref{lem:nzlemma} and \ref{lem:BBTbounds}, which are stated below as we need them.
 
The proof of Theorem~\ref{thm:main} will go in four steps. In the first two, we consider the scaling limit of a simple random walk whose jumps have the same size distribution as the increments of the random walk on the BRW at regeneration times, but with jumps at a fixed rate. We show that under a rescaling equivalent to the one of Theorem~\ref{thm:main}, this process converges to a Brownian motion with diffusion $(\kappa \Sigma)^{1/2}$ as desired. Then, we will apply a time change to have the jumps occur at random times with the same distribution as the regeneration times and show that the difference with the process of the first step vanishes in the limit. In the final step, we show that the trajectories of the random walk during excurions between regeneration times also vanishes in the limit, thus yielding the scaling limit for the process that we are after.

We assume in this section some familiarity with the theory of convergence for Markov processes. We refer the reader who is insufficiently knowledgable about this topic to \cite{Bill99,EthKur86}, where all convergence-related topics of this section that are not explicitly cited are defined and discussed.

To start, let us recall that in the setting of Theorem~\ref{thm:main} we consider $p = p_n$ with $p_n \searrow p_c$ as $n \to \infty$. Therefore, the distribution of the random walk changes with $n$. We will thus indicate this dependency on $n$ clearly in this section. 
\medskip

\textbf{First step: i.i.d. increments between regenerations.}  
For the first step of the proof, we define the process
\begin{align*}
\mathcal{W}^{(1)}_{t,n} := \frac{1}{\sqrt{n}}\sum_{k=1}^{\lfloor nt\rfloor} W_k
\end{align*}
as a random element of the Skorokhod space $\mathcal{D}_{\mathbb{R}^d}[0,T]$. 
Without loss of generality, we take $T\in \mathbb{N}$ and then, rescaling $n$ linearly and using the scale invariance of the limiting Brownian motion, may restrict the proof further to the case $T=1$. 
The increments of $\mathcal{W}_{t,n}^{(1)}$ are given by 
	\[
	W_1 := \frac{\varphi(X_{\tau_1})}{\sqrt{ (p_n-p_c)^2\E_{p_n}[\tau_1]}}, \quad \text{ and } \quad W_k := \frac{\varphi(X_{\tau_{k}})-\varphi(X_{\tau_{k-1}})}{\sqrt{(p_n-p_c)^{2}\mathbb{E}_{p_n}[\tau_2-\tau_1]}} \quad \text{ for $k \ge 2$.}
	\]

By Definition~\ref{def:regen} and Lemma~\ref{lem:stationarity}, the $W_k$ are independent and $W_2,W_3,\dots $ are identically distributed. Conditioned on $(X_n)_{n \ge 0}$, $\varphi(X_{\tau_{k}})-\varphi(X_{\tau_{k-1}})$ is the increment of $(Y_n)_{n \ge 0}$ (the random walk on $\Zd$ with step distribution $D$) after $|X_{\tau_{k}}|-|X_{\tau_{k-1}}|$ steps. Thus, by the assumption in Theorem~\ref{thm:main} that $\sum_{x \in \Zd} x D(x) =0$, we get $\mathbb{E}[W_k] = 0$ and by Lemma~\ref{lem:moments} and Lemma~\ref{lem:moments2} we get
\begin{align}\label{e:cond1}
\mathbb{E}_{p_n}\left[\left\| W_k \right\|^2 \right] 
= \frac{\mathbb{E}_{p_n}\left[|X_{\tau_{k}}|-|X_{\tau_{k-1}}| \right] \mathbb{E}_{p_n}\left[\| Y_1 \|^2 \right]}{(p_n-p_c)^{2} \mathbb{E}_{p_n}[\tau_2-\tau_1]}  \leq C 
\end{align}
for $k\geq 2$ and for the first term
\begin{align}\label{e:cond1b}
\mathbb{E}_{p_n}\left[\left\| W_1 \right\|^2 \right] 
= \frac{\mathbb{E}_{p_n}\left[|X_{\tau_{1}}|\right] \mathbb{E}_{p_n}\left[\| Y_1 \|^2 \right]}{(p_n-p_c)^{2} \mathbb{E}_{p_n}[\tau_1]}  \leq C .
\end{align}
Writing $v^{\mathsf{T}}$ for the transpose of $v$, we also have the convergence 
\begin{align}\label{e:cond2}
\lim_{n\to \infty } \frac{1}{n} \mathbb{E}_{p_n}\left[ \left(\sum_{k=1}^n W_k\right) \left(\sum_{k=1}^n W_k\right)^{\mathsf{T}} \right] & = \lim_{n\to \infty}  \frac{\mathbb{E}_{p_n}\left[|X_{\tau_{n}}|\right]}{n (p_n-p_c)^2 \mathbb{E}_{p_n}[\tau_2-\tau_1]} \mathbb{E}_{p_n}\left[Y_1Y_1^{\mathsf{T}} \right] \notag \\
& = \lim_{n\to \infty}  \frac{ \mathbb{E}_{p_n}\left[|X_{\tau_{2}}|-|X_{\tau_{1}}|\right]}{(p_n-p_c)^2 \mathbb{E}_{p_n}[\tau_2-\tau_1]} \mathbb{E}_{p_n}\left[Y_1Y_1^{\mathsf{T}} \right] \notag \\
& = \lim_{n\to \infty} \frac{v(p_n)}{(p_n-p_c)^2}  \mathbb{E}_{p_n}\left[Y_1Y_1^{\mathsf{T}} \right] \\
& = \kappa \Sigma , \notag
\end{align}
where in the last step we have used Theorem~\ref{thm:speed} and the assumption in Theorem~\ref{thm:main} that $\sum_{x \in \Zd} \|x\|^2 D(x) < \infty$, so that $\Sigma$ exists. Because $\E[W_k]=0$ and \eqref{e:cond1} and \eqref{e:cond2} hold, we may now apply the Invariance Principle for triangular arrays (see \cite{Bill99}, p. 147, for a one-dimensional version, which applies thanks to the Cram\'{e}r-Wold technique), and conclude the convergence 
\begin{align} \label{processconv1}
\left(\mathcal{W}^{(1)}_{t,n}\right)_{t\in[0,1]} \xrightarrow[n\to \infty]{\mathrm{d}} \left( (\kappa \Sigma)^{1/2}  B_t \right)_{t\in[0,1]} ,
\end{align}
under $\mathbb{P}_{p_n}$, with $(B_t)_{t \ge 0}$ a standard Brownian motion on $\Rd$. This concludes the first step.
\medskip

\textbf{Second step: The correct (random) number of jumps.} 
The process  $\mathcal{W}^{(1)}_{n} := (\mathcal{W}^{(1)}_{t,n})_{t \in [0,1]}$ is a piecewise constant function in $t$, jumping exactly $n$ times in $[0,1]$. The second step is to construct a process that will also jump at regular intervals, but has a random number of jumps, determined by the random number of regeneration times we see in $n$ steps. Given $n \ge 0$, let $k_n$ be the integer satisfying 
\begin{align} \label{randomnumberofjumps}
\tau_{k_n}\leq n < \tau_{k_n+1} , 
\end{align}
where we recall that $\tau_0=0$. Set $\nu_n :=(p_n-p_c)^{-3}k_n$ and define
\begin{align}\label{def2ndprocess}
\mathcal{W}^{(2)}_{t,n} := \mathcal{W}^{(1)}_{t,\nu_n} .
\end{align} 
This process jumps $\nu_n$ times and we next show that this number of jumps is asymptotically equal to 
\begin{align}\label{randomnumberofjumps2}
a_n :=n(p_n-p_c)^{-3}\mathbb{E}_{p_n}[\tau_2-\tau_1]^{-1} . 
\end{align}
Note that by Lemma \ref{lem:moments2}, $a_n\to \infty$. Now, by definition of $k_n$,
\begin{align}\label{randomnumberofjumps3}
\frac{a_n}{\nu_n} = \frac{n}{k_n\mathbb{E}_{p_n}[\tau_2-\tau_1]}\leq \frac{\tau_{k_n+1}}{k_n\mathbb{E}_{p_n}[\tau_2-\tau_1]} 
= \frac{1}{k_n \mathbb{E}_{p_n}[\tau_2-\tau_1]} \sum_{m=1}^{k_n+1} (\tau_m-\tau_{m-1}) .
\end{align}
Using the 1-dependence of $(\tau_{m}-\tau_{m-1})_m$ and the moment bound in Lemma \ref{lem:moments2}, Chebyshev's inequality shows that the variance of the right hand side of \eqref{randomnumberofjumps3} converges to zero as $n\to \infty$. Thus, the right hand side of \eqref{randomnumberofjumps3} converges to 1 in distribution under $\mathbb{P}_{p_n}$. Bounding $n\geq \tau_{k_n}$, the same argument shows that $ a_n/\nu_n\geq 1$ dominates a sequence of random variables converging to 1, which implies that under $\mathbb{P}_{p_n}$,
\begin{align} \label{randomnumberofjumps4}
\frac{a_n}{\nu_n} \xrightarrow[n\to \infty]{\mathrm{d}} 1 .
\end{align}
We may then apply \cite[Theorem~14.4]{Bill99} to conclude that $\mathcal{W}^{(2)}_{t,n}$ has the same limit as $\mathcal{W}^{(1)}_{t,n}$, that is, \eqref{processconv1} with $\mathcal{W}^{(1)}_{t,n}$ replaced by $\mathcal{W}^{(2)}_{t,n}$.
This concludes the second step.
\medskip

\textbf{Third step: The correct (random) jump times.} 
The process  $\mathcal{W}^{(2)}_{n} := (\mathcal{W}^{(2)}_{t,n})_{t \in [0,1]}$ still jumps with equal intervals, at times $t=\tfrac{1}{\nu_n},\tfrac{2}{\nu_n},\dots ,1$. 
The third step of the proof is to consider instead the process $\mathcal{W}^{(3)}_n := (\mathcal{W}^{(3)}_{t,n})_{t \in [0,1]}$, which has the same increments as $\mathcal{W}^{(2)}_{n}$, but jumps at times $t= \tau_1 /\tau_{\nu_n},$ $\tau_2 /\tau_{\nu_n}, \dots ,1$. That is, 
\begin{align*}
\mathcal{W}^{(3)}_{t,n} := \mathcal{W}^{(2)}_{k/\nu_n,n} = \frac{\varphi(X_{\tau_{k}})}{\sqrt{\nu_n(p_n-p_c)^{2} \mathbb{E}_{p_n}[\tau_2-\tau_1]}}  \quad \text{ for } \quad \tau_{k} \leq \lfloor t \tau_{\nu_n} \rfloor < \tau_{k+1},
\end{align*}
so $\mathcal{W}^{(3)}_{n}$ is a random time change of $\mathcal{W}^{(2)}_{n}$. We can bound their distance in the Skorokhod metric $d_S$ as
\begin{align} \label{skorokhodbound}
		d_S( \mathcal{W}^{(3)}_n,\mathcal{W}^{(2)}_n) & \leq \max_{k=1,\dots ,\nu_n} \left| \frac{\tau_k}{\tau_{\nu_n}}-\frac{k}{\nu_n} \right| \notag \\
		& \leq \max_{k=1,\dots ,\nu_n} \left| \frac{\tau_k-\mathbb{E}_{p_n}[\tau_k-\tau_1]}{\tau_{\nu_n}}\right| +  \max_{k=1,\dots ,\nu_n} \left|\frac{\mathbb{E}_{p_n}[\tau_k-\tau_1]}{\tau_{\nu_n}}-\frac{k}{\nu_n} \right|  .
\end{align}
By Lemma~\ref{lem:stationarity2} we have $\mathbb{E}_{p_n}[\tau_k-\tau_1] = k \mathbb{E}_{p_n}[\tau_2-\tau_1]$, so we may write the second term as
\begin{align*}
 \max_{k=1,\dots ,\nu_n} \left|\frac{\mathbb{E}_{p_n}[\tau_k-\tau_1]}{\tau_{\nu_n}}-\frac{k}{\nu_n} \right|  
 =  \left|  \frac{\nu_n\mathbb{E}_{p_n}[\tau_2-\tau_1]}{\tau_{\nu_n}}  -1 \right|.
\end{align*} 
Using again the moment bounds of Lemma~\ref{lem:moments2},
\begin{align}\label{skorokhodbound2}
\frac{\tau_{\nu_n}}{\nu_n \mathbb{E}_{p_n}[\tau_2-\tau_1]} = \frac{\tau_1 }{\nu_n\mathbb{E}_{p_n}[\tau_2-\tau_1]} + \frac{1}{\nu_n} \sum_{k=2}^{\nu_n} \frac{\tau_k-\tau_{k-1}}{\mathbb{E}_{p_n}[\tau_2-\tau_1]} \xrightarrow[n\to \infty]{\mathrm{d}} 1,
\end{align}
in distribution under $\mathbb{P}_{p_n}$. This shows that the second term in \eqref{skorokhodbound} vanishes as $n\to \infty$.

We write the first term in \eqref{skorokhodbound} as
\begin{align*}
\max_{k=1,\dots ,\nu_n} \left| \frac{\tau_1}{\nu_n\mathbb{E}_{p_n}[\tau_2-\tau_1]} + \frac{(\tau_k-\tau_1)-\mathbb{E}_{p_n}[\tau_k-\tau_1]}{\nu_n\mathbb{E}_{p_n}[\tau_2-\tau_1]}\right| \left| \frac{\tau_{\nu_n}}{\nu_n\mathbb{E}_{p_n}[\tau_2-\tau_1]} \right|^{-1}.
\end{align*}
The second factor converges in distribution to $1$ by \eqref{skorokhodbound2}. That the first factor also converges to $1$ follows from Lemma~\ref{lem:stationarity2} combined with an Invariance Principle for triangular arrays of 1-dependent random variables due to Chen and Romano \cite[Theorem~2.1]{CheRom99}. This implies that \eqref{skorokhodbound} vanishes in probability and thus \eqref{processconv1} also holds with $\mathcal{W}^{(1)}_{n}$ replaced by $\mathcal{W}^{(3)}_{n}$. 
This concludes the third step.
\medskip

\textbf{Fourth step: Adding the path between regenerations.} 
The fourth and last main step of the proof is now to transfer the convergence to 
\begin{align*}
\mathcal{W}^{(4)}_{t,n} := \frac{\varphi(X_{\lfloor t\tau_{\nu_n}\rfloor})}{\sqrt{\nu_n (p_n -p_c)^{2}\mathbb{E}_{p_n}[\tau_2-\tau_1]}},  
\end{align*}
i.e., to the full (rescaled) process of random walk on a BRW, where we now also incorporate the fluctuations of $(X_n)_{n \ge 0}$ between regeneration times. 
Note that by definition of $a_n$ and $\nu_n$ and by \eqref{randomnumberofjumps4}, we have
\begin{align*}
\frac{1}{\sqrt{\nu_n (p_n-p_c)^{2}\mathbb{E}_{p_n}[\tau_2-\tau_1]}} = \sqrt{\frac{a_n(p_n-p_c)}{\nu_n n}}\leq C \sqrt{\frac{p_n-p_c}{n}},
\end{align*}
so bounding
\begin{align} \label{skorokhodbound3}
d_S( \mathcal{W}^{(4)}_n,\mathcal{W}^{(3)}_n) \leq C \sqrt{\frac{p_n-p_c}{n}} \max_{k=1,\dots ,\nu_n}\ \max_{\tau_{k-1}\leq i<\tau_k} \| \varphi(X_{i}) - \varphi(X_{\tau_{k-1}})\|, 
\end{align}
it suffices to show that the right-hand side vanishes in probability. Again by the moment bounds of Lemma \ref{lem:moments2}, $\nu_n\leq Cn$ with probability converging to 1, so that we may replace the remaining $\nu_n$ in \eqref{skorokhodbound3} by $Cn$.  
We need a lemma that controls the maximal displacement of the embedded traps. As we show in Section~\ref{sec:maxdisp}, this is a direct consequence of a precise estimate by Neuman and Zheng \cite{NeuZhe17} for the maximal displacement of subcritical branching random walks. 

Recall that in the decomposition of $\Tcal$ into the backbone tree and the traps given in Section~\ref{sec:prelim0}, for any $v \in \Tcal^{\Bb}$ we have that the law of $\Tcal^{\trap}_v$ only depends on $v$ through $\deg_{\Tcal}(v)$, and that we have assumed in Theorem~\ref{thm:main} that $\deg_{\Tcal}(v) \le \Delta+1$ almost surely.

\begin{lem}[Maximal displacement inside traps]
\label{lem:nzlemma}
Under the assumptions of Theorem \ref{thm:main}, there exists a $\gamma>0$ and $C<\infty$, both independent of $p$, such that
\begin{align} \label{trapsizeassumption}
\max_{1 \le \delta \le \Delta -1} \mathbb{E}_p\left[\left. \exp \left( \gamma \sqrt{p-p_c} \sup_{w\in  \mathcal{T}^{\trap}_v} \|\varphi(w)-\varphi(v)\| \right) \right| \deg_{\Tcal^\trap_v}(v) = \delta \right] \leq C.
\end{align}
\end{lem}
We prove this lemma in Section~\ref{sec:maxdisp} below.

We also need to bound the number of visited backbone-tree vertices. 
For this we write
\begin{align*}
\operatorname{BBT}_k := \{X_{\tau_{k-1}},X_{\tau_{k-1}+1},\dots ,X_{\tau_{k}} \} \cap \mathcal{T}^{\Bb },
\end{align*}
i.e., $\operatorname{BBT}_k$ is the trace of $(X_n)_{\tau_{k-1} \le n \le \tau_k}$ restricted to the backbone tree.

\begin{lem}[Size of and maximal distance from $\operatorname{BBT}_k$] 
\label{lem:BBTbounds}
Under the assumptions of Theorem \ref{thm:main}, for any $a\geq 1$ there exists a $C_a$ such that for any $k\geq 1$ 
\begin{align}
\mathbb{E}_p[|\operatorname{BBT}_k|^a] \leq C_a(p-p_c)^{-a} , \label{backbonetracebound1}
\end{align}
and
\begin{align}
\mathbb{E}_p[\max_{v\in \operatorname{BBT}_k} \|\varphi(v)-\varphi(X_{\tau_{k-1}})\|^a] \leq C_a(p-p_c)^{-a/2} . \label{backbonetracebound2} 
\end{align}
\end{lem}
We prove this lemma in Section~\ref{sec:BBT} below.

Using the estimates of Lemmas \ref{lem:nzlemma} and \ref{lem:BBTbounds}, we can show that the right-hand side of \eqref{skorokhodbound3} vanishes. Indeed,
\begin{multline}\label{fluctuationbound}
  \sqrt{\frac{p_n-p_c}{n}} \max_{\tau_{k-1}\leq i<\tau_k} \| \varphi(X_{i}) - \varphi(X_{\tau_{k-1}})\|  \\
  \leq  \sqrt{\frac{p_n-p_c}{n}}  \left( \max_{v\in \operatorname{BBT}_k}  \|\varphi(v)-\varphi(X_{\tau_{k-1}})\| 
    + \max_{v\in \operatorname{BBT}_k }  \sup_{w\in  \mathcal{T}^{\trap}_v} \|\varphi(w)-\varphi(v)\|  \right).
\end{multline}
By \eqref{backbonetracebound2} with $a=3$, we have for any $\varepsilon>0$,
\begin{align} \label{fluctuationbound2}
\mathbb{P}_{p_n} \left( \sqrt{\frac{p_n-p_c}{n}}  \max_{v\in \operatorname{BBT}_k}  \|\varphi(v)-\varphi(X_{\tau_{k-1}})\| >\varepsilon \right) \leq \frac{1}{\varepsilon^3 n^{3/2}} C.
\end{align}
Using further that by the decomposition \eqref{treedecomposition}, the GW-trees of the traps are independent of the backbone tree,
\begin{align}
&  \mathbb{P}_{p_n} \left( \sqrt{\frac{p_n-p_c}{n}}  \max_{v\in \operatorname{BBT}_k }  \sup_{w\in  \mathcal{T}^{\trap}_v} \|\varphi(w)-\varphi(v)\|  >\varepsilon \right)\notag \\  
& \quad \leq \mathbb{E}_{p_n} [|\operatorname{BBT}_k| ]\max_{1 \le \delta \le \Delta-1} \mathbb{P}_{p_n} \left( \left. \sqrt{\frac{p_n-p_c}{n}}  \sup_{w\in  \mathcal{T}^{\trap}_v} \|\varphi(w)-\varphi(v)\|  >\varepsilon \right| \deg_{\Tcal^\trap_v} (v) = \delta\right) \notag \\
& \quad \leq  \mathrm{e}^{-\gamma\varepsilon\sqrt{n}} \, \mathbb{E}_{p_n} [|\operatorname{BBT}_k| ] \notag \\
& \qquad \times \max_{1 \le \delta \le \Delta-1} \mathbb{E}_{p_n}\left[ \left.\exp\left( \gamma \sqrt{p_n-p_c} \sup_{w\in  \mathcal{T}^{\trap}_v} \|\varphi(w)-\varphi(v)\| \right) \right| \deg_{\Tcal^\trap_v} (v) = \delta \right] \notag \\
& \quad \leq C \left( (p_n-p_c) \mathrm{e}^{\gamma\varepsilon\sqrt{n}} \right)^{-1} , \label{fluctuationbound3}
\end{align}
where the last inequality follows from \eqref{trapsizeassumption} and \eqref{backbonetracebound1}. 
From \eqref{fluctuationbound}, \eqref{fluctuationbound2}, and \eqref{fluctuationbound3} we may conclude
\begin{multline} \label{fluctuationbound4}
\mathbb{P}_{p_n} \left(  \sqrt{\frac{p_n-p_c}{n}} \max_{k=1,\dots ,Cn}\ \max_{\tau_{k-1}\leq i<\tau_k} \| \varphi(X_{i}) - \varphi(X_{k-1})\| > \varepsilon\right)\\
 \leq \frac{C}{\varepsilon^3 n^{1/2}}  + C \left( (p_n-p_c) \mathrm{e}^{\gamma\varepsilon\sqrt{n}} \right)^{-1} .
\end{multline}
Recall that in Theorem~\ref{thm:main} we assumed that for any $\delta>0$, $\mathrm{e}^{\delta \sqrt{n}}(p_n-p_c)\to \infty$. Under this assumption the right-hand side converges to zero, which implies that the right-hand side of \eqref{skorokhodbound3} vanishes in probability. 
Therefore, the convergence \eqref{processconv1} also holds with $\mathcal{W}^{(1)}_n$ replaced by $\mathcal{W}^{(4)}_n$, that is, 
\begin{align} \label{processconv2}
\left(\frac{ \varphi(X_{\lfloor t\tau_{\nu_n}\rfloor})}{\sqrt{\nu_n\mathbb{E}_{p_n}[\tau_2-\tau_1](p_n-p_c)^{2}}} \right)_{t\in[0,1]} \xrightarrow[n\to \infty]{\mathrm{d}} \left( (\kappa \Sigma)^{1/2}  B_t \right)_{t\in[0,1]}.
\end{align}
This completes the fourth step.
\medskip

All that remains to finish the proof of Theorem~\ref{thm:main} is to apply the time change $t \mapsto t n(p_n-p_c)^{-3}\tau_{\nu_n}^{-1}$, which, by the scale invariance of Brownian motion, yields   
\begin{align} \label{processconv2}
\left(\sqrt{\frac{(p_n-p_c)^{3}\tau_{\nu_n}}{n\nu_n\mathbb{E}_{p_n}[\tau_2-\tau_1](p_n-p_c)^{2}}}   \varphi(X_{\lfloor tn(p_n-p_c)^{-3}\rfloor}) \right)_{t\in[0,1]} \xrightarrow[n\to \infty]{\mathrm{d}} \left( (\kappa \Sigma)^{1/2}  B_t \right)_{t\in[0,1]} .
\end{align}
The Law of Large Numbers for $\tau_n$ then allows us to replace the prefactor by $\sqrt{(p-p_c)/n}$, which completes the proof.\qed


\section{Proofs of the escape time estimates}
\label{sec:aprioriproofs}

In this section we prove the lemmas of Section~\ref{prelim}. We assume in this section some familiarity with the theory of reversible Markov chains and electrical networks. We refer the reader who is insufficiently knowledgable about this topic to \cite{AldFil02,LyoPer16}, where all reversiblility-related topics of this section that are not explicitly cited are defined and discussed.

\subsection{Proof of Lemma \ref{lem:backtrackprob}}\label{sec:backtrack}

Let $T$ be an infinite tree and $v\in G_{[1]}(T)$. We have
\begin{align*}
P_T^v(\eta_0 <\infty) = \lim_{n\to \infty} P_T^v(\eta_0 <\eta_n) .
\end{align*}
Denote by $\mathcal{C}_G(x,A)$ the effective conductance in a graph $G = (V,E)$ with unit edge weights between a vertex $x\in V$ and a subset $A\subset V$, and write $\mathcal{R}_G(x,A) := \mathcal{C}_G(x,A)^{-1}$ for the associated effective resistance. For two disjoint sets $A,B$ and $x\notin A,B$ we have the bound
\begin{align} \label{effresistanceineq}
 P_G^x(\eta(A) <\eta(B)) \leq \frac{\mathcal{C}_G(x,A)}{\mathcal{C}_G(x,A\cup B)}\leq  \frac{\mathcal{C}_G(x,A)}{\mathcal{C}_G(x,B)} ,
\end{align}
see \cite[Exercise 2.34]{LyoPer16} or \cite[Fact 2]{BerGanPer03}. In fact, the arguments in the latter reference even show that if it is impossible that the random walk hits both $A$ and $B$ during the same excursion starting from $x$, then 
\begin{align} \label{effresistanceineq0}
 P_G^x(\eta(A) <\eta(B)) = \frac{\mathcal{C}_G(x,A)}{\mathcal{C}_G(x,A\cup B)} = \frac{\mathcal{C}_G(x,A)}{\mathcal{C}_G(x,A)+\mathcal{C}_G(x,B)} .
\end{align}
In our setting, \eqref{effresistanceineq} gives
\begin{align} \label{effresistanceineq01}
 P_T^v(\eta_0 <\eta_n) \leq \frac{\mathcal{C}_T(v,\varrho)}{\mathcal{C}_T(v,G_{[n]}(T))}.
\end{align}
The Series Law for conductances implies that for any $v \in G_{[1]}(T)$ we have $\mathcal{C}_T(v,\varrho)=\lfloor L/(p-p_c)\rfloor^{-1}$. Rayleigh's Monotonicity Principle, moreover, implies that for any $v \in G_{[1]}(T)$ we have
\begin{align*}
\mathcal{C}_T(v,G_{[n]}(T))\ge \mathcal{C}_{T_v}(v,G_{[n-1]}(T_v)) .
\end{align*} 
Observe that $\mathcal{C}_{T_v}(v,G_{[n-1]}(T_v))$ (and any bound we formulate from here on) depends only on the subtree rooted at $v$, and is thus independent of the tree segment $\Tcal_{[0,1]}$ under $\bar{\mathbf{P}}_p$. This implies
\begin{align}\label{effresistanceineq1}
\mathbb{P}_p^v(\eta_0 <\infty \mid \mathcal{T}_{[0,1]}=T_{[0,1]}) \leq \lim_{n\to \infty}  \frac{p-p_c}{L} \bar{\mathbf{E}}_p\left[ \mathcal{R}_{\Tcal_v}(v,G_{[n-1]}(\Tcal_v)) \right].
\end{align}

Now given an infinite tree $T$ and $v \in G_{[1]}(T)$, we prune $T_v$ as follows: first, remove any vertex $w \in T_v$ that does not have an infinite line of descent to obtain $T^{\Bb}_v$, and second, at every vertex $w\in T^{\Bb}_v$ with more than two children in $T^{\Bb}_v$, keep the first two in Ulam-Harris ordering and delete all other children (and their subtrees). Call the resulting tree $\tilde T_v$. Since we only removed edges, Rayleigh's Monotonicity Principle once more implies
\begin{align}\label{effresistanceineq2}
\mathcal{R}_{T_v}(v,G_{[n]}(T))\le \mathcal{R}_{\tilde T_v}(v,G_{[n]}(\tilde T_v)) .
\end{align}
Now let $\tilde G_{[n]}(\tilde T_v)$ denote the set of vertices $w \in \tilde T_v$ such that $\operatorname{deg}_{\tilde T_v}(w) = 3$ and such that the unique path $v =v_0, v_1, \dots, v_m =w$ in $\tilde T_v$ satisfies $|\{v_i \,:\, \deg_{\tilde{T}_v}(v_i) =3\}|=n$ (i.e., there are exactly $n$ degree-three vertices between $v$ and $w$ in $\tilde T_v$).
For $\bar{\mathbf{P}}$-almost all trees, $\tilde G_{[n]}(\tilde \Tcal_v)$ is nonempty for all $n$. Moreover, for any $T$, 
\begin{align} \label{effresistanceineq3}
\lim_{n\to \infty} \mathcal{R}_{\tilde T_v}(v,G_{[n]}(\tilde T_v)) = \lim_{n\to \infty} \mathcal{R}_{\tilde T_v}(v,\tilde G_{[n]}(\tilde T_v)) ,
\end{align}
because this limit does not depend on the choice of the sequence of exhausting subgraps of $\tilde T_v$ (see \cite[Exercise 2.4]{LyoPer16}). 

Under $\bar{\mathbf{P}}_p$, the subtree of $\tilde \Tcal_v$ of all paths connecting $v$ to $\tilde G_{[n]}(\tilde \Tcal_v)$ now has a particularly easy structure: vertices with two children are connected by line segments (a sequence of vertices with one child) up to the $n$-th ``generation'' of branching points. By the Series Law, the effective resistances of these line segments are just their respective lengths, so they are independent Geo$(1-\hat p_1)$ random variables (which take values on $\{1,2,\dots\}$), where $\hat p_1$ is the probability that the root of a Galton-Watson tree with generating function $\hat f(s)$ as described in \eqref{treedecomposition} has exactly one child. By Lemma~\ref{lem:gf}, $\hat p_1 = \hat f'(0) = f'(q_p) = 1-c_1(p-p_c)(1+o(1))$, so the mean length of a path segment between two branching points in $\tilde \Tcal_v$ is at most $\lambda (p-p_c)^{-1}$ for some $\lambda >0$.

Suppose $\tilde G_{[1]}(\tilde T_v)= \{w\}$ and call $\tilde T_w^{(1)}, \tilde T_w^{(2)}$ the two subtrees of $\tilde T_v$ rooted at $w$. An application of the Series Law and the Parallel Law gives rise to the recursion
\begin{align} \label{effresistancerec}
\mathcal{R}_{\tilde T_v}(v,\tilde G_{[n]}(\tilde T_v)) = d_{\tilde T_v} (v,w) + \left(\frac{1}{\mathcal{R}_{\tilde T_w^{(1)}}(w,\tilde G_{[n-1]}(\tilde T_w^{(1)}))}+\frac{1}{\mathcal{R}_{\tilde T_w^{(2)}}(w,\tilde G_{[n-1]}(\tilde T_w^{(2)}))} \right)^{-1}.
\end{align}
Now write $\tilde \Rcal_{[n]} := \Rcal_{\tilde \Tcal_v}(v, \tilde G_{[n]}(\tilde \Tcal_v))$ for the effective resistance between $v$ and $\tilde G_{[n]}(\tilde \Tcal_v)$, and write $\tilde\Rcal_{[n-1]}^{(1)}$ and $\tilde\Rcal_{[n-1]}^{(2)}$ for independent copies of $\tilde\Rcal_{[n-1]}$.
Then, under $\bar{\mathbf{P}}_p$, we have the following equality in distribution:
\begin{align} \label{effresistancerec2}
\tilde \Rcal_{[n]} \stackrel{\mathrm{d}}{=} S + \left(\frac{1}{\tilde\Rcal_{[n-1]}^{(1)}}+\frac{1}{\tilde\Rcal_{[n-1]}^{(2)}}\right)^{-1},
\end{align}
where $S \sim $ Geo$(1-\hat p_1)$.
Bounding the harmonic mean by the arithmetic mean, i.e., using that $(\frac{1}{x} + \frac{1}{y})^{-1} \le \tfrac14 (x+y)$ for $x,y > 0$, this yields the stochastic domination
\begin{align} \label{effresistancerec3}
\tilde \Rcal_{[n]} \preccurlyeq S + \tfrac{1}{4} ( \tilde\Rcal_{[n-1]}^{(1)}+\tilde\Rcal_{[n-1]}^{(2)} ).
\end{align}
Iterating \eqref{effresistancerec3}, we obtain
\begin{align} \label{effresistanceineq4}
\tilde \Rcal_{[n]} \preccurlyeq \sum_{i=0}^{n-1} 4^{-i} \left( \sum_{j=1}^{2^i} S_{i,j} \right) ,
\end{align}
where $(S_{i,j})_{i,j\geq 1}$ is an array of i{.}i{.}d.\ copies of $S$.
Taking the expectation, 
\begin{align} \label{effresistanceineq5}
\lim_{n\to \infty} \bar{\mathbf{E}}_p[\tilde\Rcal_{[n]}] 
\leq \lim_{n\to \infty} \sum_{i=0}^{n-1} 2^{-i} \E_p[S] = 2 \E_p[S] \leq \frac{2\lambda}{p-p_c} .
\end{align}
Combining \eqref{effresistanceineq1}, \eqref{effresistanceineq2}, \eqref{effresistanceineq3} and \eqref{effresistanceineq5}, we obtain
\begin{align*}
 \mathbb{P}^v_p(\eta_0 <\infty \mid \mathcal{T}_{[0,1]}=T_{[0,1]}) \leq \frac{p-p_c}{L} \frac{2\lambda}{p-p_c} = \frac{2\lambda}{L},
\end{align*} 
and this upper bound is smaller than $\tfrac{2}{3}$ when we choose $L \ge L_0 := 3\lambda$. \qed

\subsection{Proof of Lemma \ref{lem:backtrackprob2}}\label{sec:backtrack2}

Define the event $\mathcal{G}:=\{|G_{[1]}(\mathcal{T})|>1\}$
and, for $\alpha < \tfrac{1}{4}$, let $\mathcal{F}(\alpha):=\{|G'_{\lfloor \sqrt{\alpha}L/(p-p_c) \rfloor }(\mathcal{T})|=1\}$ denote the  event that the backbone tree does not branch before the $\lfloor \sqrt{\alpha}L/(p-p_c) \rfloor $-th generation. 
Splitting the expectation along $\mathcal{G}$ and $\mathcal{F}(\alpha)$, we bound
\begin{align}\label{escapebound}
  \bar{\mathbf{E}}_p \bigg[ \sum_{v\in G_{[1]}(\mathcal{T})} \mathbbm{1}_{\{ P^{v}_{\mathcal{T}\setminus \mathcal{T}_{v}}(\eta_{0}=\infty) < \alpha  \}}\bigg]  
 	&\leq \bar{\mathbf{P}}_p(\mathcal{G}^c) + \bar{\mathbf{E}}_p \left[ \mathbbm{1}_{\mathcal{F}(\alpha)^c} |G_{[1]}| \right]  \\
	&\quad + \bar{\mathbf{E}}_p \left[ \mathbbm{1}_{\mathcal{G} \cap \mathcal{F}(\alpha)} \sum_{v\in G_{[1]}(\mathcal{T})} \bar{\mathbf{P}}_p( P^{v}_{\mathcal{T}\setminus \mathcal{T}_{v}}(\eta_{0}=\infty) < \alpha \mid \mathcal{T}_{[0,1]} ) \right] . \notag
\end{align}
By Lemma~\ref{lem:bb} the first term on the right-hand side is bounded from above by $a_2$, as desired.
It remains to show that the second and third term on the right-hand side above both vanish as $\alpha \to 0$.

We start with the second term. Let $B_k$ be the event that the first branching of the backbone tree occurs in generation $k$, that is, $|G'_1|=\dots =|G'_{k-1}|=1$, but $|G'_k|>1$. At time $B_k$, the backbone branch splits into at most $\Delta$ branches. Let $G_{[1]}^{(i)}, 1\leq i \leq \Delta$ be i{.}i{.}d.\ copies of $G_{[1]}$, then we can stochastically dominate $|G_{[1]}|$ by $|G^{(1)}_{[1]}|+\dots +|G^{(\Delta)}_{[1]}|$, on any of the events $B_k$, $k\leq \lfloor \sqrt{\alpha}L/(p-p_c) \rfloor$. This implies
\begin{align*}
\bar{\mathbf{E}}_p \left[ \mathbbm{1}_{\mathcal{F}(\alpha)^c} |G_{[1]}| \right] 
& = \sum_{k=1}^{\lfloor \sqrt{\alpha}L/(p-p_c) \rfloor } \bar{\mathbf{E}}_p \left[ \mathbbm{1}_{B_k} |G_{[1]}| \right] \\
& \leq \sum_{k=1}^{\lfloor \sqrt{\alpha}L/(p-p_c) \rfloor } \bar{\mathbf{E}}_p \left[ \mathbbm{1}_{B_k} \sum_{i=1}^{\Delta} |G_{[1]}^{(i)}| \right]  = \Delta \bar{\mathbf{P}}_p(\mathcal{F}(\alpha)^c) \bar{\mathbf{E}}_p \left[ |G_{[1]}| \right] , 
\end{align*}
where we used the independence of $B_k$ and the $|G_{[1]}^{(i)}|$ and Wald's identity for the last equality. 
By Lemma \ref{lem:bb}, $\bar{\mathbf{E}}[|G_{[1]}(\mathcal{T})| ]\leq a_3$ and, by 
 an argument similar to \eqref{e:unifprob}, there exists a constant $c>0$ such that $\bar{\mathbf{P}}_p(\mathcal{F}(\alpha)^c) \le C \sqrt{\alpha}$ when $\alpha$ is sufficiently close to $0$, so we conclude
 \begin{equation}\label{e:midbd}
 	 \bar{\mathbf{E}}_p \left[ \mathbbm{1}_{\mathcal{F}(\alpha)^c} |G_{[1]}| \right] \le C \sqrt{\alpha} a_3^2.
	\end{equation}

Now we show that the third term on the right-hand side of \eqref{escapebound} also tends to $0$ as $\alpha \to 0$.
Let $T\in \mathcal{G}\cap \mathcal{F}(\alpha)$ and 
denote by $v_0$ the unique vertex of $G'_{\lfloor \sqrt{\alpha}L/(p-p_c) \rfloor }(T)$. Then, for any $v\in G_{[1]}$, 
\begin{align*}
P^{v}_{T\setminus T_{v}}(\eta_{0}=\infty) \geq P^{v_0}_{T\setminus T_{v}}(\eta_{0}=\infty) = \frac{\mathcal{C}_{T\setminus T_{v}}(v_0,\infty)}{\mathcal{C}_{T\setminus T_{v}}(v_0,\varrho)+\mathcal{C}_{T\setminus T_{v}}(v_0,\infty)}
\end{align*}
by \eqref{effresistanceineq0}. By the Series Law and the fact that $T \in \mathcal{G}\cap \mathcal{F}(\alpha)$ we have $\mathcal{C}_{T\setminus T_{v}}(v_0,\varrho) = \left\lfloor \sqrt{\alpha}L / (p-p_c)\right\rfloor^{-1}$. It follows that $P^{v_0}_{T\setminus T_{v}}(\eta_{0}=\infty) < \alpha$
can only hold if
\begin{align*}
\mathcal{R}_{T\setminus T_{v}}(v_0,\infty) > \frac{1-\alpha}{\sqrt{\alpha}} \frac{L}{p-p_c}.
\end{align*}
This implies that for any vertex $w\in G_{[1]}\setminus \{ v \}$,
\begin{align*}
\bar{\mathbf{P}}_p\left( P^{v}_{\mathcal{T}\setminus \mathcal{T}_{v}}(\eta_{0}=\infty) < \alpha  \mid \mathcal{T}_{[0,1]} \right) 
& \leq \bar{\mathbf{P}}_p\left( \left. \mathcal{R}_{T\setminus T_{v}}(v_0,\infty) > \frac{1-\alpha}{\sqrt{\alpha}} \frac{L}{p-p_c} \,\right|\, \mathcal{T}_{[0,1]} \right) \\
& \leq  \bar{\mathbf{P}}_p\left(\left. \frac{L}{p-p_c} + \mathcal{R}_{T_w}(w,\infty) > \frac{1-\alpha}{\sqrt{\alpha}} \frac{L}{p-p_c} \,\right|\, \mathcal{T}_{[0,1]} \right) \\
& \leq 4 \sqrt{\alpha}\frac{p-p_c}{L}\bar{\mathbf{E}}_p \left[\mathcal{R}_{T_w}(w,\infty) \right],
\end{align*}
where we have used the independence of $\mathcal{T}_w$ and $\mathcal{T}_{[0,1]}$ and the fact that $\alpha<1/4$. We established in \eqref{effresistanceineq5} above that $\bar{\mathbf{E}}_p \left[\mathcal{R}_{T_w}(w,\infty) \right]\leq C(p-p_c)^{-1}$. Inserting this and \eqref{e:midbd} into \eqref{escapebound} to obtain
\begin{align}
 \bar{\mathbf{E}}_p\left[ \sum_{v\in G_{1}(\mathcal{T})} \mathbbm{1}_{\{ P^{v}_{\mathcal{T}\setminus \mathcal{T}_{v}}(\eta_{0}=\infty) < \alpha  \}}\right]  \leq a_2 + C \sqrt{\alpha}, 
\end{align}
which completes the proof. \qed

\subsection{Proof of Lemma \ref{lem:onlychild}}\label{sec:onlychild}

The main idea in this proof is to decompose the random walk on $\mathcal{T}_{[0,2]}$ into random walks on $\mathcal{T}_{[0,1]}$ and walks on the disjoint trees of $\mathcal{T}_{[1,2]}$. Let $(Y_n)_{n\geq 1}$ be a simple random walk on $\mathcal{T}_{[0,1]}$ started at $v$, and write $\mathcal{H}=\sigma(Y_1,Y_2,\dots ,Y_{\eta_0} )$ for the history of this process until it hits the root $\varrho$. 
For $i\geq 1$, let $h_i$ be the $i$-th time that $(Y_n)_{n\geq 1}$ visits a vertex in $G_{[1]}\cup \{\varrho \}$ and let $V_i = Y_{h_i}$. Furthermore, let $H$ denote the number of visits of $(Y_n)_{n \ge 1}$ to $G_{[1]}\cup \{\varrho \}$ until its first visit to the root, i.e., $H= \inf\{ i:\, V_i = \varrho \}$. We can decompose the walk on $\Tcal_{[0,2]}$ into the walk $(Y_n)_{n \ge 1}$ and the walks on the branches of $\Tcal_{[0,2]} \setminus \Tcal_{[0,1]}$, and write
\begin{align}\label{nextlevel}
& \mathbb{P}_p^v( \eta_2<\eta_0, \sib(X_{\eta_2})=0 \mid \mathcal{T}_{[0,1]} = T_{[0,1]} )  \notag \\
& = \mathbb{E}_p\left[ \left. \mathbb{E}_p \left[ \sum_{m=1}^{H-1} P_\mathcal{T}^{V_m}(\eta_2<\eta_1^{\mathrm{p}}) \mathbbm{1}_{\{|G_{[1]}(\mathcal{T}_{V_m})|=1 \}} \left. \prod_{i=1}^{m-1} P_\mathcal{T}^{V_i}(\eta_1^{\mathrm{p}} <\eta_2) \right| \mathcal{H} \right] \right| \mathcal{T}_{[0,1]} = T_{[0,1]} \right] ,
\end{align} 
where we write $\eta_1^{\mathrm{p}}$ for the hitting time of a parent vertex of $G_{[1]}$ (i.e., the hitting time of level $\lfloor L (p-p_c)^{-1}\rfloor -1$). Observe that conditioned on $V_i$,  $P_\mathcal{T}^{V_i}(\eta_2 < \eta_1^{\mathrm{p}})$ and  $ \mathbbm{1}_{\{|G_{[1]}(\mathcal{T}_{V_i})|=1 \}}$ are independent of $\mathcal{T}_{[0,1]}$. By \eqref{effresistanceineq0} we have
\begin{align} \label{hittingprobnextlevel}
 P_\mathcal{T}^{V_m}(\eta_2<\eta_1^{\mathrm{p}}) \mathbbm{1}_{\{|G_{[1]}(\mathcal{T}_{V_m})|=1 \}} & = 
\frac{\mathcal{C}_{\mathcal{T}_{V_m}}(V_m,G_{[1]}(\mathcal{T}_{V_i}))}{1+\mathcal{C}_{\mathcal{T}_{V_m}}(V_m,G_{[1]}(\mathcal{T}_{V_i}))}  \mathbbm{1}_{\{|G_{[1]}(\mathcal{T}_{V_m})|=1 \}} \notag \\ 
 & \geq c(p-p_c)\mathbbm{1}_{\{|G_{[1]}(\mathcal{T}_{V_m})|=1 \}},
\end{align}
and similarly 
\begin{equation} \label{hittingprobnextlevel2}
P_\mathcal{T}^{V_i}(\eta_1^{\mathrm{p}} <\eta_2)  = \frac{1}{1+\mathcal{C}_{\mathcal{T}_{V_m}}(V_m,G_{[1]}(\mathcal{T}_{V_i}))}
\geq 1- c|G_{[1]}(\mathcal{T}_{V_i})|(p-p_c).
\end{equation}

Define for $m<H$ and $v\in G_{[1]}$
\begin{align}
H_m(v) := | \{k< m \, : \, Y_{h_k}=v \} | , 
\end{align}
which counts the number of visits of $(Y_n)_{n\geq 1}$ to vertex $v$ before the $m$-th visit to a vertex in $G_{[1]}$. Reordering the product in \eqref{nextlevel} according to the vertices and applying the bound in  \eqref{hittingprobnextlevel} gives the following lower bound,
\begin{multline}
\mathbb{E}_p \left[ \sum_{m=1}^{H-1} P_\mathcal{T}^{V_m}(\eta_2<\eta_1^{\mathrm{p}}) \mathbbm{1}_{\{|G_{[1]}(\mathcal{T}_{V_m})|=1 \}} \left. \prod_{i=1}^{m-1} P_\mathcal{T}^{V_i}(\eta_1^{\mathrm{p}} <\eta_2) \right| \mathcal{H} \right] \\
 \geq c (p-p_c)  \mathbb{E}_p \left[ \left. \sum_{m=1}^{H-1} \mathbbm{1}_{\{|G_{[1]}(\mathcal{T}_{V_m})|=1 \}} \prod_{v\in G_{[1]}} P_\mathcal{T}^{v}(\eta_1^{\mathrm{p}} <\eta_2)^{H_m(v)}  \right| \mathcal{H} \right].
\end{multline}
Now note that $H$, $H_m(v)$ and $V_m$ are measurable with respect to $\mathcal{H}$, and the probabilities $P_\mathcal{T}^{v}(\eta_1^{\mathrm{p}} <\eta_2)$ are independent of $\mathcal{H}$ and independent under $\bar{\mathbf{P}}_p(\,\cdot \mid \mathcal{T}_{[0,1]})$ for different $v$, so
\begin{multline} \label{hittingprobnextlevel3}
\mathbb{E}_p \left[ \sum_{m=1}^{H-1} P_\mathcal{T}^{V_m}(\eta_2<\eta_1^{\mathrm{p}} ) \mathbbm{1}_{\{|G_{[1]}(\mathcal{T}_{V_m})|=1 \}} \left. \prod_{i=1}^{m-1} P_\mathcal{T}^{V_i}(\eta_1^{\mathrm{p}} <\eta_2) \right| \mathcal{H} \right] \\
 \geq c (p-p_c)  \sum_{m=1}^{H-1} \mathbb{E}_p[ \mathbbm{1}_{\{|G_{[1]}(\mathcal{T}_{V_m})|=1 \}}P_\mathcal{T}^{V_m}(\eta_1^{\mathrm{p}} <\eta_2)^{H_m(V_m)} \mid \mathcal{H}] \\
 \times \prod_{\substack{v\in G_{[1]} : \\ v\neq V_m}} \mathbb{E}_p [ P_\mathcal{T}^{v}(\eta_1^{\mathrm{p}} <\eta_2)^{H_m(v)} \mid \mathcal{H}] .
\end{multline}
Applying now the bound \eqref{hittingprobnextlevel2}, we see that for $p$ sufficiently close to $p_c$, \eqref{hittingprobnextlevel3} is bounded from below by 
\begin{align}
&  c (p-p_c)  \sum_{m=1}^{H-1} \bar{\mathbf{P}}_p (|G_{[1]}|=1 )(1-c(p-p_c))^{H_m(V_m)} \prod_{\substack{v\in G_{[1]}: \\ v\neq V_m}}  \big(1-c\bar{\mathbf{E}}_p[|G_{[1]}|](p-p_c)\big)^{H_m(v)} \notag \\
& \geq c (p-p_c) \bar{\mathbf{P}}_p(|G_{[1]}|=1 ) \sum_{m=1}^{H-1} \prod_{v\in G_{[1]}} \big(1-c\bar{\mathbf{E}}_p[|G_{[1]}|](p-p_c)\big)^{H_m(v)} \notag \\
& = c (p-p_c) \bar{\mathbf{P}}_p(|G_{[1]}|=1 ) \sum_{m=1}^{H-1} \big(1-c\bar{\mathbf{E}}_p[|G_{[1]}|](p-p_c)\big)^{m}  \\
& =  c \bar{\mathbf{E}}_p[|G_{[1]}|]^{-1} \bar{\mathbf{P}}_p(|G_{[1]}|=1 ) \left( \big(1-c\bar{\mathbf{E}}_p[|G_{[1]}|](p-p_c)\big)-\big(1-c\bar{\mathbf{E}}_p[|G_{[1]}|](p-p_c)\big)^{H}\right) .\notag
\end{align}

By Lemma \ref{lem:bb}, $\bar{\mathbf{P}}_p(|G_{[1]}(\mathcal{T})|=1)\geq a_1$ and $\bar{\mathbf{E}}_p[|G_{[1]}(\mathcal{T})| ]\leq a_3$. So to show a uniform lower bound for \eqref{nextlevel}, it remains to show that
\begin{align} \label{nextlevel2}
\mathbb{E}_p\left[ \left. \big(1-c\bar{\mathbf{E}}[|G_{[1]}|](p-p_c)\big)^{H} \right| \mathcal{T}_{[0,1]} = T_{[0,1]} \right] < 1- a
\end{align}
for an $a>0$ independent of $p$. While the exact distribution of $H$ depends on $\mathcal{T}_{[0,1]}$, we may give an easy lower bound. Recall that $H-1$ is the number of visits of $(Y_n)_{n\geq 1}$ to $G_{[1]}$ until it hits the root, starting at some $v\in G_{[1]}$. For any tree $T$ with $v\in G_{[1]}(T)$ we have
\begin{align}
P_T^v(\eta_0<\eta_1^+) \leq \mathcal{C}_T(v,\varrho) \leq c(p-p_c) ,
\end{align}
where $\eta_1^+$ denotes the first hitting time of $G_{[1]}$ after time $0$.
This bound implies that $H$ stochastically dominates a Geometric random variable $H'$ with success probability $c(p-p_c)$ and generating function 
\begin{align*}
\mathbb{E}_p[\theta^{H'}] = \frac{c(p-p_c)\theta}{1-\theta+c\theta(p-p_c)} 
\end{align*}
for $0\leq\theta\leq 1$, which in turn implies that
\begin{align}
&  \mathbb{E}_p\left[ \left. \big(1-c\bar{\mathbf{E}}_p[|G_{[1]}|](p-p_c)\big)^{H} \right| \mathcal{T}_{[0,1]} = T_{[0,1]} \right] 
\leq \mathbb{E}_p\left[  \big(1-c\bar{\mathbf{E}}_p[|G_{[1]}|](p-p_c)\big)^{H'}  \right]  \notag \\
& \qquad \qquad   = 1-\frac{c\bar{\mathbf{E}}_p[|G_{[1]}|]}{c\bar{\mathbf{E}}_p[|G_{[1]}|] + c\big(1-c\bar{\mathbf{E}}_p[|G_{[1]}|](p-p_c)\big)} .
\end{align}
The last term is bounded away from 0 uniformly in $p$, which shows \eqref{nextlevel2} and completes the proof. \qed


\section{Moment bounds on regeneration distances: proof of Lemma \ref{lem:moments}}\label{sec:regdist}  
This proof is inspired by a similar moment estimate of \cite{Demetal02} and a regularity estimate for trees from \cite{GriKes01}. 
By Lemmas~\ref{lem:bb}, \ref{lem:backtrackprob}, and \ref{lem:stationarity2}, 
\begin{align*}
\mathbb{E}_p[(|X_{\tau_{2}}|-|X_{\tau_1}|)^q] & =  \mathbb{E}_p^v[|X_{\tau_1}|^q \mid G_{[1]}(\mathcal T)=\{v\} , \eta_0=\infty ]\\ & \leq  C \mathbb{E}_p[|X_{\tau_1}|^q] ,
\end{align*}
so it suffices to find a moment estimate for $|X_{\tau_1}|$. Recall that by Definition~\ref{def:regen}, the first regeneration occurs at one of the potential regeneration times $S_k$. More precisely, the regeneration time is set to equal the last $S_k$ that is finite, so that
\[
	\begin{split}
		\mathbb{E}_p[|X_{\tau_1}|^q] 
		& \leq \sum_{k\geq 1} \mathbb{E}_p[|X_{S_k}|^q\mathbbm{1}_{\{ S_k<\infty \} }] \\
		& \leq \sum_{k\geq 1} \mathbb{E}_p[|X_{S_k}|^{2q}\mathbbm{1}_{\{ S_k<\infty \} }]^{1/2} \mathbb{P}_p(S_k<\infty )^{1/2} .
	\end{split}
\]
For $S_k$ to be finite, it has to occur at least $k$ times that the walker sees a new part of the tree but then backtracks a generation, so by Lemma \ref{lem:backtrackprob},
\begin{align*}
\mathbb{P}_p(S_k<\infty ) \leq \Big( \sup_{T:v\in G_{[1]}(T)} \mathbb{P}_p^v \big(\eta_0<\infty \mid \mathcal{T}_{[0,1]} = T_{[0,1]} \big)\Big)^k \leq \left(\tfrac{2}{3}\right)^k .
\end{align*}
By Definition~\ref{def:regen}, we can also write 
\begin{align} \label{recursionS_k}
\frac{|X_{S_k}|}{\lfloor L / (p-p_c) \rfloor} = M_k =: N_{k-1} + H_{k} , \qquad M_1 =: H_1
\end{align}
so that $H_k= M_k-N_{k-1}$. That is, $H_k$ counts the number of $L$-levels that are visited by the walk after it surpasses the previous maximum generation $N_{k-1}$, until it arrives at the next possible regeneration point, i.e., at a vertex $v$ in an $L$-level with $\sib(v)=0$. Then $H_k$ can be stochastically dominated by $2+\tilde{H}_k$, where $\tilde{H}_k$ is a Geometric random variable with success probability 
\begin{align*}
\inf_{T:v\in G_{[1]}(T)} \mathbb{P}_p^v( \eta_2<\eta_0,  \sib(X_{\eta_2})=0 \mid \mathcal{T}_{[0,1]} = T_{[0,1]} ) \geq a_3, 
\end{align*}
where the lower bound holds by Lemma \ref{lem:onlychild}. Since Geometric random variables have all moments bounded, it also holds for any $q >0$ that $\E_p[H_k^{2q}] < \infty$.

Write $\tilde N_k := N_k-M_k$. Iterating the recursion \eqref{recursionS_k}, and noting that $\{S_k<\infty\}\subset \{R_{k-1}<\infty\}$, we arrive at
\begin{align*}
(p-p_c)^q \mathbb{E}_p[|X_{\tau_1}|^q] \leq & L^q \sum_{k\geq 1}\mathbb{E}_p\left[ \left( H_1+ \sum_{i=1}^{k-1}(\tilde{N}_i+H_{i+1})\right)^{2q} \mathbbm{1}_{\{ R_{k-1}<\infty \} } \right]^{1/2}\left(\tfrac{2}{3}\right)^{k/2} \\
\leq & L^q \left(\sum_{k\geq 1} (2k)^{2q-1} \left( \sum_{i= 1}^k \mathbb{E}_p[H_i^{2q}] + \sum_{i= 1}^{k-1} \mathbb{E}_p[\tilde{N}_i^{2q}\mathbbm{1}_{\{ R_i<\infty \} }]\right)\right)^{1/2}\left(\tfrac{2}{3}\right)^{k/2} .
\end{align*} 
Since the $(2q)-$th moment of $H_i$ is uniformly bounded, the proof will be completed once we show that
\begin{align} \label{momboundN}
\mathbb{E}_p\left[ \mathrm{e}^{s \tilde N_i}\mathbbm{1}_{\{ R_i<\infty \} }\right] \leq C<\infty,
\end{align}
for some $s >0$. Observe that $\tilde N_i$ counts how many new $L$-levels are reached by the walker at $X_{\eta_{M_i}}$ before backtracking to its $L$-ancestor. Since for different $i$ these excursions happen in disjoint parts of the tree, the $\tilde N_i$ are in fact i{.}i{.}d.\ under $\mathbb{P}_p$. It thus suffices to bound
\begin{align*}
{\mathbb{E}}_p\left[ \mathrm{e}^{s \tilde N_1}\mathbbm{1}_{\{ R_1<\infty \} }\right] = \sum_{n\geq 1} \mathrm{e}^{s n} \mathbb{P}_p(\tilde N_1=n, R_1<\infty) .
\end{align*}

On the event $\{R_1<\infty\}$, a large value for $\tilde N_1$ implies that the random walk backtracks a long distance towards the root. We will bound the probability by decomposing this trajectory into level-sized chunks. Since the segment of the backbone tree between $X_{\eta_{M_i}}$ and its $L$-ancestor is by definition of $\eta_{M_i}$ isomorphic to a line graph of length $\lfloor L / (p-p_c)\rfloor$, the first step in this decomposition is to bound
\[
	\mathbb{P}_p(\tilde N_1=n, R_1<\infty) \leq  \bar{\mathbf{E}}_p\left[ \left. \sum_{v\in G_{[n]}(\mathcal{T})} P^{v_1}_\mathcal{T}(X_{\eta_n}=v) P_\mathcal{T}^v(\eta_0<\infty) \right| |G_{[1]}(\mathcal{T})|= 1 \right] ,
\]
where $v_1$ is the unique vertex of $G_{[1]}(\mathcal{T})$. Note that Lemma~\ref{lem:backtrackprob} is not sufficient to bound $P_\mathcal{T}^v(\eta_0<\infty)$, because the lemma gives an annealed bound, whereas we need a quenched bound. 
We instead use Lemmas~\ref{lem:bb} and \ref{lem:backtrackprob2} to show that the event
\begin{align*}
B_n(\alpha,\beta) := \left\{ \sum_{i=1}^n \mathbbm{1}_{\{ P^{v_i}_{\mathcal{T}\setminus \mathcal{T}_{v_i}}(\eta_{i-1}=\infty) \geq \alpha \} } \geq \beta n  \text{ for all } v\in G_{[n]}(\mathcal{T}) \right\} ,
\end{align*}
has a high probability provided $\alpha$ and $\beta$ are small enough, where, for $v=v_n\in G_{[n]}(T)$, we denote by $v_i$ its ancestor in $G_{[i]}(T)$ for $0\leq i < n$. 
If we show this, then we are done, because, for $\mathcal{T} \in B_n(\alpha,\beta)$ and $v\in G_{[n]}(\mathcal{T})$, 
\begin{align} \label{longbacktrack2}
P_\mathcal{T}^v(\eta_0<\infty) \leq \prod_{i=1}^n P^{v_i}_{\mathcal{T}\setminus \mathcal{T}_{v_i}}(\eta_{i-1}<\infty) \leq (1-\alpha)^{\beta n},
\end{align}
so that
\begin{equation}\label{e:goodbadsplit}
	\begin{split}
		\mathbb{P}_p(\tilde N_1=n, R_1<\infty) &\leq  (1-\alpha)^{\beta n} \bar{\mathbf{E}} \bigg[ \sum_{v \in G_{[n]}(\Tcal)} P^{v_1}_\Tcal (X_{\eta_n} = v) \, \bigg| \, |G_{[1]}(\mathcal{T})|= 1 \bigg]\\
		&\qquad + \bar{\mathbf{P}}_p\big(B_n(\alpha,\beta)^c \mid |G_{[1]}(\mathcal{T})|= 1  \big)\\
		& = (1-\alpha)^{\beta n} + \bar{\mathbf{P}}_p\big(B_n(\alpha,\beta)^c \mid |G_{[1]}(\mathcal{T})|= 1 \big),
	\end{split}
\end{equation}
where for the equality we have used that conditionally on $\sib(v_1)=0$, the probabilities add to $1$ because $G_{[n]}(\Tcal)$ is a cutset on $\Tcal$ separating the root from infinity.
To show that $B_n(\alpha,\beta)^c$ also has exponentially small probability, we introduce
\begin{align*}
A_n(\alpha,v)  := \sum_{i=1}^n \mathbbm{1}_{\{ P^{v_i}_{\mathcal{T}\setminus \mathcal{T}_{v_i}}(\eta_{i-1}=\infty) \geq \alpha \}} ,
\end{align*}
and
\begin{align*}
Z_n(\alpha, \theta) := \sum_{v\in G_n(\mathcal{T})} \mathrm{e}^{-\theta A_n(\alpha, v) } .
\end{align*}
We want to show that $Z_n(\alpha, \theta)$ decays exponentially for $\theta$ large enough. Note that by asking for the event $\{\eta_{i-1}=\infty\}$ only on the cropped tree $ \mathcal{T}\setminus \mathcal{T}_{v_i}$ we have independence of $P_{\Tcal \setminus \Tcal_{v_i}}^{v_i}(\eta_{i-1} = \infty)$ for different $i$ under the environment law,
which allows us to calculate recursively
\[
	\begin{split}
		\bar{\mathbf{E}}_p[Z_n(\alpha, \theta) & \mid \sib(v_1)=0] \\ 
		& =\bar{\mathbf{E}}_p \Bigg[  \sum_{v\in G_{[n-1]}(\mathcal{T})}  \mathrm{e}^{-\theta A_{n-1}(\alpha, v) } 
		\bar{\mathbf{E}}_p\bigg[ \sum_{w\in G_{[1]}(\mathcal{T}_v)} \mathbbm{1}_{\{ P^{w}_{\mathcal{T}\setminus \mathcal{T}_{w}}(\eta_{i-1}		=\infty) < \alpha  \}} \\
		& \qquad \qquad \qquad\qquad\qquad+ \mathrm{e}^{-\theta} \mathbbm{1}_{\{ P^{w}_{\mathcal{T}\setminus \mathcal{T}_{w}}(\eta_{i-1}=\infty) \geq \alpha  \}}\bigg| |G_{[n-1]}(\mathcal{T}) \bigg] \Bigg|\,  |G_{[1]}(\mathcal{T})|= 1 \Bigg] \\
		& =  \bar{\mathbf{E}}_p[Z_{n-1}(\alpha, \theta) \mid |G_{[1]}(\mathcal{T})|= 1] \zeta(\alpha, \theta) = \zeta(\alpha, \theta)^{n-1} ,
	\end{split}
\]
where 
\begin{align*}
\zeta(\alpha, \theta) := \bar{\mathbf{E}}_p \left[ \sum_{v\in G_{[1]}(\mathcal{T})} \mathbbm{1}_{\{ P^{v}_{\mathcal{T}\setminus \mathcal{T}_{v}}(\eta_{0}=\infty) < \alpha  \}} + \mathrm{e}^{-\theta} \mathbbm{1}_{\{ P^{v}_{\mathcal{T}\setminus \mathcal{T}_{v}}(\eta_{0}=\infty) \geq \alpha  \}} \right] .
\end{align*}
In the last step in the recursion we have used that $\bar{\mathbf{E}}_p [Z_1(\alpha, \theta)|\,|G_{[1]}(\mathcal{T})|= 1 ]=1$, since $|G_{[1]}(\mathcal{T})|= 1$ implies that $P^{v_1}_{\mathcal{T}\setminus \mathcal{T}_{v_1}}(\eta_{0}=\infty) =0 $.

We proceed by bounding $\zeta(\alpha, \theta)$ as
\begin{align} \label{zetabound}
\zeta(\alpha, \theta) & \leq  \bar{\mathbf{E}}_p \left[ \sum_{v\in G_{1}(\mathcal{T})} \mathbbm{1}_{\{ P^{v}_{\mathcal{T}\setminus \mathcal{T}_{v}}(\eta_{0}=\infty) < \alpha  \}}\right] + \mathrm{e}^{-\theta} \bar{\mathbf{E}}_p [|G_{[1]}(\mathcal{T})| ]. 
\end{align}
By Lemma \ref{lem:bb} we have $\bar{\mathbf{E}}_p [|G_{[1]}(\mathcal{T})| ] \leq a_3$.  
With the bound of Lemma \ref{lem:backtrackprob2}, this means that we can bound
\begin{align*}
\zeta(\alpha, \theta) \leq a_2+h(\alpha)+\mathrm{e}^{-\theta}a_3=: \gamma(\alpha, \theta) ,
\end{align*}
with $\gamma(\alpha, \theta)<1$ if $\alpha$ is small enough and $\theta$ is large enough. Fix such a sufficient choice of $\alpha$ and $\theta$, then, by Markov's inequality,
\begin{align*}
 \bar{\mathbf{P}}_p(B_n(\alpha,\beta)^c|\,|G_{[1]}(\mathcal{T})|= 1  ) 
 & =  \bar{\mathbf{P}}_p \left( \left. \min_{v\in G_{[n]}(\mathcal{T})} A_n(\alpha, v)<\beta n \right| \,|G_{[1]}(\mathcal{T})|= 1 \right) \\
 & = \bar{\mathbf{P}}_p \left( \left. \mathrm{e}^{-\theta \min_v A_n(\alpha, v)} >\mathrm{e}^{-\theta \beta n} \right| \,|G_{[1]}(\mathcal{T})|= 1 \right) \\
 & \leq \mathrm{e}^{n \beta \theta}  \bar{\mathbf{E}}_p \left[ \left. \mathrm{e}^{-\theta \min_v A_n(\alpha,v)} \right| \,|G_{[1]}(\mathcal{T})|= 1 \right] \\
 & \leq \mathrm{e}^{n\beta \theta} \bar{\mathbf{E}}_p [Z_n(\alpha,\theta)|\,|G_{[1]}(\mathcal{T})|= 1 ] \\
 & \leq \mathrm{e}^{n\beta \theta} \gamma(\alpha,\theta)^{n-1},
\end{align*}
and this bound decays exponentially in $n$ if $\beta $ is small enough. Inserting this bound into \eqref{e:goodbadsplit} we see that \eqref{momboundN} holds for some $s >0$ sufficiently small, which completes the proof. 
\qed


\section{Moment bounds on the regeneration times: proof of Lemma \ref{lem:moments2}}\label{sec:regtime}

We start by establishing moment bounds for the time spent in the finite trees $\Tcal^{\trap}_v$, which as the name suggests, act as traps for the random walk. Let $H_v :=\inf\{n\geq 0 \,:\, X_n=v\}$ be the hitting time of vertex $v$ and $H_v^+=\inf\{n> 0 \,:\, X_n=v\}$. When $v=\varrho$ we will suppress the subscript. If $\pi_T$ is the invariant distribution for random walk on a tree $T$, then the well-known formula 
\begin{align}\label{returntime}
E_T^v[H_v^+] = \pi_T(v)^{-1}
\end{align}
holds. For the second moment we have the following identity, see \cite[(2.21)]{AldFil02},
\begin{align}\label{returntime2}
E_T^v[(H_v^+)^2] = \pi_T(v)^{-1}\left( 2 \sum_{w\in T} \pi_T(w)E_T^w[H_v] +1\right) .
\end{align}
Furthermore, we may bound $E_T^w[H_v]$ by the \emph{commute time} $E_T^w[H_v]+E_T^v[H_w]$. The Commute Time Identity of \cite{Chaetal96} applied to simple random walk on a finite tree $T$ gives
\begin{align}
E_T^w[H_v]+E_T^v[H_w] = 2 (|T|-1) d_T(v,w),
\end{align}
with $|T|$ the vertex cardinality of the tree and $d_T(v,w)$ the graph distance between $v$ and $w$ on $T$. This implies
\begin{align}
E_T^v[(H_v^+)^2]\leq 4 \pi_T(v)^{-1} (|T|-1) \sum_{w\in T} \pi_T(w) d_T(v,w) .
\end{align}
Moreover, for simple random walk on a finite connected tree $T$, 
\[
	\pi_T(v)= \frac{\deg_T(v)}{\sum_{w \in T} \deg_T(w)} = \frac{\deg_T(v)}{2(|T|-1)}.
\]
So we arrive at
\begin{align} \label{returntime3}
E_T^v[(H_v^+)^2]\leq 4 \operatorname{deg}_T(v)^{-1} (|T|-1) \sum_{w\in T} \operatorname{deg}_T(w) d_T(v,w) .
\end{align} 
Now, from here on, consider instead of $T$ the tree $\Tcal^{\trap}_v$ located at $v \in \Tcal^{\Bb}$, i.e., a tree rooted at $v$ where $v$ has degree $U_v$ (given by \eqref{e:Uvdef}) and where to the $i$th edge, $i \in \{1,\dots,U_v\}$, coming out of the root, there is an independent tree $\Tcal^*_i$ with generating function $f^*(s)$ as defined in \eqref{treedecomposition} and with root $\varrho_i$ attached to the other end of the $i$th edge. We know that $\Tcal^{\trap}_v$ is almost surely finite. Write $H^\trap_v$ for $H^+_v$ of random walk started at $v$, restricted to $\Tcal^\trap_v$, with the convention that $H^\trap_v=0$ if $\Tcal^{\trap}_v=\{v\}$, that is, when $U_v=0$. Then, by \eqref{returntime},
\begin{equation}
	\begin{split}
		\mathbb{E}_p[H^\trap_v\mid U_v] &= \mathbf{E}_p[\pi_{\Tcal^\trap_v}(v)^{-1} \mid U_v ] \indi_{\{U_v \neq 0\}}\\
		& = 2\mathbf{E}_p\left[\left. U_v^{-1} \left(\left(\sum_{i=1}^{U_v} |  \mathcal{T}_i^*|\right)-1\right)\, \right|\, U_v\right] \indi_{\{U_v \neq 0\}} \\
		&  = 2\mathbf{E}_p [|\Tcal^*|-U_v^{-1}] \indi_{\{U_v \neq 0\}} .  
	\end{split}
\end{equation}
Applying Lemma~\ref{lem:gf} we obtain
\begin{align} \label{traptimelowerbound}
	\mathbb{E}_p[H^\trap_v\mid U_v] = 2(1-\mu^*)^{-1}\indi_{\{U_v \neq 0\}} +O(1)= C(p-p_c)^{-1}\indi_{\{U_v \neq 0\}} +o((p-p_c)^{-1}).
\end{align}

From \eqref{returntime3} we get for the second moment
\begin{equation}\label{traptimesecmom}\begin{split}
		\mathbb{E}_p[(H^\trap_v)^2|U_v] &\leq  4\mathbf{E}_p\left[\left.  U_v^{-1} \left(\sum_{i=1}^{U_v}|\mathcal{T}^*_i|\right)\left( \sum_{i=1}^{U_v} \sum_{w\in \mathcal{T}^*_i} \operatorname{deg}_{\Tcal^\trap_v}(w) d_{\Tcal^\trap_v} (w,v) \indi_{\{U_v \neq 0\}}\right)\right| U_v \right]\\
		&\leq  4\sum_{i=1}^\Delta \mathbf{E}_p\left[|\mathcal{T}^*_i| \sum_{w\in \mathcal{T}^*_i} (\operatorname{deg}_{\Tcal^*_i}(w)+1) (d_{\Tcal^*_i} (w,\varrho_i)+1) \right] \\
		& \qquad + 4\sum_{i\neq j} \mathbf{E}_p\left[|\mathcal{T}^*_i| \sum_{w\in \mathcal{T}^*_j} (\operatorname{deg}_{\Tcal^*_j}(w)+1) (d_{\Tcal^*_j} (w,\varrho_j)+1)  \right]	 \\	
		&= 4\Delta  \mathbf{E}_p\left[|\mathcal{T}^*| \sum_{w\in \mathcal{T}^*} (\operatorname{deg}_{\Tcal^*}(w)+1) (d_{\Tcal^*} (w,\varrho)+1) \right] \\
		& \qquad + 4\Delta(\Delta-1) \mathbf{E}_p\left[|\mathcal{T}^*|\right] \mathbf{E}_p\left[ \sum_{w\in \mathcal{T}^*} (\operatorname{deg}_{\Tcal^*}(w)+1) (d_{\Tcal^*} (w,\varrho)+1) \right] \\
	 & =: (I) + (I\! I) ,
\end{split}
\end{equation}
where we have used $U_v\leq \Delta$, the independence of the different GW-trees in a trap and in the second step we added one to the degree because $\deg_{\Tcal^\trap_v}(\varrho_i) = \deg_{\Tcal^*_i}(\varrho_i)+1$.
Now write $Z_n^*$ for the size of the $n$th generation of $\Tcal^*$ and observe that $|\Tcal^*| =\sum_{w \in \Tcal^*} 1 =\sum_{n \ge 0}Z_n^*$, and that $\sum_{w \in \Tcal^*} \deg_{\Tcal^*}(w) = \sum_{n \ge 0} (Z_n^* + Z_{n+1}^*)$ (each vertex has degree equal to the number of its offspring plus one, the latter accounting for its ancestor), so that we may bound
\[\begin{split}
	(I) &\le 4 \Delta \mathbf{E}_p\left[\left(\sum_{n \ge 0}Z_n^*\right)\left(\sum_{n \ge 0} (2 Z_n^* + Z_{n+1}^*)(n+1)\right)\right]\\
	 &\le C \mathbf{E}_p \left[\left(\sum_{n \ge 0}Z_n^*\right)\left(\sum_{n \ge 0} n Z_n^* \right)\right],
	 \end{split}
\]
where in the last step we used the assumption that $\Delta < \infty$.
By conditioning on the earlier generations we can write the right-hand side as
\[
	\begin{split}
		(I) & \le C \sum_{0 \le m<n} \mathbf{E}_p[  Z^*_m \mathbf{E}_p [  nZ^*_n \mid Z_m^*]] + C \sum_{n \ge 0} \mathbf{E}_p[ n (Z^*_n)^2]\\
		&\qquad +C \sum_{m>n \ge 0} \mathbf{E}_p[ \mathbf{E}_p[  Z^*_m \mid Z^*_n]  nZ^*_n ] \\
		& =  C \sum_{0 \le m<n} n \mathbf{E}_p[  (Z^*_m)^2 ]\mathbf{E}_p [  Z^*_{n-m}] + C \sum_{n \ge 0} n \mathbf{E}_p[ (Z^*_n)^2]\\
		&\qquad +C \sum_{m>n \ge 0} n\mathbf{E}_p[ Z^*_{m-n}] \mathbf{E}_p[  (Z^*_n)^2 ] .  
	\end{split}
\]
With the inequality $\mathbf{E}_p[(Z^*_n)^2]\leq C n (\mu^*)^{2n}$, this simplifies to
\begin{align}\label{traptimeupperbound}
(I) & \leq C \sum_{0\le m<n} n m (\mu^*)^{n+m} + C\sum_{n \ge 0} n^2 (\mu^*)^{2n} + C \sum_{m>n \ge 0} n^2 (\mu^*)^{n+m} \notag \\
& \leq  C  \left( \sum_{n\geq 0} n(\mu^*)^{n} \right)^2 + C \sum_{m \ge n\geq 0} n^2 (\mu^*)^{n+m} \notag \\
& = \frac{C(\mu^*)^2(2+\mu^*)^2}{(1-\mu^*)^4} \le  \frac{C}{(p-p_c)^4},
\end{align}
where the final inequality is due to Lemma~\ref{lem:gf}.

An easier computation shows that
\begin{equation}\label{traptimeupperbound2}
	(I\! I) \le C \E_p[\Tcal^*] \E_p\left[\sum_{n \ge 0} n Z_n^*\right] \le \frac{C}{(p-p_c)^3}.
\end{equation}
Combining the bounds \eqref{traptimeupperbound} and \eqref{traptimeupperbound2}, we may thus conclude that
\begin{equation} \label{traptimeupperbound3}
	\mathbb{E}_p[(H^\trap_v)^2\mid U_v] \le \frac{C}{(p-p_c)^4}.
\end{equation}

We now need to combine the time spent in the traps with the time the random walk spends on the backbone tree between regeneration times. For this, let $(X_m^{\Bb })_{m \ge 0}$ be the random walk restricted to the backbone tree without holding times, i.e., writing $t(m)$ for the time that the random walk $(X_n)_{n \ge 0}$ visits a vertex of $\Tcal^{\Bb}$, distinct from the previously visited vertex, for the $m$th time, and we let $X_m^{\Bb} := X_{t(m)}$.
If $v\in \mathcal{T}^{\Bb }$, then we write $\ell_n(v)$ for the local time of $(X_m^{\Bb })_{m \ge 0}$ at $v$ until time $n$. 
During a visit of $(X_m^{\Bb })_{m \ge 0}$ to $v\in \Tcal^{\Bb}$, the total time $(X_n)_{n \ge 0}$ spends in $\mathcal{T}^\trap_v$ is given by
\begin{align} \label{e:totaltraptime}
L_v = 1+\sum_{i=1}^{Y_v} H^\trap_{v,i} ,
\end{align} 
where for each $v\in \mathcal{T}^{\Bb}$, $H^\trap_{v,1},H^\trap_{v,2},\dots $ are independent random variables that, conditioned on $\mathcal{T}^\trap_v$, have the same distribution as $H^\trap_v$ under $P_{\mathcal{T}^\trap_v}^v$. That is, they count the number of steps spent in the trap until returning to $v$. Furthermore, conditioned on $\Tcal$, we let $Y_{v}$ be a Geo$(1-U_v / \deg_\Tcal(v))$ random variable (of the type that takes values in $\{0,1,2,\dots\}$), independent of $H^\trap_{v,i}$ for all $i$. So $Y_{v}$ counts the number of times that the random walk $(X_n)_{n \ge 0}$ at a visit to $v$ enters and exits the trap before moving on to a different vertex on the backbone tree.

If $\tau$ is a random time measurable with respect to $\sigma(X_m^{\Bb}, m\geq 0)$, then these definitions allow us to write
\begin{equation}\label{e:Etau}
	\E_p[\tau] = \E_p\left[\sum_{v \in \Tcal^{\Bb}} \sum_{i=1}^{\ell_\tau(v)} L_v^{(i)} \right],
\end{equation}
where $L_v^{(i)}$ are i.i.d. copies of $L_v$, with the convention that if $\ell_\tau(v) =0$, then the sum equals $0$.

We first apply this general formula to determine a lower bound for the inter-regeneration times. Recall from Definition~\ref{def:regen} the regeneration structure and the various definitions. We bound
\begin{equation}\label{e:taulb}
	\begin{split}
\mathbb{E}_p[\tau_2-\tau_1] & \ge \E_p[(\tau_2 - \tau_1) \indi_{\{|G_{[\Lambda_1 +1]}| =1\}}]\\
	& \geq \mathbb{E}_p[(\eta_{\Lambda_1+1}-\tau_1) \indi_{\{|G_{[\Lambda_1+1]}| =1\}}]  \\
	& \ge a_1 \E_p[\eta_{\Lambda_1+1} -\tau_1 \mid |G_{[\Lambda_1+1]}|=1 ] \\
	& = a_1 \E_p[\eta_{\Lambda_1+1} -\tau_1 \mid |G_{[\Lambda_1]}| = |G_{[\Lambda_1+1]}|=1, \eta_{\Lambda_1 +1} < \eta_{\Lambda_1 -1}]\\
	& = a_1 \E^w_p[\eta_2 - \eta_1 \mid G_{[1]}=\{w\}, |G_{[2]}|=1, \eta_2 < \eta_0] \\
& = a_1 \mathbb{E}_p^w\left[\left.\sum_{v\in \mathcal{T}^{\Bb}_{[0,2]}} \sum_{i=1}^{\ell_{\eta_2-\eta_1}(v)} L_v^{(i)} \,\right|\, G_{[1]}=\{w\}, |G_{[2]}|=1, \eta_2 < \eta_0 \right],
	\end{split}
\end{equation}
where in the third step we have used Lemma~\ref{lem:bb}, in the fourth step we have conditioned on further events that are guaranteed to occur by the definition of the regeneration times, and in the fifth step we have used that $\tau_1$ is a regeneration time to justify the time shift.

Conditionally on the event $\{|G_{[1]}|= |G_{[2]}|=1\}$, the graph $\Tcal^{\Bb}_{[0,2]}$ is isomorphic to a line graph of length $2\lfloor L / (p-p_c)\rfloor$, so every vertex, except for the root, and possibly the vertex at the second $L$-level, has degree $2$ in $\Tcal^{\Bb}_{[0,2]}$. 
Conditionally also on $U_v$, the $Y_{v,i}$ are distributed as independent Geo$(1-U_v /(U_v+2))$ random variables. 
From \eqref{e:totaltraptime} we get
\begin{align} \label{e:totaltraptimelb}
& \mathbb{E}_p^v[ L_v \mid U_v , \deg_{\Tcal^\Bb}(v) = 2 ] \notag \\
&\qquad = 1+\mathbb{E}_p^v\big[ Y_v \mid U_v , \deg_{\Tcal^\Bb}(v) = 2 \big] \mathbb{E}_p^v\big[  H^\trap_{v,1} \mid U_v , \deg_{\Tcal^\Bb}(v) = 2 \big] \notag \\
& \qquad = 1+ \tfrac12 U^{(2)}_v \mathbb{E}_p^v[  H^\trap_{v,1} \mid U_v ] \geq  c (p-p_c)^{-1} U^{(2)}_v ,
\end{align}
where $U^{(2)}_v$ is distributed as $U_v$ conditionally on $\deg_{\Tcal^\Bb}(v)=2$, and we have used \eqref{traptimelowerbound} for the last step.
Thus, for any $v \in \Tcal^{\Bb}_{[0,2]} \setminus \{\varrho, X_{\eta_2}\}$, 
\begin{multline}\label{e:condbdYH}
	\E_p^w\left[\left.\sum_{i=1}^{\ell_{\eta_2-\eta_1}(v)} L_v^{(i)} \,\right|\, U_v, \ell_{\eta_2-\eta_1}(v), G_{[1]}=\{w\}, |G_{[2]}|=1 ,\eta_2<\eta_{0}\right] \\
	\geq \ell_{\eta_2 -\eta_1}(v) c (p-p_c)^{-1} U^{(2)}_v.
\end{multline}
Using \eqref{e:Uvdef} and Lemma~\ref{lem:gf}, we can bound
\[
	\E_p[U^{(2)}_v] \ge  \P(U_v^{(2)} > 0) = 1-f'_p(0)/f'_p(q_p) \ge c.
\]
Combining this bound with \eqref{e:condbdYH}, inserting that into \eqref{e:taulb}, and writing $\eta_j^{\Bb}$ for the hitting time of $G_{[j]}$ for $(X_m^{\Bb})_{m \ge 0}$, we obtain
\[
	\E_p[\tau_2-\tau_1] \ge \frac{c}{p-p_c} \E_p^w \left[\eta_2^{\Bb} - \eta_1^{\Bb} -2 \mid G_{[1]} =\{w\}, |G_{[2]}|=1, \eta_2^{\Bb} < \eta_1^{\Bb} \right] \ge \frac{c'}{(p-p_c)^3},
\]
where, in the last step, we have used that $(X_m^{\Bb})_{m \ge 0}$ restricted to $\Tcal^{\Bb}_{[0,2]}$, conditionally on $\{|G_{[1]}| = |G_{[2]}|=1\}$, is equivalent to a simple random walk on a line graph of length $2\lfloor L / (p-p_c) \rfloor$.
This proves the lower bound in Lemma \ref{lem:moments2}. 
\medskip

For the upper bounds, we first notice that by the same reasoning as in Lemma \ref{lem:moments} it suffices to show the bound for $\tau_1$. Let us again first consider $L_v$ as in \eqref{e:totaltraptime}, for which we can bound
\begin{align} \label{e:totaltraptimeub}
\mathbb{E}_p^v[L_v^2 \mid  \Tcal^{\Bb}] & \leq 1+ \mathbb{E}_p^v[Y_v^2 \mid  \Tcal^{\Bb}]  \mathbb{E}_p^v[(H^\trap_{v,i})^2 \mid  \Tcal^{\Bb}] .
\end{align}
For the Geometric random variable $Y_v$, the parameter satisfies $1-U_v/\deg_\Tcal(v)\geq 1/\Delta $, which implies 
\begin{align} \label{e:totalmoment1}
\mathbb{E}_p^v[Y_v^2 \mid  \Tcal^{\Bb}] \leq \Delta^2 ,
\end{align}
while for the second expectation in \eqref{e:totaltraptimeub} we may apply \eqref{traptimeupperbound3} to conclude 
\begin{align} \label{e:totaltraptimeub2}
\mathbb{E}_p^v[L_v^2 \mid  \Tcal^{\Bb}] \leq C(p-p_c)^{-4} .
\end{align}

We can use this upper bound to state a bound on the second moment of a generic random time $\tau$ that is measurable with respect to $\sigma(X_n^{\Bb}, n\geq 0)$,
\begin{align}
\mathbb{E}_p[\tau^2] & = \mathbb{E}_p\left[ \left( \sum_{v\in \mathcal{T}^{\Bb}} \sum_{i=1}^{\ell_{\tau}(v)}L_v^{(i)} \right)^2 \right] \notag \\
& =  \mathbb{E}_p\left[  \sum_{v\in \mathcal{T}^{\Bb}} \left(\sum_{i=1}^{\ell_{\tau}(v)}L_v^{(i)} \right)^2 \right] 
 + \mathbb{E}_p\left[  \sum_{v\neq w \in \mathcal{T}^{\Bb}} \left(\sum_{i=1}^{\ell_{\tau}(v)}L_v^{(i)} \right)\left(\sum_{i=1}^{\ell_{\tau}(w)}L_w^{(i)} \right)\right] \notag \\
& =  \mathbb{E}_p\left[  \sum_{v\in \mathcal{T}^{\Bb}} \ell_{\tau}(v) \mathbb{E}^v[L_v^2 \mid  \Tcal^{\Bb}] + \ell_{\tau}(v)(\ell_{\tau}(v)-1)\mathbb{E}_p^v[L_v \mid  \Tcal^{\Bb}]^2 \right] \notag \\
& \qquad + \mathbb{E}_p\left[  \sum_{v\neq w \in \mathcal{T}^{\Bb}} \ell_{\tau}(v) \ell_{\tau}(w)\mathbb{E}_p^v[L_v \mid  \Tcal^{\Bb}]\mathbb{E}_p^w[L_w \mid  \Tcal^{\Bb}]\right] \notag \\
& \leq C (p-p_c)^{-4} \mathbb{E}_p[\tau^{\Bb}]  + C (p-p_c)^{-2} \mathbb{E}_p[(\tau^{\Bb})^2]  ,
\end{align}
where, for the third equality, we have used that traps at distinct vertices of the backbone tree are independent, and for the final inequality we have used \eqref{e:totaltraptimeub}. Applying this upper bound to $\tau_1^\Bb$, which denotes the number of steps of $(X_n^\Bb)_n$ until $\tau_1$, the upper bounds in Lemma~\ref{lem:moments2} will therefore follow once we show
\begin{align} \label{momentboundBB}
\mathbb{E}_p[(\tau_1^{\Bb })^2]\leq C(p-p_c)^{-4} .
\end{align}

Since $|X_n^{\Bb }|$ stochastically dominates a simple random walk on $\mathbb{N}_0$ reflected at the origin, $\eta_m^{\Bb }$ is stochastically dominated by the time the simple random walk hits the vertex $m \lfloor L/(p-p_c)\rfloor$. We have $ \mathbb{E}_p[(\eta_m^{\Bb})^k]\leq Cm^{2k}(p-p_c)^{-2k} $ for any $k\geq 1$, so
\begin{align} \label{momentboundBB2}
\mathbb{E}_p[(\tau_1^{\Bb})^q] & = \sum_{m=1}^\infty  \mathbb{E}_p[(\eta_m^{\Bb})^q\mathbbm{1}_{\{  X_{\tau_1} \in G_{[m]} \} }]\notag \\
& \leq  \sum_{m=1}^\infty  \mathbb{E}_p[(\eta_m^{\Bb})^{2q}]^{1/2} \mathbb{P}_p ( X_{\tau_1} \in G_{[m]})^{1/2} \notag \\
& \leq  C (p-p_c)^{-2q} \sum_{m=1}^\infty  m^{2q} \mathbb{P}_p ( X_{\tau_1} \in G_{[m]})^{1/2} .
\end{align}
By the moment bound of Lemma \ref{lem:moments}, $\mathbb{P}_p ( (p-p_c)|X_{\tau_1}| \geq m) \leq Cm^{-4q-4}$. This implies \eqref{momentboundBB}, which concludes the proof. \qed


\section{The maximal displacement of a trap: proof of Lemma \ref{lem:nzlemma}}\label{sec:maxdisp}

The claimed moment bound follows from the following bound for the projected processes, 
\begin{align} \label{nzestimate}
 \max_{1 \le \delta \le \Delta-1} \mathbb{E}_p\left[ \left. \exp \left( \gamma \sqrt{p-p_c} \sup_{w\in  \mathcal{T}^{\trap}_\varrho} \varphi(w)\cdot \mathbf{e} \right) \right| \deg_{\Tcal^\trap_\varrho} (\varrho) = \delta \right] \leq C ,
\end{align}
for constants $\gamma,C$, arbitrary unit vectors $\mathbf{e}$ with $\|\mathbf{e}\|=1$ and where, without loss of generality, we consider a trap at the root, such that $\varphi(\varrho)=0$.

Recall from the decomposition in \eqref{treedecomposition} that $\mathcal{T}^{\trap}_\varrho$ consists of at most $\Delta-1$ GW-trees attached by a single edge to $\varrho$, each tree having generating function $f^*$ and mean number $\mu^*_p$ of offspring. Thus, it suffices to show \eqref{nzestimate} with $\mathcal{T}^{\trap}_\varrho$ replaced by a single tree $\mathcal T^*$ with generating function $f^*$. Let $W$ denote a $\Zd$-valued random variable with distribution $D$, and write $K(s) = \mathbb{E}_p[s^{W \cdot \mathbf{e}} ]$ for the probability generating function of $W \cdot \mathbf{e}$. For $x>1$, let $\theta(x)$ be the unique solution in $(1,\infty)$ to 
\begin{align} \label{nzestimate2}
K(\theta(x)) = x.
\end{align}
Neuman and Zheng \cite[Theorem~1.2]{NeuZhe17} show that
\begin{align}  \label{nzestimate1}
\limsup_{n\to \infty} \theta ((\mu^*_p)^{-1})^n \mathbb{P}_p\bigg( \sup_{w\in  \mathcal{T}^*} \varphi(w)\cdot \mathbf{e} \geq n \bigg) \leq 1 
\end{align}
(it is here that we need the assumption $\mathbb{E}_p[\mathrm{e}^{c||W||}]=\sum_{x \in \Zd} \e^{c \|x\|}D(x) < \infty$ for all $c>0$). 
The assertion of Lemma \ref{lem:nzlemma} thus easily follows if it holds that
\begin{align} \label{nzestimate3}
\theta ((\mu^*_p)^{-1})^{-n(p-p_c)^{-1/2}} \leq  \e^{-cn}  ,
\end{align}
for $n$ sufficiently large, uniformly in $p$ close enough to $p_c$. We proceed by showing that this is the case. 

From Lemma \ref{lem:gf} we obtain
\begin{align} \label{nzestimate4}
(\mu^*_p)^{-1} = 1+c(p-p_c)(1+o(1)) .
\end{align}
Since $\mathbb{E}_p[W] = \sum_{x \in \Zd} x D(x) =0$ by the assumption in Theorem~\ref{thm:main} we have $\E_p[W \cdot \mathbf{e}] =0$, which implies $K'(1)=0$, so that expansion around $s=1$ gives 
\begin{align}
K(s) \leq 1+c(s-1)^2 ,
\end{align}
for $s$ close enough to 1 (again using that $\mathbb{E}_p[\mathrm{e}^{c||W||}]<\infty$). This in turn implies that for $x$ close enough to 1, 
\begin{align} \label{nzestimate5}
\theta(x)\geq 1+c(x-1)^{1/2} .
\end{align}
Combining \eqref{nzestimate4} and \eqref{nzestimate5}, as well as $(\mu^*_p)^{-1}\geq 1+c(p-p_c)$, we arrive at
\begin{align} \label{nzestimate6}
\theta ((\mu^*_p)^{-1})^{-n(p-p_c)^{-1/2}} \leq \left(1+c(p-p_c)^{1/2} \right)^{-n(p-p_c)^{-1/2}} \le \e^{-c n} ,
\end{align}
for $p$ close enough to $p_c$. This implies \eqref{nzestimate3}, which finished the proof.
\qed

\section{The trace of the walk on the backbone: proof of Lemma \ref{lem:BBTbounds}}\label{sec:BBT}

We will show a general moment bound for the trace until the second regeneration, that is, a moment bound for the cardinality of
\begin{align} \label{BBTbound1}
\operatorname{BBT}_1 \cup \operatorname{BBT}_2 = \{X_0,X_1,\dots ,X_{\tau_{2}} \} \cap \mathcal{T}^{\Bb } 
=  \{X_0^\Bb,X^\Bb_1,\dots ,X^\Bb_{\tau_{2}^\Bb} \} . 
\end{align}
The bound on $\E_p[|\operatorname{BBT}_k|]$ for any $k\geq 1$ then follows from Lemma \ref{lem:stationarity2}. 
We denote by $\mathcal{B}=\mathcal{B}(\mathcal{T})$ the set of branch points of the backbone tree, that is,
\begin{align*}
\mathcal{B}(\mathcal{T}) := \{\varrho\} \cup \{ v\in \mathcal{T}^\Bb \, : \, \operatorname{deg}_{\Tcal^\Bb} (v) \ge 3 \} .
\end{align*}
The bound for the trace on the backbone tree will follow from a bound on the number of visited branch points of the backbone. Between those branchpoints, the random walk has to travel across a section of the tree isomorphic to a line segment. For $v\in \mathcal{T}^\Bb$, let $\anc_{\mathcal{B}(\mathcal{T})}(v)$ denote the first vertex on the shortest path between $v$ and $\varrho$ that is in $\mathcal{B}(\mathcal{T})$. Given $v \in \mathcal{B}(\Tcal)$, we call the path in $\Tcal$ between $v$ and $\anc_{\mathcal{B}(\Tcal)}(v)$ the \emph{backbone branch to $v$}.

Fix a small $\delta>0$ (to be determined below) and let $\mathcal{S}(\mathcal{T})$ be the set of all vertices $v \in \mathcal{T}^\Bb$ for which there exists a path connecting $v$ to the root $\varrho$, such that for any vertex $w$ on this path,  
\begin{align*}
d_\Tcal ( w, \anc_{\mathcal{B}(\mathcal{T})}(w)) \leq \lfloor \delta /(p-p_c)\rfloor,
\end{align*}
i.e., $v$ is in $\mathcal{S}(\Tcal)$ if the path from $v$ to $\anc_{\mathcal{B}(\Tcal)}(v)$ and all backbone branches along the path to $\varrho$ have length at most $\lfloor \delta /(p-p_c)\rfloor$.

The motivation for these definitions is as follows: Until the random walk exits $\mathcal{S}(\mathcal{T})$, the trace of the random walk on the backbone is a subset of $\mathcal{S}(\mathcal{T})$, and in order to exit $\mathcal{S}(\Tcal)$ the random walk has to cross a backbone branch of length at least $\lfloor \delta /(p-p_c)\rfloor$. We can find a lower bound for the time it takes to traverse such a backbone branch, which implies an upper bound for the size of the trace on the backbone until that time. When the random walk exits $\mathcal{S}(\mathcal{T})$ and crosses a long backbone branch, it enters a new, independent subtree, and in this new tree we can iterate this upper bound.

So set $\mathcal{S}_1 :=\mathcal{S}(\mathcal{T})$ and let $E_1 :=H_{\mathcal{S}_1^c}^\Bb$ be the exit time of the set $\mathcal{S}_1$ for the random walk $(X_n^\Bb)_{n\geq 1}$ restricted to the backbone tree. Define recursively for $k>1$, 
\begin{align} \label{backboneshrubs}
\mathcal{S}_k := \mathcal{S}_{k-1} \cup \mathcal{S}(\mathcal{T}_{X_{E_{k-1}}}) \qquad \text{ and } \qquad E_k := H^\Bb_{\mathcal{S}_k^c} . 
\end{align} 
Then $E_k$ stochastically dominates the sum $\tilde E_k = \tilde E_1^{(1)} + \dots + \tilde E_1^{(k)}$ of $k$ i{.}i{.}d.\ copies of 
\begin{align}\label{linesegmentcrossings}
\tilde E_1 = \inf\{n\geq 0 \, : \, \, Y_n = \lfloor \delta /(p-p_c)\rfloor \} , 
\end{align}
 where $(Y_n)_{n\geq 1}$ is a simple random walk on $\{0,\dots , \lfloor \delta /(p-p_c)\rfloor \}$ starting in and reflected at 0. Moreover, 
 $|\mathcal{S}_k|$ stochastically dominates $|\{X_0,\dots ,X_{E_k}\}\cap \mathcal{T}^\Bb|$. Since at each time $E_i$ the set $\mathcal{S}(\mathcal{T}_{X_{E_{i}}})$ is independent of the tree explored so far, we can bound
\begin{align} \label{backboneshrubs2}
\mathbb{E}_p[|\{X_0,\dots ,X_{E_k}\}\cap \mathcal{T}^\Bb|^q]\leq  \mathbb{E}_p[|\mathcal{S}_k|^q] \leq \mathbb{E}_p \left[ \left( \sum_{i=1}^k \tilde{\mathcal{S}}^{(i)} \right)^q \right] \leq k^q \mathbb{E}_p[|\mathcal{S}_1|^q ]  , 
\end{align}
where we wrote $\tilde{\mathcal{S}}^{(i)}$ for i{.}i{.}d.\ copies of $\mathcal{S}_1$. 

The next step is to show that, for any $q\geq 1$, 
\begin{align} \label{backboneshrubsmoment}
 \mathbb{E}_p[|\mathcal{S}_1|^q ] \leq C_q (p-p_c)^{-q} . 
\end{align}
Recall that $\Delta$ is the maximal number of children in the original Galton-Watson tree. By definition of $\mathcal{S}_1$, 
\begin{align} \label{backboneshrubsmoment2}
|\mathcal{S}_1| \leq \Delta \cdot |\mathcal{S}_1 \cap \mathcal{B}(\mathcal{T}) | \cdot \lfloor \delta /(p-p_c)\rfloor ,  
\end{align}
where we have used that $\Delta$, the maximal number of offspring in the unpercolated tree, is finite. By construction, the set $|\mathcal{S}_1 \cap \mathcal{B}(\mathcal{T}) |$ is stochastically dominated by the total progeny of a Galton-Watson process, denoted by $\tilde{\mathcal{Z}}$, with offspring distribution given by the law of 
\begin{align}
\tilde \xi := \sum_{i=1}^{\Delta} \mathbbm{1}_{\{\ell_i \leq \lfloor \delta /(p-p_c)\rfloor \} } .
\end{align} 
Here, $\ell_i$ are independent random variables corresponding to the length of a backbone branch, so that each $\ell_i$ is an independent Geometric random variable, having mean $c(p-p_c)^{-1}$ (see Remark \ref{rem:treeshape}). Then $\tilde \xi$ is a Binomial random variable with $\Delta$ trials and probability of success 
\begin{align}
\mathbb{P}_p(\ell_i \leq \lfloor \delta /(p-p_c)\rfloor) \leq 1-(1-c(p-p_c))^{\lfloor \delta /(p-p_c)\rfloor} \leq 1-\mathrm{e}^{-c\delta} ,
\end{align}
where we have used Lemma \ref{lem:gf}. This means that we can choose $\delta>0$ small enough so that uniformly in $p$,
\begin{align}
 \mathbb{E}_p [\tilde \xi^{2q} ]\leq c<1 ,
\end{align}
i.e., $\tilde{\mathcal{Z}}$ is a subcritical Galton-Watson process. Writing $\tilde{\mathcal{Z}}_n$ for the size of the $n$-th generation of $\tilde{\mathcal{Z}}$ and using  $\mathbb{E}_p[\tilde{\mathcal{Z}}_n^{2q} ] \leq \mathbb{E}_p[\tilde{\mathcal{Z}}_{n-1}^{2q} ]  \mathbb{E}_p [\tilde \xi^{2q} ] \leq \mathbb{E}_p [\tilde \xi^{2q} ]^n$, we can bound the moments of the total progeny $|\tilde{\mathcal{Z}}|$ as 
\begin{align} \label{backboneshrubsmoment3}
 \mathbb{E}_p [|\tilde{\mathcal{Z}}|^q ] & = \sum_{N=0}^\infty \mathbb{E}_p \left[ \left(\sum_{n= 0}^N \tilde{\mathcal{Z}}_n\right)^q \mathbbm{1}_{\{\tilde{\mathcal{Z}}_{N}\neq 0 ,\, \tilde{\mathcal{Z}}_{N+1}= 0 \} }\right] \notag \\
 & \leq \sum_{N=0}^\infty \sum_{n= 0}^N N^{q-1} \mathbb{E}_p \left[  \tilde{\mathcal{Z}}_n^q \mathbbm{1}_{\{\tilde{\mathcal{Z}}_{N}\neq 0 \} }\right] .
\end{align} 
 Applying the Cauchy-Schwarz inequality and using $ \mathbb{P}_p(\tilde{\mathcal{Z}}_{N}\neq 0) \leq \mathbb{E}_p [\tilde{\mathcal{Z}}_{N}]= \mathbb{E}_p [\tilde \xi ]^N$, 
\begin{align}\label{backboneshrubsmoment3b}
 \mathbb{E}_p [|\tilde{\mathcal{Z}}|^q ]  & \leq \sum_{N=0}^\infty \sum_{n= 0}^N N^{q-1} \mathbb{E}_p [  \tilde{\mathcal{Z}}_n^{2q}]^{1/2} \mathbb{P}_p(\tilde{\mathcal{Z}}_{N}\neq 0 )^{1/2} \notag \\
 & \leq \sum_{N=0}^\infty N^{q-1} \left( \sum_{n= 0}^\infty \mathbb{E}_p [\tilde \xi^{2q} ]^{n/2} \right) \mathbb{P}_p(\tilde{\mathcal{Z}}_{N}\neq 0 )^{1/2} \notag \\
 & \leq (1- \mathbb{E}_p [\tilde \xi^{2q} ]^{1/2})^{-1} \sum_{N=0}^\infty N^{q-1} \mathbb{E}_p [\tilde \xi ]^{N/2} .
\end{align}
Since $\mathbb{E}_p [\tilde \xi ] \le \E_p[\tilde\xi^{2q}]<1$ uniformly in $p$, \eqref{backboneshrubsmoment3} is bounded uniformly in $p$ also. This implies that by the Cauchy-Schwarz inequality, \eqref{backboneshrubs2} and \eqref{backboneshrubsmoment},
\begin{align} \label{backboneshrubsmoment4} 
\mathbb{E}_p [|\mathcal{S}_1 \cap \mathcal{B}(\mathcal{T}) |^q ]\leq \mathbb{E}_p [|\tilde{\mathcal{Z}}|^q ]\leq C_q ,
\end{align}
 and, by \eqref{backboneshrubsmoment2}, \eqref{backboneshrubsmoment4} shows that \eqref{backboneshrubsmoment} holds.

So far, we have shown that \eqref{backboneshrubs2} is bounded by $C_q k^q (p-p_c)^{-q}$. In order to turn this into a bound for \eqref{BBTbound1}, we note that
\begin{align} \label{BBTbound2}
\mathbb{E}_p [ |\{X_0^\Bb,X^\Bb_1,\dots ,X^\Bb_{\tau_{2}^\Bb} \}|^{q/2} ] & \leq \sum_{k=1}^\infty  \mathbb{E}_p [ |\{X_0^\Bb,X^\Bb_1,\dots ,X^\Bb_{E_k} \}|^{q/2} \mathbbm{1}_{\{ E_{k-1}\leq \tau_2^\Bb < E_k\} } ] \notag \\
& \leq C_q (p-p_c)^{-q} \sum_{k=1}^\infty k^q \mathbb{P}_p(E_{k-1}\leq \tau_2^\Bb < E_k )^{1/2} \notag \\ 
& \leq C_q (p-p_c)^{-q} \sum_{k=1}^\infty k^q \mathbb{P}_p(E_{k-1}\leq \tau_2^\Bb )^{1/2} ,
\end{align}
where we set $E_0=0$. It remains to show that the sum is bounded.

Consider first the random variable $\tilde E_k$, which is stochastically dominated by $E_k$ and is a sum of $k$ i{.}i{.}d.\ random variables each distributed as $\tilde E_1$ in \eqref{linesegmentcrossings}. Letting the random walk $(Y_n)_{n \ge 1}$ in the definition \eqref{linesegmentcrossings} run for $(p-p_c)^{-2}$ steps and using the scaling limit of a simple random walk, it is easy to see that $(p-p_c)^{2}\tilde E_1$ dominates a Bernoulli distributed random variable with a success probability independent of $p$. 
This implies that by Cram\'er's Theorem, there is a $\gamma>0$ such that
\begin{align} \label{BBTbound3}
\mathbb{P}_p(\tilde{E}_k \leq k (p-p_c)^{-2} \gamma ) \leq C\mathrm{e}^{-ck} . 
\end{align}
Since $E_k$ dominates $\tilde E_k$, we also get
\begin{align} \label{BBTbound4}
\mathbb{P}_p({E}_k \leq k (p-p_c)^{-2} \gamma ) \leq C\mathrm{e}^{-ck} . 
\end{align}
For this choice of $\gamma$, we use that
\begin{align} \label{BBTprobsplit}
\mathbb{P}_p(E_{k}\leq \tau_2^\Bb ) \leq \mathbb{P}_p({E}_k \leq k (p-p_c)^{-2} \gamma ) + \mathbb{P}_p(\tau_2^\Bb \geq k (p-p_c)^{-2} \gamma ) .
\end{align}
On the other hand, the arguments in \eqref{momentboundBB2} with $\tau_1^\Bb$ replaced by $\tau_2^\Bb$ show that $(p-p_c)^{4q+8}\mathbb{E}_p[(\tau_2^\Bb)^{2q+4}] \leq C<\infty$, so Markov's inequality implies
\begin{align}\label{BBTbound5}
\mathbb{P}_p(\tau_2^\Bb \geq k (p-p_c)^{-2} \gamma ) \leq Ck^{-2q-4} .
\end{align}
The bounds \eqref{BBTbound4} and \eqref{BBTbound5} thus together with \eqref{BBTprobsplit} imply 
\begin{align} \label{BBTbound6}
\mathbb{P}_p(E_{k}\leq \tau_2^\Bb ) \leq Ck^{-2q-4} ,
\end{align}
which, when inserted into \eqref{BBTbound2}, proves the first assertion of Lemma \ref{lem:BBTbounds}. 
\medskip

To prove \eqref{backbonetracebound2}, the second assertion of Lemma \ref{lem:BBTbounds}, we first observe that the number of visited branchpoints until time $E_k$ satisfies
\begin{align}
\mathbb{E}_p[|\{X_0,\dots ,X_{E_k}\}\cap \mathcal{T}^\Bb\cap \mathcal{B}(\mathcal T)|^q] \leq k^q \mathbb{E}_p[|\mathcal{S}_1\cap \mathcal{B}(\mathcal T)|^q ] ,
\end{align}
and similarly to \eqref{backboneshrubs}, and by \eqref{backboneshrubsmoment4}, this moment is bounded by $C_q k^q$. If $\mathcal{N}_{E_k}$ denotes the number of distinct backbone branches visited by $(X_n^\Bb)_{n\geq0}$ before time $E_k$, then $\mathcal{N}_{E_k}$ is bounded by the number of visited branchpoints multiplied by the bound $\Delta$ for the number of offspring, which implies
\begin{align}
\mathbb{E}_p[(\mathcal{N}_{E_k})^q] \leq C_q \Delta^q k^q.
\end{align}
By the same arguments as in the preceding paragraph, we may then go on to conclude that 
\begin{align} \label{BBTbound7}
\mathbb{E}_p[(\mathcal{N}_{\tau_2^\Bb})^q] \leq C_q \Delta^q . 
\end{align}

Recall that the embedding $\varphi$ of the tree $\Tcal$ into $\mathbb{Z}^d$ is generated by the random walk one-step distribution $D$. Following the branching random walk interpretation, we can think of the process as a random walker that splits into $\xi$ random particles at each step. Therefore, we arrive at an upper bound for $\|\varphi(v)-\varphi(X_{\tau_{1}})\|$ by considering the maximal displacement among $\mathcal{N}_{\tau_2^\Bb}$ independent $|\operatorname{BBT}_2|$-step random walks started at $0$ with step distribution $D$. Let $(W^{(i)}_n)_{n\geq 0}$ for $ 1 \le i \leq \mathcal{N}_{\tau_2^\Bb}$ denote these independent random walks. Then 
\begin{align} \label{BBTbound8}
\mathbb{E}_p \big[\max_{v\in \operatorname{BBT}_2} \|\varphi(v)-\varphi(X_{\tau_{1}})\|^q \big] 
& \leq \mathbb{E}_p \Big[ \max_{1\leq i \leq \mathcal{N}_{\tau_2^\Bb} } \max_{0\leq k \leq m \leq |\operatorname{BBT}_2|} \|W_k^{(i)}- W_m^{(i)}\|^q \Big] \notag \\
& \leq \mathbb{E}_p \left[ \mathcal{N}_{\tau_2^\Bb} \right] \mathbb{E}_p \Big[ \max_{0\leq k \leq m \leq |\operatorname{BBT}_2|} \|W_k^{(1)}- W_m^{(1)}\|^q \Big] \notag \\
& \leq C \mathbb{E}_p \left[ \mathcal{N}_{\tau_2^\Bb} \right] \mathbb{E}_p \left[ |\operatorname{BBT}_2|^{q/2} \right]
\end{align}
where we have bounded the maximum over $i$ by the sum and used Wald's identity. For the last inequality we have used the independence between the increments of the embedding and $|\operatorname{BBT}_2|$, and that $\mathbb{E}_p \left[ \max_{0\leq k \leq m \leq N} \|W_k^{(1)}- W_m^{(1)}\|^q \right]\leq C N^{q/2}$ by the assumption that $D$ has finite exponential moments. By \eqref{BBTbound7}, the first expectation in \eqref{BBTbound8} is bounded by a constant, and by \eqref{backbonetracebound1} the second expectation is bounded by $C (p-p_c)^{q/2}$. This completes the proof of Lemma \ref{lem:BBTbounds}. 
\qed

\section{The speed of the walk: proof of Theorem \ref{thm:speed} and Lemma \ref{lem:gf}}\label{sec:speed}

In this section we prove Theorem \ref{thm:speed} and Lemma \ref{lem:gf}. We start with Lemma \ref{lem:gf}, and also derive a result that will be crucially used in the proof of Theorem \ref{thm:speed}. In the following, let 
	\begin{align*}
	m_{p,k} := f_p^{(k)}(1) = \mathbf{E}_p[\xi_p(\xi_p-1)\cdots (\xi_p-k+1)]
	\end{align*}
denote the $k$-th factorial moment of $\xi_p$. We abbreviate $m_{k}= m_{1,k}$, noting that $m_{p,k}=p^km_k$. Recall that $q_p$ is the extinction probability, i.e., the unique fixed point of $f_p$ in $(0,1)$.
\medskip

\noindent
{\bf Proof of Lemma \ref{lem:gf}.} We start by deriving the generating function $f_p$, for which we use that, conditioned on $\{\xi=n\}$, the number of offspring $\xi_p$ in the pruned tree is distributed as a Binomial with $n$ trials and success probability $p$. Therefore, 
\begin{align}
f_p(s) & = \mathbf{E}_p[\xi_p^s] = \sum_{k=0}^\infty \mathbf{E}_p[\xi_p^s \mid  \xi=n ] p_k \notag \\
& =   \sum_{k=0}^\infty (ps+(1-p))^k p_k = f(ps+(1-p)) ,
\end{align}
where $f$ is the generating function of the offspring distribution $(p_k)_{k\geq 0}$. Thus, we obtain that
	\eqn{
	f_p'(0)=pf'(1-p)\geq c_0>0,
	}
uniformly for $p\geq p_c$, as required. Furthermore, setting $s=1$ yields $f_p'(1)=p f'(1)=pm_1=p/p_c$. This proves the two equalities in \eqref{equal-lem2.1}.

We proceed with the three asymptotic relations in \eqref{asy-lem2.1}, where we start with the first. 

Since $q_p$ is the fixed point of $f_p$, we have
	\eqn{
	1-q_p=1-f_p(q_p)=1-f(pq_p+1-p)=1-f(1-p(1-q_p)).
	}
Taking the Taylor expansion of the right-hand side in $1-q_p$ and using that $q_p \to 1$ as $p \searrow p_c$ yields
	\eqn{
	1-q_p=pf'(1)(1-q_p)-\tfrac{1}{2}p^2 f''(1)(1-q_p)^2 +\tfrac{1}{6} p^3 f'''(1) (1-q_p)^3 +o((1-q_p)^{3}).
	}
Here we have used the assumption that $\mathbf{E}_p[\xi^3]<\infty$.
Using that $f^{(k)}(1)=m_k$ and dividing through by $1-q_p>0$ we obtain
	\eqn{
	\label{Taylor-exp-second-order}
	\tfrac{1}{2}p^2 m_2(1-q_p)=(pm_1-1) +\tfrac{1}{6} p^3 m_3 (1-q_p)^2 +o((1-q_p)^{2}).
	}
We will use these asymptotics in the proof of Theorem~\ref{thm:speed} below.
Replacing the last two terms by $O((1-q_p)^2)$, the first asymptotics in \eqref{asy-lem2.1} follows, as
	\eqn{
	\label{(1-qp)-asy}
	1-q_p=\frac{2}{m_2p^2} (pm_1-1)+O((1-q_p)^2)=\frac{2m_1^3}{m_2} (p-p_c)+O((1-q_p)^2),
	}
 with $c_2=2m_1^3/m_2$, since $p=p_c+O(p-p_c)=1/m_1+O(p-p_c)$.
	
To prove the second asymptotics in \eqref{asy-lem2.1}, we note that 
\[
	f'_p(q_p)=  pf'(q_p)=pf'(1-p(1-q_p))=pf'(1)-p^2f''(1)(1-q_p)+O((1-q_p)^2),
\]
so that, by \eqref{(1-qp)-asy},
	\eqan{
	f'_p(q_p)&=p m_1-(1-q_p)p^2m_2+O((1-q_p)^2) \notag \\
	&=1+m_1(p-p_c)-2m_1(p-p_c)+O((p-p_c)^2)\notag \\
	&=1-(p-p_c)m_1+O((p-p_c)^2).
	}
The third asymptotics in \eqref{asy-lem2.1} similarly follows, since 
\[\begin{split}
	\hat f_p''(0) &= (1-q_p) f''_p(q_p)= (1-q_p) p^2f''(1-p(1-q_p)) \\
	&= (1-q_p) p^2_cf''(1)(1+ o(1))=2m_1 (p-p_c) (1+o(1)).
\end{split}\]
\qed
\medskip

We proceed with the proof of Theorem \ref{thm:speed}, for which we rely on an explicit formula for the effective speed due to Lyons, Pemantle and Peres \cite[Page 601]{LyoPemPer96a},
	\begin{align*}
	v(p) = \sum_{k=0}^\infty \frac{k-1}{k+1} p_k(p) \frac{1-q_p^{k+1}}{1-q_p^2},
	\end{align*}
with $p_k(p) = \mathbf{P}_p(\xi_p=k)$. Using Lemma \ref{lem:gf}, we expand this expression for $p$ close to $p_c$.

\proof[Proof of Theorem \ref{thm:speed}] 
We again expand in terms of $1-q_p$. We rewrite
	\eqn{
	(1+q_p) v(p) = \sum_{k=0}^\infty \frac{k-1}{k+1} p_k(p) \sum_{i=0}^k q_p^i
	=\sum_{k=0}^\infty \frac{k-1}{k+1} p_k(p) \sum_{i=0}^k (1-(1-q_p))^i.
	}
Substituting the expansion
	\eqn{
	(1-(1-q_p))^i=1-i(1-q_p)+\tfrac{1}{2} i(i-1) (1-q_p)^2 +o\big(i^{2}(1-q_p)^{2}).
	}
gives three terms and an error term. Using $p_c=1/m_1=1/\mathbf{E}_p[\xi]$, the first term equals
	\eqan{
	\sum_{k=0}^\infty \frac{k-1}{k+1} p_k(p) (k+1)&=\mathbf{E}_p[\xi_p -1]=p\mathbf{E}_p[\xi]-1
	=m_1 (p-p_c) .
	}
Using $\sum_{i=0}^k i=\tfrac{1}{2} k(k+1)$, the second equals
	\eqn{
	-\tfrac{1}{2}(1-q_p) \sum_{k=0}^\infty k(k-1)p_k(p) =-\tfrac{1}{2}(1-q_p) m_{p,2}=-\tfrac{1}{2}(1-q_p) p^2 m_2.
	}
Using $\sum_{i=0}^k i(i-1)=\tfrac{1}{3} (k+1)k(k-1)$, the third equals
	\eqan{
	&\tfrac{1}{2} (1-q_p)^2 \sum_{k=0}^\infty \frac{k-1}{k+1} p_k(p) \sum_{i=2}^k i(i-1)\\
	&\quad=\tfrac{1}{2} (1-q_p)^2  \sum_{k=0}^\infty \frac{k-1}{k+1} p_k(p)\tfrac{1}{3}(k+1)k(k-1)\nonumber\\
	&\quad=\tfrac{1}{6} (1-q_p)^2  \Big(\sum_{k=0}^\infty k(k-1)(k-2)p_k(p)+\sum_{k=0}^\infty k(k-1)p_k(p)\Big)\nonumber\\
	&\quad=\tfrac{1}{6} (1-q_p)^2\big(m_{p,3}+m_{p,2}\big).\nonumber
	}
Using $\mathbf{E}_p[\xi^{3}]<\infty$ and $|(k-1)/(k+1)|\leq 1$, the error term can be estimated as
	\eqan{
	&o((1-q_p)^{2}) \sum_{k=0}^\infty p_k(p) \sum_{i=0}^k i^{2}
	=o((1-q_p)^{2}).
	}
We conclude that
	\eqan{
	(1+q_p) v(p)&=	m_1 (p-p_c)-\tfrac{1}{2}(1-q_p) p^2 m_2\\
	&\qquad +\tfrac{1}{6} (1-q_p)^2\big(m_{p,3}+m_{p,2}\big)+o((1-q_p)^{2}).\nonumber
	}
Substituting \eqref{Taylor-exp-second-order} for the second term, we obtain
	\eqan{
	m_1 (p-p_c)-\tfrac{1}{2}(1-q_p) p^2 m_2
	&=-\tfrac{1}{6} p^3 m_3 (1-q_p)^2+o((p-p_c)^{2}),
	}
which cancels up to leading order with the term $\tfrac{1}{6} (1-q_p)^2m_{p,3}$.  (We have no intuitive explanation for why the third moment drops out.)

By \eqref{(1-qp)-asy}, we arrive at
	\eqan{
	(1+q_p) v(p) &=\tfrac{1}{6}(1-q_p)^2 m_{p,2}+o\big((p-p_c)^{2}\big)\\
	&=\tfrac{2}{3} (p-p_c)^2 \frac{m_1^4}{m_2}+o\big((p-p_c)^{2}\big).\nonumber
	}
Using that $1+q_p=2+O(p-p_c)$ completes the proof.\hfill $\Box$

\bigskip

\textbf{Acknowledgements.} The work of JN was supported by the Deutsche Forschungsgemeinschaft (DFG) through grant NA 1372/1.  
RvdH was supported by NWO through VICI-grant 639.033.806. The work of RvdH and TH is also supported by the Netherlands Organisation for Scientific Research (NWO) through Gravitation-grant NETWORKS-024.002.003.

\bibliographystyle{alpha}
\bibliography{bib_nearcritGWT}

\newcommand{\etalchar}[1]{$^{#1}$}
\begin{thebibliography}{DMFGW89}

\bibitem[AF02]{AldFil02}
David Aldous and Jim Fill.
\newblock Reversible {M}arkov chains and random walks on graphs.
\newblock \textit{available at
  }\href{https://www.stat.berkeley.edu/~aldous/RWG/book.pdf}{\texttt{https://www.stat.berkeley.edu/$\sim$aldous/RWG/book.pdf}},
  2002.

\bibitem[Ald91]{Aldo91}
David Aldous.
\newblock The continuum random tree. {I}.
\newblock {\em Ann. Probab.}, 19(1):1--28, 1991.

\bibitem[ALW17]{athreya2017invariance}
Siva Athreya, Wolfgang L{\"o}hr, and Anita Winter.
\newblock Invariance principle for variable speed random walks on trees.
\newblock {\em Ann. Probab.}, 45(2):625--667, 2017.

\bibitem[AN72]{AthNey72}
Krishna~B. Athreya and Peter~E. Ney.
\newblock {\em Branching processes}.
\newblock Springer-Verlag, New York-Heidelberg, 1972.
\newblock Die Grundlehren der mathematischen Wissenschaften, Band 196.

\bibitem[BACF16a]{BenCabFri16b}
G{\'e}rard Ben~Arous, Manuel Cabezas, and Alexander Fribergh.
\newblock Scaling limit for the ant in a simple labyrinth.
\newblock \textit{arXiv:1609.03980}, 2016.

\bibitem[BACF16b]{BenCabFri16a}
G{\'e}rard Ben~Arous, Manuel Cabezas, and Alexander Fribergh.
\newblock Scaling limit for the ant in high-dimensional labyrinths.
\newblock \textit{arXiv:1609.03977}, 2016.

\bibitem[BAF16]{BenFri16}
G\'erard Ben~Arous and Alexander Fribergh.
\newblock Biased random walks on random graphs.
\newblock In {\em Probability and statistical physics in {S}t. {P}etersburg},
  volume~91 of {\em Proc. Sympos. Pure Math.}, pages 99--153. Amer. Math. Soc.,
  Providence, RI, 2016.

\bibitem[BB07]{BerBis07}
Noam Berger and Marek Biskup.
\newblock Quenched invariance principle for simple random walk on percolation
  clusters.
\newblock {\em Probab. Theory Related Fields}, 137(1-2):83--120, 2007.

\bibitem[BGP03]{BerGanPer03}
Noam Berger, Nina Gantert, and Yuval Peres.
\newblock The speed of biased random walk on percolation clusters.
\newblock {\em Probab. Theory Related Fields}, 126(2):221--242, 2003.

\bibitem[Bil99]{Bill99}
Patrick Billingsley.
\newblock {\em Convergence of probability measures}.
\newblock Wiley Series in Probability and Statistics: Probability and
  Statistics. John Wiley \& Sons, Inc., New York, second edition, 1999.
\newblock A Wiley-Interscience Publication.

\bibitem[Bis11]{biskup2011recent}
Marek Biskup.
\newblock Recent progress on the random conductance model.
\newblock {\em Probability Surveys}, 8, 2011.

\bibitem[BJKS08]{BarJarKumSla08}
M.T. Barlow, A.A. J{\'a}rai, T.~Kumagai, and G.~Slade.
\newblock Random walk on the incipient infinite cluster for oriented
  percolation in high dimensions.
\newblock {\em Comm. Math. Phys.}, {\bf 278}(2):385--431, (2008).

\bibitem[CR99]{CheRom99}
Hui Chen and Joseph~P. Romano.
\newblock An invariance principle for triangular arrays of dependent variables
  with application to autocovariance estimation.
\newblock {\em Canad. J. Statist.}, 27(2):329--343, 1999.

\bibitem[Cro09]{Croy09}
David~A. Croydon.
\newblock Hausdorff measure of arcs and {B}rownian motion on {B}rownian spatial
  trees.
\newblock {\em Ann. Probab.}, 37(3):946--978, 2009.

\bibitem[CRR{\etalchar{+}}97]{Chaetal96}
Ashok~K. Chandra, Prabhakar Raghavan, Walter~L. Ruzzo, Roman Smolensky, and
  Prasoon Tiwari.
\newblock The electrical resistance of a graph captures its commute and cover
  times.
\newblock {\em Comput. Complexity}, 6(4):312--340, 1996/97.

\bibitem[CS15]{CheSak15}
Lung-Chi Chen and Akira Sakai.
\newblock Critical two-point functions for long-range statistical-mechanical
  models in high dimensions.
\newblock {\em Ann. Probab.}, {\bf 43}(2):639--681, (2015).

\bibitem[dG76]{deGe76}
Pierre~Gilles de~Gennes.
\newblock Percolation: un concept unificateur.
\newblock {\em La Recherche}, 7:919--927, 1976.

\bibitem[DGPZ02]{Demetal02}
Amir Dembo, Nina Gantert, Yuval Peres, and Ofer Zeitouni.
\newblock Large deviations for random walks on {G}alton-{W}atson trees:
  averaging and uncertainty.
\newblock {\em Probab. Theory Related Fields}, 122(2):241--288, 2002.

\bibitem[DMFGW85]{DeMetal85}
Anna De~Masi, Pablo~A. Ferrari, Sheldon Goldstein, and W.~David Wick.
\newblock Invariance principle for reversible {M}arkov processes with
  application to diffusion in the percolation regime.
\newblock In {\em Particle systems, random media and large deviations
  ({B}runswick, {M}aine, 1984)}, volume~41 of {\em Contemp. Math.}, pages
  71--85. Amer. Math. Soc., Providence, RI, 1985.

\bibitem[DMFGW89]{DeMetal89}
Anna De~Masi, Pablo~A. Ferrari, Sheldon Goldstein, and W.~David Wick.
\newblock An invariance principle for reversible {M}arkov processes.
  {A}pplications to random motions in random environments.
\newblock {\em J. Statist. Phys.}, 55(3-4):787--855, 1989.

\bibitem[DS98]{DerSla98}
Eric Derbez and Gordon Slade.
\newblock The scaling limit of lattice trees in high dimensions.
\newblock {\em Comm. Math. Phys.}, 193(1):69--104, 1998.

\bibitem[EK86]{EthKur86}
Stewart~N. Ethier and Thomas~G. Kurtz.
\newblock {\em Markov processes}.
\newblock Wiley Series in Probability and Mathematical Statistics: Probability
  and Mathematical Statistics. John Wiley \& Sons Inc., New York, (1986).

\bibitem[Eva93]{Evan93}
Steven~N. Evans.
\newblock Two representations of a conditioned superprocess.
\newblock {\em Proc. Roy. Soc. Edinburgh Sect. A}, 123(5):959--971, (1993).

\bibitem[FvdH17]{FitHof17}
Robert Fitzner and Remco van~der Hofstad.
\newblock Mean-field behavior for nearest-neighbor percolation in {$d>10$}.
\newblock {\em Electron. J. Probab.}, 22:Paper No. 43, 65, 2017.

\bibitem[GGN17]{GanGuoNag2017}
Nina Gantert, Xiaoqin Guo, and Jan Nagel.
\newblock Einstein relation and steady states for the random conductance model.
\newblock {\em Ann. Probab.}, 45(4):2533--2567, 2017.

\bibitem[GK01]{GriKes01}
Geoffrey Grimmett and Harry Kesten.
\newblock Random electrical networks on complete graphs ii: Proofs.
\newblock \textit{arXiv preprint math/0107068}, 2001.

\bibitem[GMP12]{GanMatPia2012}
Nina Gantert, Pierre Mathieu, and Andrey Piatnitski.
\newblock Einstein relation for reversible diffusions in a random environment.
\newblock {\em Communications on Pure and Applied Mathematics}, 65(2):187--228,
  2012.

\bibitem[Guo16]{Guo2016}
Xiaoqin Guo.
\newblock Einstein relation for random walks in random environment.
\newblock {\em Ann. Probab.}, 44(1):324--359, 2016.

\bibitem[Hof06]{Hofs06a}
Remco~van~der Hofstad.
\newblock Infinite canonical super-{B}rownian motion and scaling limits.
\newblock {\em Comm. Math. Phys.}, 265(3):547--583, (2006).

\bibitem[HS90]{HarSla90a}
Takashi Hara and Gordon Slade.
\newblock Mean-field critical behaviour for percolation in high dimensions.
\newblock {\em Comm. Math. Phys.}, {\bf 128}(2):333--391, (1990).

\bibitem[HS00]{HarSla00a}
Takashi Hara and Gordon Slade.
\newblock The scaling limit of the incipient infinite cluster in
  high-di\-men\-sio\-nal percolation. {I}. {C}ritical exponents.
\newblock {\em J. Statist. Phys.}, {\bf 99}(5-6):1075--1168, (2000).

\bibitem[Hul15]{Huls15}
Tim Hulshof.
\newblock The one-arm exponent for mean-field long-range percolation.
\newblock {\em Electron. J. Probab.}, 20:no. 115, 26, 2015.

\bibitem[HvdH17]{HeyHof17}
Markus Heydenreich and Remco van~der Hofstad.
\newblock {\em Progress in high-dimensional percolation and random graphs}.
\newblock Springer, 2017.

\bibitem[HvdHH14]{HeyHofHul14a}
Markus Heydenreich, Remco van~der Hofstad, and Tim Hulshof.
\newblock Random walk on the high-dimensional {IIC}.
\newblock {\em Comm. Math. Phys.}, {\bf 329}(1):57--115, (2014).

\bibitem[HvdHS08]{HeyHofSak08}
Markus Heydenreich, Remco van~der Hofstad, and Akira Sakai.
\newblock Mean-field behavior for long- and finite range {I}sing model,
  percolation and self-avoiding walk.
\newblock {\em J. Stat. Phys.}, {\bf 132}(6):1001--1049, (2008).

\bibitem[JM05]{JanMar05}
Svante Janson and Jean-Fran{\c{c}}ois Marckert.
\newblock Convergence of discrete snakes.
\newblock {\em J. Theoret. Probab.}, {\bf 18}(3):615--647, (2005).

\bibitem[Kes86]{Kest86b}
Harry Kesten.
\newblock Subdiffusive behavior of random walk on a random cluster.
\newblock {\em Ann. Inst. H. Poincar\'e Probab. Statist.}, 22(4):425--487,
  (1986).

\bibitem[Kes95]{Kest95}
Harry Kesten.
\newblock Branching random walk with a critical branching part.
\newblock {\em J. Theoret. Probab.}, {\bf 8}(4):921--962, (1995).

\bibitem[KN09]{KozNac09}
Gady Kozma and Asaf Nachmias.
\newblock The {A}lexander-{O}rbach conjecture holds in high dimensions.
\newblock {\em Invent. Math.}, 178(3):635--654, 2009.

\bibitem[KV86]{KipVar86}
Claude Kipnis and S.~R.~Srinivasa Varadhan.
\newblock Central limit theorem for additive functionals of reversible {M}arkov
  processes and applications to simple exclusions.
\newblock {\em Comm. Math. Phys.}, 104(1):1--19, 1986.

\bibitem[LP16]{LyoPer16}
Russell Lyons and Yuval Peres.
\newblock {\em Probability on trees and networks}, volume~42 of {\em Cambridge
  Series in Statistical and Probabilistic Mathematics}.
\newblock Cambridge University Press, New York, 2016.

\bibitem[LPP95]{LyoPemPer96a}
Russell Lyons, Robin Pemantle, and Yuval Peres.
\newblock Ergodic theory on {G}alton-{W}atson trees: speed of random walk and
  dimension of harmonic measure.
\newblock {\em Ergodic Theory Dynam. Systems}, 15(3):593--619, 1995.

\bibitem[LPP96]{LyoPemPer96b}
Russell Lyons, Robin Pemantle, and Yuval Peres.
\newblock Biased random walks on {G}alton-{W}atson trees.
\newblock {\em Probab. Theory Related Fields}, 106(2):249--264, 1996.

\bibitem[Lyo92]{Lyon92}
Russell Lyons.
\newblock Random walks, capacity and percolation on trees.
\newblock {\em Ann. Probab.}, 20(4):2043--2088, 1992.

\bibitem[MP07]{MatPia07}
Pierre Mathieu and Andrey Piatnitski.
\newblock Quenched invariance principles for random walks on percolation
  clusters.
\newblock {\em Proc. R. Soc. Lond. Ser. A Math. Phys. Eng. Sci.},
  463(2085):2287--2307, 2007.

\bibitem[NZ17]{NeuZhe17}
Eyal Neuman and Xinghua Zheng.
\newblock On the maximal displacement of subcritical branching random walks.
\newblock {\em Probab. Theory Related Fields}, 167(3-4):1137--1164, 2017.

\bibitem[PZ08]{PerZei2008}
Yuval Peres and Ofer Zeitouni.
\newblock A central limit theorem for biased random walks on {G}alton--{W}atson
  trees.
\newblock {\em Probability Theory and Related Fields}, 140(3-4):595--629, 2008.

\bibitem[Sah94]{Sahi94}
Muhammad Sahimi.
\newblock {\em Applications of percolation theory}.
\newblock Taylor \& Francis, (1994).

\bibitem[Sla02]{Slad02}
Gordon Slade.
\newblock Scaling limits and super-{B}rownian motion.
\newblock {\em Notices Amer. Math. Soc.}, {\bf 49}(9):1056--1067, (2002).

\bibitem[SS04]{SidSzn04}
Vladas Sidoravicius and Alain-Sol Sznitman.
\newblock Quenched invariance principles for walks on clusters of percolation
  or among random conductances.
\newblock {\em Probab. Theory Related Fields}, 129(2):219--244, 2004.

\bibitem[Sta85]{Stau85}
Dietrich Stauffer.
\newblock {\em Introduction to percolation theory}.
\newblock Taylor \& Francis, Ltd., London, 1985.

\bibitem[SZ99]{SznZer1999}
Alain-Sol Sznitman and Martin Zerner.
\newblock A law of large numbers for random walks in random environment.
\newblock {\em Ann. Probab.}, pages 1851--1869, 1999.

\end{thebibliography}

\end{document}